\documentclass[a4paper,10pt,reqno]{amsart}

\usepackage{xcolor,calc, blindtext}
\usepackage{amsmath,amssymb,mathtools,amsthm,mathrsfs,bm}
\numberwithin{equation}{section}

\usepackage{YDLthesissymbols}

\usepackage[colorlinks]{hyperref}

\usepackage{setspace}
\usepackage[margin=1.2in]{geometry} 
\usepackage{graphicx}     
 
\usepackage{fancyhdr}  
\usepackage{appendix}
\usepackage{subfiles}
\usepackage{bbm}
\usepackage{tikz-cd,hf-tikz}
\usepackage[utf8]{inputenc}
\usepackage[T1]{fontenc}
\usepackage[normalem]{ulem}
\usepackage{indentfirst} 
\usepackage{enumerate} 
\usepackage{enumitem} 
\usepackage{lscape}
\usepackage{etoolbox}
\usepackage{stmaryrd}
\usepackage{tabularx}
\usepackage{import}
\usepackage{multirow}

\usepackage[french,english]{babel}

\allowdisplaybreaks

\begin{document}
\title{Approche non-invariante de la correspondance de Jacquet-Langlands : analyse spectrale}
\author{Yan-Der LU}
\date{\today}

{\selectlanguage
{english}
\begin{abstract}
This is the second article in a two-part series presenting a new proof comparing the non-invariant trace formula for a general linear group with that of one of its inner forms. In this article, we focus on the spectral side of the trace formula. We complete the proof of the global Jacquet-Langlands correspondence using the non-invariant trace formula and examine its arithmetic implications. Furthermore, we define the notion of non-invariant spectral transfer of a test function and show that it coincides with the non-invariant geometric transfer introduced in our first article. This provides a positive answer to a conjecture of Arthur and extends a well-known theorem of Kazhdan within our framework.
\end{abstract}}

\selectlanguage{french}
\maketitle

\tableofcontents

\section{Introduction}
\subsection{Préface}
Nous avons entamé l'étude de la correspondance de Jacquet-Langlands globale à travers la formule des traces non-invariante dans un travail précédent \cite{YDL23b}, où l'analyse géométrique a été menée à terme. Le présent article a pour objectif de mener l'analyse spectrale et de compléter la démonstration de la correspondance de Jacquet-Langlands globale. Nous étudierons également les identités de transfert, tant au niveau local que global.

\subsection{Résultats principaux}
Dans la suite, $F$ sera un corps de nombres, et $G$ un groupe réductif sur $F$. On note $\g$ l'algèbre de Lie de $G$. Soit $v$ une place de $F$, on note $F_v$ le complété local de $F$ en $v$. On note $G_v$ le changement de base $G\times_FF_v$. On note $F_\infty$ le produit des $F_v$ avec $v$ parcourt l'ensemble des places archimédiennes de $F$. Soient $\A_F$ l'anneau des adèles de $F$, et $\A_{F,\fin}$ l'anneau des adèles finis de $F$. On a $\A_F=\A_{F,\fin}\times F_\infty$.

On définit dans un premier temps un espace des fonctions tests pour la formule des traces non-invariante. Soit $\mathcal{U}(\mathfrak{g}(F_\infty)_\C)$ l'algèbre enveloppante de $\mathfrak{g}(F_\infty)_\C\eqdef \mathfrak{g}(F_\infty)\otimes_{\R}\C$. 

\begin{definition}[M\"{u}ller-Finis-Lapid, définition \ref{def:espaceMuLaFi}]~{}
\begin{enumerate}
    \item Pour $J$ un sous-groupe ouvert compact de $G(\A_{F,\fin})$, soit $\mathcal{C}(G(\A_F),J)$ l'espace des fonctions lisses invariantes à droite par $J$ sur $G(\A_F)$, telles que les fonctions elle-mêmes ainsi que toutes leurs dérivées appartiennent à $L^1(G(\A_F))$. C'est un espace de Fréchet muni des semi-normes $(\|-\|_X)_{X\in\mathcal{U}(\mathfrak{g}(F_\infty)_\C)}$ avec
    \[\|f\|_X\eqdef\| f\ast X\|_{L^1(G(\A_F))},\,\,\,\,f\in \mathcal{C}(G(\A_F),J).\]
    \item On note $\mathcal{C}(G(\A_F))$ la limite inductive des $\mathcal{C}(G(\A_F),J)$ avec $J$ parcourt l'ensemble des sous-groupes ouverts compacts de $G(\A_{F,\fin})$, il est muni de la topologie limite inductive.
\end{enumerate}

\end{definition}

Par un sous-groupe de Levi d'un groupe réductif on entend le facteur de Levi d'un sous-groupe parabolique. Soient $ M \subseteq H$ des sous-groupes de $G$. On note $\L^H(M)$ l’ensemble des sous-groupes de Levi de $G$ inclus dans $H$ et contenant $M$ ; on note $\P^H(M)$ l’ensemble des sous-groupes paraboliques de $G$ inclus dans $H$, dont $M$ est un facteur de Levi ; on note $\F^H(M)$ l’ensemble des sous-groupes paraboliques de $G$ inclus dans $H$ et contenant $M$. On fixe $M_0$ un sous-groupe de Levi minimal de $G$. Lorque $M = M_0$, on omettra souvent $(M_0)$ dans les trois notations précédentes, i.e. $\L^H=\L^H(M_0)$, $\P^H=\P^H(M_0)$ et $\F^H=\F^H(M_0)$. Soit $M\in \L^G$. Pour un élément $P\in \F^G(M)$, on écrit $M_P$ le facteur de Levi de $P$ qui contient $M$. 

On note $X^\ast(G)$ le groupe des caractères de $G$ définis sur $F$ et $a_G$ (resp. $a_G^\ast$) l’espace vectoriel réel $\Hom(X^\ast(G),\R)$ (resp. $X^\ast(G)\otimes_\Z\R$). Les espaces $a_G$ et $a_G^\ast$ sont duaux l'un de l'autre. pour tout $M,L\in \L^G$ avec $M\subseteq L$, la restriction définit une inclusion naturelle $a_L^\ast\subseteq a_M^\ast$. On note $a_M^L$ le noyau de $a_L^\ast$ dans $a_M$, puis $(a_M^L)^\ast\subseteq a_M^\ast$ son dual. On a $a_M=a_L\oplus a_M^L$ et $a_M^\ast= a_L^\ast \oplus (a_M^L)^\ast$.

Pour $M\in \L^G$, on définit le groupe de Weyl relatif de $(G,M)$ par
\[W^G(M)=W_M^G\eqdef\text{Norm}_{G(F)}(M)/M(F),\]
avec $\text{Norm}_{G(F)}(M)$ le normalisateur de $M$ dans $G(F)$. Pour $M=M_0$ on abrège $W_{M_0}^G$ en $W_0^G$. Le groupe $W^G(M)$ agit par conjugaison sur $\L^G(M)$, $\P^G(M)$, et $\F^G(M)$. On écrit 
\[W^G(M)_\reg=\{s\in W^G(M)\mid \text{l'ensemble des points fixes de $a_M$ par $s$ est $a_G$}\}.\]


On note $J_\spec^G(-)$ le côté spectral de la formule des traces non-invariante (cf. théorème \ref{YDLspecthm:non-invtraceformula}). Puis on définit les caractères pondérés globaux $J_M^G(\pi,-)$ ainsi que $\J_M^G(\pi,-)=\int_{ i(a_M^G)^\ast}J_M^G(\pi_\lambda
,-)\,d\lambda$ (discussion à partir de l'équation \eqref{eqdef:weightedchara}). Notons $G(\A_F)^1$ le noyau de l'application de Harish-Chandra $H_G:G(\A_F)\to a_G$. On établit le développement fin du côté spectral pour l'espace des fonctions tests $\CMFL$. Posons $\Pi_{\disc}(G(\A_F)^1)$ l’ensemble des classes d’isomorphismes des sous-représentations irréductibles de la représentation régulière à droite de $G(\A_F)^1$ sur $ L^2(G(F)\backslash G(\A_F)^1)
$, puis $\Pi_\disc(G)$ l'ensemble constitué des sous-quotients irréductibles des représentations induites
\[\I_P^G(\sigma),\,\,\,\,M\in\L^G,P\in\P^G(M),\sigma\in\Pi_\disc(M(\A_F)^1)\]
de $G(\A_F)^1$, où $\sigma$ vérifie la condition suivante : il existe $s\in W^G(M)_\reg$ tel que $s\cdot\sigma=\sigma$. On voit $\Pi_\disc(G)$ comme un espace topologique discret dénombrable. On regarde ensuite l'ensemble
\[\coprod_{M\in\L^G}\coprod_{P\in\P^G(M)}\{\I_{P}^G(\pi_\lambda)\mid \pi\in\Pi_\disc(M),\lambda\in i(a_M^G)^\ast\}.\]
On définit une relation d'équivalence $\sim$ sur cet ensemble par $\I_{P}^G(\pi_\lambda)\sim \I_{P'}^G(\pi_{\lambda'}')$ si et seulement si $M_P=M_{P'}$, $\pi=\pi'$ et $\lambda=\lambda'$. Introduisons le quotient
\begin{equation}
\Pi(G)=\coprod_{M\in\L^G}\coprod_{P\in\P^G(M)}\{\I_{P}^G(\pi_\lambda)\mid \pi\in\Pi_\disc(M),\lambda\in i(a_M^G)^\ast\}/\sim.    
\end{equation}
On tient à préciser que les éléments de $\Pi(G)$ sont les classes des symboles formels $\I_{P}^G(\pi_\lambda)$. 

\begin{theorem}[{{développement fin, théorème \ref{thm:devlopfinspec}}}]~{}
\begin{enumerate}
    \item Pour tout $M\in\L^G$ et $\pi\in\Pi(M)$, les expressions $J_M^G(\pi,-)$ et $\J_M^G(\pi,-)$ sont des distributions (i.e. fonctionnelles continues) sur $\CMFL$.
    \item Il existe, pour tout $M\in \L^G$ et $\pi\in \Pi(M)$, un coefficient $a^M(\pi)\in\C$ qui ne dépend que de $M$ et $\pi$ et non de $G$, tel que pour tout $f\in \CMFL$, on ait
    \[J_\spec^G(f)=\sum_{M\in\L^G}|W_0^M||W_0^G|^{-1}\int_{\Pi(M)}a^M(\pi)\J_M^G(\pi,f)\,d\pi.\]
    Toutes les expressions intervenant dans la définition de l'intégrale de droite sont absolument convergentes, et cette intégrale elle-même l'est.
\end{enumerate}
\end{theorem}

Nous procédons ensuite à l'examen de la correspondance de Jacquet-Langlands, en débutant par l'étude locale. Soit donc $v$ une place de $F$. Soit $\Pi_{2}(G(F_v))$ (resp. $\Pi_{\temp}(G(F_v))$ ; $\Pi_\unit(G(F_v))$) l'ensemble des classes d’isomorphismes des représentations irréductibles de carré intégrable (resp. tempérées ; unitaires). On a $\Pi_{2}(G(F_v))\subseteq \Pi_{\temp}(G(F_v))\subseteq \Pi_{\unit}(G(F_v))$. On note $G_{\rss}$ l'ensemble des éléments semi-simples réguliers de $G$. Si $\pi\in \Pi_\unit(G(F_v))$, on écrit $\Theta_\pi$ pour le caractère de Harish-Chandra de $\pi$. C'est une fonction sur $G_\rss(F_v)$ invariante par $G(F_v)$-conjugaison. 

Soit $D_v$ une algèbre simple centrale sur $F_v$ de degré $d_v$. Pour tout entier strictement positif $m$, on a un torseur intérieur $\eta: \GL_{m,D_v\otimes_{F_v}\overline{F_v}}\to \GL_{md_v,\overline{F_v}}$, avec $\overline{F_v}$ une clôture algébrique de $F_v$. Par abus de notation on écrira aussi $\eta: \GL_{m,D_v}\to \GL_{md_v,F_v}$. La notion de transfert pour les classes de conjugaison est définie comme dans le texte précédent (\cite[définition 2.4]{YDL23b}). Nous en rappelons ici la définition. Soient $\o_v$ (resp. $\o_v'$) une classe de $\GL_{m}(D_v)$-conjugaison (resp. $\GL_{md_v}(F_v)$-conjugaison) dans $\gl_{m}(D_v)$ (resp. $\gl_{md_v}^\ast(F_v)$). On dit que $\mathfrak{o}_v'$ se transfère en $\mathfrak{o}_v$ si $\eta((\Ad \GL_{m}(D_v\otimes_{F_v}\overline{F_v}))\mathfrak{o}_v)=(\Ad \GL_{md_v}(\overline{F_v}))\mathfrak{o}_v'$. On écrit $\o_v'=\o_v^\ast$ et $\o_v\arr \o_v^\ast$, et on dit que $\o_v^\ast$ se transfère à $\GL_{md_v,F_v}$.

On écrit $\mathcal{R}(\GL_m(D_v))$ pour le groupe de Grothendieck des représentations admissibles de longueur finie de $\GL_m(D_v)$. Posons $\mathcal{R}(D_v)=\bigoplus_{k\in \N}\mathcal{R}(\GL_m(D_v))$ avec pour convention $\GL_0(D_v)=\mathcal{R}(\GL_0(D_v))=\Z$. C'est un anneau commutatif avec l'induction parabolique comme loi multiplicative que l'on notera $\times$. 

Divers travaux en théorie des représentations ont permis d'établir la correspondance locale de Jacquet-Langlands :


\begin{theorem}[{{correspondance locale de Jacquet-Langlands, théorème \ref{thm:appLJlocal}, corollaire \ref{coro:localJLtemp}}}]
Il existe un morphisme d'anneaux $\bLJ_v:\mathcal{R}(F_v)\to\mathcal{R}(D_v)$ vérifiant 
\begin{enumerate}
    \item $\bLJ_v=0$ sur $\bigoplus_{k\in \N,d_v\nmid k}\mathcal{R}(\GL_k(F_v))$, et $\bLJ(\mathcal{R}(\GL_{md_v}(F_v)))\subseteq\mathcal{R}(\GL_{m}(D_v))$ pour tout $m\in \mathbb{N}$ ;
    \item $\Theta_{\pi_v'}(\o_v^\ast)=e(\GL_{m,D_v})\Theta_{\bLJ_v(\pi_v')}(\o_v)$, pour tout $m\in \N$, $\pi_v'\in \mathcal{R}(\GL_{md_v}(F_v))$, $\o_v\subseteq \GL_{m,\rss}(D_v)$, et $\o_v^\ast\subseteq \GL_{md_v,\rss}^\ast(F_v)$ avec $\o_v\arr \o_v^\ast$. Ici $e(\GL_{m,D_v})\eqdef (-1)^{m(d_v-1)}$ est le signe de Kottwitz de $\GL_{m,D_v}$ ; 
    \item $\bLJ_v$ est surjectif ;
    \item $\bLJ_v$ se restreint à une bijection $\Pi_{2}(\GL_{md_v}(F_v))\rightarrow \Pi_{2}(\GL_{m}(D_v))$ pour tout $m\in \N$;
    \item $\bLJ_v(\Pi_{\unit}(\GL_{md_v} (F_v)))\subseteq \{0\}\sqcup\pm \Pi_{\unit}(\GL_m(D_v))$ pour tout $m\in \N\setminus\{0\}$ ;
    \item $\bLJ_v(\pi_{v,\lambda^\ast}')=\bLJ_v(\pi_v')_\lambda$ pout tous $m\in\N$, $\pi_v'\in \mathcal{R}(\GL_{md_v}(F_v))$ et $\lambda^\ast\in a_{\GL_{md_v,F_v}}^\ast$, ici $\lambda\in a_{\GL_{m,D_v}}^\ast$ est l'image réciproque de $\lambda^\ast$ par l'isomorphisme $a_{\GL_{m,D_v}}^\ast\xrightarrow{\sim}a_{\GL_{md_v,F_v}}^\ast$ induit par le torseur intérieur associé $\eta: \GL_{m,D_v}\to \GL_{md_v,F_v}$.
\end{enumerate}
C'est l'unique morphisme d'anneaux vérifiant les points 2, 4 et 6. Soit $m\in \mathbb{N}\setminus\{0\}$. Alors, par restriction, nous avons $\bLJ_v:\Pi_{\temp}(\GL_{md_v}(F_v))\rightarrow \{0\}\sqcup \Pi_{\temp}(\GL_{m}(D_v))$, et son image contient $\Pi_{\temp}(\GL_{m}(D_v))$. De plus, toute fibre sur $\Pi_{\temp}(\GL_{m}(D_v))$ de cette restriction est un singleton. 
\end{theorem}

On pose
\[\Pi_{\unit}^{d_v}(\GL_{md_v}(F_v))\eqdef \{\pi_v'\in \Pi_{\unit}(\GL_{md_v}(F_v))\mid \bLJ_v(\pi_v')\not=0\}.\]
On écrit $|\bLJ_v|:\Pi_{\unit}^{d_v}(\GL_{md_v}(F_v))\rightarrow \Pi_{\unit}(\GL_m(D_v))$ l'application qui à $\pi_v'$ associe la vraie représentation parmi $\{\bLJ_v(\pi_v'),-\bLJ_v(\pi_v')\}$.

On se donne un caractère additif non-trivial $\psi_v :F_v\to \C$. Pour $\pi_v$ une représentation irréductible admissible de $\GL_m(D_v)$ on note $L_v(\pi_v,s)$ sa fonction $L$ locale, $\epsilon_v(\pi_v,\psi_v,s)$ son facteur epsilon local, et $\gamma_v(\pi_v,\psi_v,s)$ son facteur gamma local. Par le biais de la théorie de Godement-Jacquet (rappelée en sous-section \ref{subsec:reviewGJtheory}), on établit : 
\begin{proposition}[{{proposition \ref{prop:Lepsilongammalocal}}}]~{}
\begin{enumerate}
    \item Soit $\pi_v'\in \Pi_{\unit}^{d_v}(\GL_{md_v}(F_v))$, alors
    \begin{align*}
    L_v(\pi_v',s)\mid L_v(|\bLJ_v|(\pi_v'),s),\,\,\epsilon_v(\pi_v',\psi_v,s)&=\epsilon_v(|\bLJ_v|(\pi_v'),\psi_v,s) \\
    \,\,\text{et}\,\,\gamma_v(\pi_v',\psi_v,s)&=\gamma_v(|\bLJ_v|(\pi_v'),\psi_v,s).\    
    \end{align*}
    Plus précisément le quotient de fonctions $L$ est un polynôme de $q_v^{-s}$ de terme constant 1, avec $q_v$ le cardinal du corps résiduel de $F_v$.
    \item Soit $\pi_v^\ast\in \Pi_{\unit}^{d_v}(\GL_{md_v}(F_v))\cap \Pi_\temp(\GL_{md_v}(F_v))$, alors
    \[L_v(\pi_v^\ast,s)=L_v(\bLJ_v(\pi_v^\ast),s)\]
\end{enumerate}
\end{proposition}

Abordons à présent la correspondance de Jacquet-Langlands globale. Soit donc  $\eta:\GL_{m,D}=G\rightarrow G^\ast=\GL_{n,F}$ un torseur intérieur, avec $D$ une algèbre simple centrale sur $F$ de degré $d=n/m$. Nous définissons la notion de transfert pour les sous-groupes de Levi, les sous-groupes paraboliques, ainsi que pour les classes de conjugaison comme dans le texte précédent (\cite[définitions 2.4, 2.7]{YDL23b}). On suit aussi les conventions de \cite[sous-section 2.4]{YDL23b} : on note, en ajoutant une étoile $^\ast$, en exposant les objets relativement à $G^\ast$ qui se transfèrent à $G$ ; pour les objets généraux relativement à $G^\ast$ on les note en ajoutant un signe prime $'$ en exposant. Par exemple en écrivant $M^\ast\in \L^{G^\ast}$ il est sous-entendu qu'il existe $M\in\L^G$ tel que $M^\ast$ se transfère en $M$. Alors que l'écriture $M'\in \L^{G^\ast}$ n'impose pas de condition particulière sur $M'$. Le contexte devrait permettre d’enlever toute ambiguïté. Soit $M'\in\L^{G^\ast}$. Quand $M'$ se transfère à $G$ on écrit $M'=M^\ast$ et $M'\arr M^\ast$, avec $M\in\L^G$ le transfert de $M^\ast$. Les mêmes notations valent pour les sous-groupes paraboliques.

On note $\S(\g(\A_F))$ (resp. $\S(\g(F_v))$) l'espace de Schwartz-Bruhat de $\g(\A_F)$ (resp. $\g(F_v)$). On pose $\nu_\g:\g\to \gl_{1,F}$ la norme réduite sur $\g$. On note $|-|_{\A_F}$ (resp. $|-|_{v}$) la valeur absolue usuelle sur $\A_F$ (resp. $F_v$). Pour $f\in \S(\g(\A_F))$ (resp. $\S(\g(F_v))$) et $s\in\{-\frac{n-1}{2}\}\cup\{s'\in\C\mid \Re(s')>-\frac{n-1}{2}\}$, posons $f_s$ la fonction $X\in \g(\A_{F})\mapsto f(X)|\nu_{\g}(X)|_{\A_{F}}^{s+\frac{n-1}{2}}\in \C$ (resp. $X\in \g(F_v)\mapsto f(X)|\nu_{\g}(X)|_{v}^{s+\frac{n-1}{2}}\in \C$). Définissons $\cZ(\g(\A_{F}),s)\eqdef \{f_s \mid f\in \S(\g(\A_{F}))\}$ et $\cZ(\g(F_v),s)\eqdef \{f_s \mid f\in \S(\g(F_v))\}$) (définitions \ref{def:espacedeGJ}, \ref{def:GJuniontopo}). Ces espaces nous permettent de relier l'espace de Schwartz-Harich-Chandra et l'espace $\CMFL$ :

\begin{proposition}[{{propositions \ref{prop:GJinSHC}, \ref{prop:espacedeGJsubseteqMuFiLa}}}]~{}
\begin{enumerate}
    \item Soit $v$ une place de $F$. On écrit $G_v=\GL_{m_v,D_v}$ avec $D_v$ une algèbre simple centrale sur $F_v$ de degré $d_v$, et $m_vd_v=n$.  On note $\mathcal{HC}(G(F_v))$  l'espace de Schwartz-Harich-Chandra de $G(F_v)$ (cf. définition \ref{YDLspecdef:SCHspace}). On a une inclusion continue, où l'espace de gauche est muni de la topologie donnée dans la définition \ref{def:GJuniontopo},
    \[
    \bigcup_{s\in\C : \Re(s)>\frac{-d_v+1}{2}}\cZ(\g(F_v),s)|_{G(F_v)}\subseteq \mathcal{HC}(G(F_v)).    
    \]
    \item On a une inclusion continue, où l'espace de gauche est muni de la topologie donnée dans la définition \ref{def:GJuniontopo},
    \[\bigcup_{s\in \C:\Re(s)>\frac{n+1}{2}}\cZ(\g(\A_F),s)|_{G(\A_F)}\subseteq \CMFL.\]
\end{enumerate}  
\end{proposition}

Nous voulons avoir recours au développement fin du côté spectral de la formule des traces pour la déduction de la correpondance de Jacquet-Langlands globale. Il n'est pas possible d'utiliser une fonction de classe Schwartz-Bruhat comme fonction test pour des raisons de convergence. Il nous faut faire appel à l'espace $\mathcal{C}(G(\A_F))$. \`{A} la suite de l'analyse des côtés géométriques exposée dans le texte précédent \cite{YDL23b}, nous concluons que $J_\spec^G(f_s)=J_\spec^{G^\ast}(f_s^\ast)$ pour tout $\S(\g(\A_F))\ni f\underset{\geom}{\arr} f^\ast \in \S(\g^\ast(\A_F))$ et $s\in\C$ avec $\Re(s)>\frac{n+1}{2}$ (définition \ref{def:transfergeomgeom'spec}, proposition \ref{prop:innerformspecsideequality}). Par la suite, en maniant un argument d'isolement des représentations du essentiellement à Lapid-Rogawski (lemme \ref{lem:isolement}), et en suivant les démarches de Badulescu (\cite{Badu08}), nous établissons la correspondance de Jacquet-Langlands pour le spectre discret. Enfin, en exploitant des résultats de la classification des représentations tels que l'irréductibilité de l'induite parabolique d'une représentation unitaire irréductible (théorème \ref{thm:irredunitparaind}) et une généralisation du théorème fort de multiplicité un (proposition \ref{prop:rigidificationL2spec}), nous généralisons cette correspondance à tout spectre $L^2$, mettant en lumière ses propriétés arithmétiques :

\begin{theorem}[{{théorème \ref{thm:globalJLL2}, corollaire \ref{coro:LandgammaunderglobalJL}}}]
Notons $\Pi_{L^2}(G(\A_F)^1)$ l’ensemble des classes d’isomorphismes des sous-quotients irréductibles de la représentation régulière à droite de $G(\A_F)^1$ sur $ L^2(G(F)\backslash G(\A_F)^1)$.
\begin{enumerate}
    \item Il existe une unique application $\bG:\Pi_{L^2}(G(\A_F)^1)\rightarrow \Pi_{L^2}(G^\ast(\A_F)^1)$ vérifiant que, pour tout $\pi\in \Pi_{L^2}(G(\A_F)^1)$, si $\pi^\ast=\bG(\pi)$ alors $|\bLJ_v|(\pi_v^\ast)=\pi_v$ pour toute place $v\in \V_F$. L'application $\bG$ est injective, et son image, notée $\Pi_{L^2}^D(G^\ast(\A_F)^1)$, est constituée des représentations $\otimes_{v\in\V_F}\pi_v'\in \Pi_{L^2}(G^\ast(\A_F)^1)$ telles que $\pi_v'\in \Pi_{\unit}^{d_v}(G^\ast(F_v))$ pour tout $v\in \V_F$.
    \item Pour $\pi\in \Pi_{L^2}(G(\A_F)^1)$ on note $L(\pi,s)$ sa fonction L globale, et $\epsilon(\pi,s)$ son facteur epsilon global. On a 
    \[L(
\bG(\pi),s)\mid L(\pi,s)\,\,\,\,\text{et}\,\,\,\,\epsilon(\bG(\pi),s)=\epsilon(\pi,s).\]
\end{enumerate}    
\end{theorem}

En supposant que la conjecture de Ramanujan-Petersson soit vérifiée pour les formes intérieures d'un groupe général linéraire, on pourra établir l'égalité $L(\bG(\pi),s)=L(\pi,s)$ pour toutes $\pi$ et $\bG(\pi)$ cuspidales (point 2 de la proposition \ref{prop:Lepsilongammalocal}).

Comparons maintenant tous les termes des développements fins des formules des traces non-invariantes. Pour $M\in\L^G$ on note $\bG^M:\Pi_{L^2}(M(\A_F)^1)\to \Pi_{L^2}(M^\ast(\A_F)^1)$ la correspondance de Jacquet-Langlands globale de $M(\A_F)^1$. Fixons un caractère additif non-trivial $\psi:F\backslash\A_F\to \C$. On définit un sous-espace $\S(\g(\A_F),\psi) \subseteq \S(\g(\A_F))$ en imposant des certaines conditions aux places archimédiennes (cf. numéro \ref{subsubsec:FTnon-inv}). 

\begin{theorem}[théorème {{\ref{thm:idenendoglobal}}}]~{}
\begin{enumerate}
    \item Soient $\L^G\ni M\arr M^\ast\in\L^{G^\ast}$, $\P^G(M)\ni P\arr P^\ast\in\P^{G^\ast}(M^\ast)$, $\L^G(M)\ni L\arr L^\ast\in \L^{G^\ast}(M^\ast)$, $\pi\in \Pi_\disc(M)$, et $\lambda\in i(a_M^L)^\ast$. On note $\pi^\ast=\bG^M(\pi)\in \Pi_\disc(M^\ast)$ et $\lambda^\ast=\eta(\lambda)\in i(a_{M^\ast}^{L^\ast})^\ast$. Alors nous avons
    \[a^L(\I_{P\cap L}^L(\pi_\lambda))=a^{L^\ast}(\I_{P^\ast\cap L^\ast}^{L^\ast}(\pi_{\lambda^\ast}^\ast)).\]
    \item Soient $
\S(\g(\A_F),\psi)\ni f\underset{\geom}{\arr} f^\ast\in \S(\g^\ast(\A_F),\psi)$, $L'\in\L^{G^\ast}$, $Q'\in \F^{G^\ast}(L')$, $\lambda'\in ia_{L'}^\ast$ et $\pi'\in\Pi_{L^2}(L'(\A_F)^1)$. Nous avons, pour tout $s\in\C$ avec $\Re(s)>\frac{n+1}{2}$, 
    \begin{equation*}
    J_{L'}^{Q'}(\pi_{\lambda'}',f_s^\ast) =\begin{cases}J_L^Q(\pi_\lambda,f_s) &, \text{ si $(L',Q',\lambda',\pi')=(L^\ast,Q^\ast,\lambda^\ast,\pi^\ast\eqdef\bG^L(\pi)$)} \\
    0 &, \text{ sinon}.
    \end{cases}    
    \end{equation*}
    Plus précisément, « sinon » signifie que $(L',Q',\lambda')=(L^\ast,Q^\ast,\lambda^\ast)$ et $\pi'\not\in \Pi_{L^2}^{D}(L^\ast(\A_F)^1)$, ou que $L'$ ne se transfère pas.
\end{enumerate}
\end{theorem}

Dans la dernière partie de l'article, notre objectif est de comprendre la notion de transfert des fonctions test. Soit $v$ une place non-archimédienne de $F$.

\begin{definition}[définition \ref{def:loctemptrans}]
Soit $L_v\in \L^{G_v}$. On dit qu'un élément $\pi_v'\in \Pi_\temp(L_v^\ast(F_v))$ se transfère à $L_v$ si $\pi_v\eqdef\bLJ_v(\pi_v')\in \Pi_\temp(L_v(F_v))$. On écrit $\pi_v'=\pi_v^\ast$ et $\pi_v\arr \pi_v^\ast$. 
\end{definition}

Soient $L_v\in \L^{G_v}$, $Q_v\in \F^{G_v}(L_v)$, et $\o_v$ une classe de $L_v(F_v)$-conjugaison dans $\mathfrak{l}_v(F_v)$. On pose $J_{L_v}^{Q_v}(\o_v,-)$ l'intégrale orbitale pondérée associée (cf. \cite[numéro 3.3.2]{YDL23b}). C'est une distribution sur $\S(\g(F_v))$.  

\begin{definition}[définition \ref{def:transfergeomgeom'spec}]~{}
\begin{enumerate}
    \item Pour tout $L_v\in \L^{G_v}$, $Q_v\in \F^{G_v}(L_v)$ et $\o_v$, une classe de $L_v(F_v)$-conjugaison dans $\mathfrak{l}_v(F_v)$ avec $(\Ad G_v)\o_v$ une classe de conjugaison semi-simple régulière, on appelle $(L_v,Q_v,\o_v)$ un triplet géométrique de $G_v$.
    \item Soient $f_v\in \S(\g(F_v))$ et $f_v^\ast\in \S(\g^\ast(F_v))$. Si
    \begin{align*}
    J_{L_v'}^{Q_v'}(\mathfrak{o}_v',f_v^\ast)=\begin{cases}
    J_{L_v}^{Q_v}(\mathfrak{o}_v,f_v) & \text{si $(L_v',Q_v',\mathfrak{o}_v')= (L_v^\ast,Q_v^\ast,\mathfrak{o}_v^\ast)$ ;} \\
    0 & \text{si $Q_v'$ ne se transfère pas.}
    \end{cases}
    \end{align*}
    pour tout triplet géométrique $(L_v',Q_v',\mathfrak{o}_v')$, alors on note $f_v\underset{\geom}{\arr} f_v^\ast$.
    
    \item Pour tout $L_v\in \L^{G_v}$, $Q_v\in \F^{G_v}(L_v)$ et $\pi_v\in \Pi_\temp(L_v(F_v))$, on appelle $(L_v,Q_v,\pi_v)$ un triplet spectral de $G_v$.
    \item Soient $f_v\in \S(\g(F_v))$ et $f_v^\ast\in \S(\g^\ast(F_v))$, si 
    \begin{align*}
    J_{L_v'}^{Q_v'}(\pi_v',f_{v,s}^\ast)=\begin{cases}
    e(G_v) J_{L_v}^{Q_v}(\pi_v,f_{v,s}) & \text{si $(L_v',Q_v',\pi_v')=(L_v^\ast,Q_v^\ast,\pi_v^\ast)$ ;} \\
    0 &\text{si $Q_v'$ ne se transfère pas.}
    \end{cases}
    \end{align*}
    pour tout $s\in \C$ tel que $\Re(s)>0$, et tout triplet spectral $(L_v',Q_v',\pi_v')$, alors on note $f_v\underset{\spec}{\arr} f_v^\ast$.
\end{enumerate}
\end{definition}

La notion $\underset{\geom}{\arr}$ est celle que nous avons employée tout au long du texte précédent et jusqu'à maintenant. On compare ces deux notions.

\begin{theorem}[{{théorèmes \ref{thm:transfertspecmeaning}, \ref{thm:transfergeom'specequiv}}}]Soient $f_v\in \S(\g(F_v))$ et $f_v^\ast\in\S(\g^\ast(F_v))$. Alors $f_v\underset{\geom}{\arr} f_v^\ast$ si et seulement si $f_v\underset{\spec}{\arr} f_v^\ast$.   
\end{theorem}

On remarque que l'implication directe de cette équivalence a déjà été conjecturée par Arthur dans \cite[conjecture 10.1]{Art96}, dans un contexte différent. En effet, il compare les caractères pondérés invariants des représentations tempérées d'un groupe sur un corps $p$-adique avec les caractères pondérés stables des représentations tempérées d'un groupe endoscopique elliptique. Notre résultat fournit ainsi une évidence de cette conjecture dans le cas particulier d'un groupe général linéaire et d'une de ses formes intérieures. Il est opportun de noter que cette conjecture d'Arthur se trouve uniquement dans son manuscrit original et non dans la version publiée du texte. L'auteur tient à remercier Arthur de l'avoir autorisé à citer le manuscrit original. 

\begin{corollary}[{{corollaire \ref{coro:kazhdan0noninvgrp}}}] ~{}
Soit $f_v\in C_c^\infty(G(F_v))$. On a l'équivalence entre
    \begin{itemize}
    \item $J_{L_v}^{Q_v}(\o_v,f_v)=0$ pour tout triplet géométrique $(L_v,Q_v,\o_v)$ avec $\o_v\subseteq L_v(F_v)$.
    \item $J_{L_v}^{Q_v}(\o_v,f_v)=0$ pour tout $L_v\in \L^{G_v}$, $Q_v\in \F^{G_v}(L_v)$ et $\o_v$ classe de $L_v(F_v)$-conjugaison dans $\mathfrak{l}_v(F_v)$.
    \item $J_{L_v}^{Q_v}(\pi_v,f_{v})=0$ pour tout triplet spectral $(L_v,Q_v,\pi_v)$.
    \item $J_{L_v}^{Q_v}(\pi_v,f_{v})=0$ pour tout $L_v\in \L^{G_v}$, $Q_v\in \F^{G_v}(L_v)$ et $\pi_v\in \Pi_{\unit}(L_v(F_v))$.
        \end{itemize}
\end{corollary}

Ceci constitue une généralisation non-invariante d'un résultat connu de Kazhdan (\cite[théorème 0, (b) $\Leftrightarrow$ (c) $\Leftrightarrow$ (d)]{Kaz86}).

\subsection{Plan de l'article}
Le contenu du présent texte est organisé de la manière suivante : dans la section \ref{sec:preliminaires}, nous introduisons les notions de base et démontrons le développement fin du côté spectral pour l'espace $\CMFL$ applicable à tout groupe réductif. Dans la section \ref{YDLspecsec:convention} nous fixons des conventions. La section \ref{sec:analyseespaceGJ} aborde la théorie de Godement-Jacquet et l'isomorphisme de Satake, tout en discutant le choix de l'espace des fonctions test locales et globales pour la formule des traces. Ensuite, dans la section \ref{sec:corrJL}, nous examinons la correspondance de Jacquet-Langlands ainsi que les implications arithmétiques qui en découlent. Pour conclure, dans la section \ref{sec:Idenspec}, nous établissons l'identité globale de transfert non-invariante attendue, suivie d’une comparaison des deux notions de transfert des fonctions tests évoquées ci-dessus.

\subsection*{Remerciements}
L'auteur tient particulièrement à remercier Pierre-Henri Chaudouard, son directeur de thèse. Sans ses connaissances, conseils, temps et patience, cet article n'aurait pas été possible. L'auteur tient également à remercier l'école doctorale 386 Sciences Mathématiques de Paris Centre ainsi que l'université Paris-Cité pour avoir financé ce projet de thèse.

\section{Préliminaires}\label{sec:preliminaires}

\subsection{Préliminaires}
On s'efforce de conserver les mêmes notations que dans le texte précédent \cite{YDL23b}. Comme on n'entreprend pas tout de suite la preuve de la correspondance de Jacquet-Langlands, on fixe dans cette section $G$ un groupe réductif défini sur $F$ un corps local ou global de caractéristique 0. Fixons $M_0$ un sous-groupe de Levi minimal de $G$. Si $F$ est un corps global, soit $\V_F$ (resp. $\V_{\fin}$ ; $\V_\infty$) l’ensemble des places (resp. places finies ; places archimédiennes) de $F$ ; si $F$ est un corps local, soit $\V_F$ le singleton dont l'unique élément est la place qui munie $F$. Pour $v\in \V_F$, on note $F_v$ le complété local de $F$ en $v$.

Si $F$ est un corps global, alors on note $\A_F$ son anneau des adèles, $\A_{F,\fin}$ son anneau des adèles finis, et $F_\infty\eqdef \prod_{v\in \V_\infty} F_v$. On a $\A_F= \A_{F,\fin}\times F_\infty$. 

Lorsque $F$ est un corps local, on note $|-|_F$ la valeur absolue usuelle sur $F$ ; lorsque $F$ est un corps global, on note $|-|_{\A_F}$ la valeur absolue usuelle sur $\A_F$, et on abrège $|-|_{F_v}$ en $|-|_v$ pour tout $v\in \V_F$. 

On note $X^\ast(G)$ le groupe des caractères de $G$ définis sur $F$ et $a_G$ (resp. $a_G^\ast$) l’espace vectoriel réel $\Hom(X^\ast(G),\R)$ (resp. $X^\ast(G)\otimes_\Z\R$). Soient $ M \subseteq H$ des sous-groupes de $G$. On note $\L^H(M)$ l’ensemble des sous-groupes de Levi de $G$ inclus dans $H$ et contenant $M$ ; on note $\P^H(M)$ l’ensemble des sous-groupes paraboliques de $G$ inclus dans $H$, dont $M$ est un facteur de Levi ; on note $\F^H(M)$ l’ensemble des sous-groupes paraboliques de $G$ inclus dans $H$ et contenant $M$. Lorque $M = M_0$, on omettra souvent $(M_0)$ dans les trois notations précédentes. 

Soit $M\in\L^G$. Soit $H$ un sous-groupe de $G$ stable par conjugaison par $A_M$. On note $\Sigma(\mathfrak{h}; A_M)$ l’ensemble des racines de $A_M$ sur $\mathfrak{h}$. Soient $P,Q\in\F^G(M)$ avec $P\subseteq Q$. Notons $M_P$ l'unique facteur de Levi de $P$ contenant $M$, et $N_P$ le radical unipotent de $P$. On note $\Delta_P^Q$ l’ensemble des racines simples dans $\Sigma(\m_Q\cap \n_P; A_{M_P})$. Les ensembles $\Sigma(\g; A_{M_P})$ et $\Sigma(\mathfrak{p}; A_{M_P})$ sont canoniquement inclus dans $a_{M_P}^\ast$. \`{A} chaque racine $\beta\in \Sigma(\g; A_{M})$ est associée une coracine $\beta^\vee\in a_{M}$. On définit alors dans $a_{M_P}$ le sous-ensemble $(\Delta_P^Q)^\vee$ des coracines des éléments de $\Delta_P^Q$. Soit $a_P^Q$ (resp. $(a_p^Q)^\ast$) le sous-espace vectoriel engendré par $(\Delta_P^Q)^\vee$ (resp. $\Delta_P^Q$).

\subsection{Normalisations des mesures de la section \ref{sec:preliminaires}}\label{subsec:normalisationmesureofspec1}

Pour $H$ un groupe séparé localement compact muni d'une mesure de Haar et $X\subseteq H$ une partie mesurable, on note $\vol(X;H)$ le volume de $X$ dans $H$ et $\vol(H)$ le volume total de $H$.

Soit dans ce numéro $F$ un corps de nombres. Posons $K=\prod_{v\in\V_F} K_v$ avec $K_v$ un sous-groupe compact maximal de $G(F_v)$ en bonne position par rapport à $M_{0,v}\eqdef M_0\times_FF_v$ (cf. \cite[page 10]{Art91}).

Dans la suite le symbole $v$ désigne un élément général de $\V_F$. Fixons $M_{v,0}$ un sous-groupe de Levi minimal de $G_v$, incluant dans $M_{0,v}\eqdef M_{0}$. On choisit un produit scalaire $W^G_{M_0}$-invariant (resp. $W^{G_v}_{M_{v,0}}$-invariant) sur $a_{M_0}$ (resp. $a_{M_v,0}$). On prend sur tout sous-espace de $a_{M_0}$ (resp. $a_{M_v,0}$) la mesure engendrée par ce produit scalaire. On aimerait préciser qu'a priori, la mesure sur $a_M$ n'a pas de lien avec celle sur $a_{M_v}$ pour $M\in \L^G$. 

Soit $M\in\L^G$. On note $A_{M,\Q}$ le sous-tore central $\Q$-déployé maximal dans $\Res_{F/\Q}M$ puis $A_{M,\infty}\eqdef A_{M,\Q}(\R)^\circ$ la composante neutre de $A_{M,\Q}(\R)\subseteq A_{M,\Q}(F_\infty)$, ici $\R$ s'injecte de façon diagonale dans $F_\infty$. On définit l'application de Harish-Chandra $H_M : M(\A_F)\rightarrow a_{M}$ par $e^{\langle H_M(m),\chi\rangle}=|\chi(m)|_{\A_F}$ pour tout $\chi\in X^\ast(M)$. Posons $M(\A_F)^1$ le noyau de l'application $H_M$. On a $M(\A_F)=M(\A_F)^1\times A_{M,\infty}$. Constatons que la restriction de $H_M$ sur $A_{M,\infty}$ donne un isomorphisme de groupes  $A_{M,\infty}\xrightarrow{\sim}a_M$, on se procure ainsi une mesure sur $A_{M,\infty}$ par cet isomorphisme.

Pour $P\in \P^G(M_0)$ ou $P\in \P^{G_v}(M_{v,0})$, on écrit $dH$ pour la mesure sur $a_P$ et $d\lambda$ pour la mesure duale sur $ia_P^\ast$, dans le sens où
\[\int_{ia_P^\ast}\int_{a_P}f(H)e^{-
\lambda(H)}\,dH\,d\lambda=f(0)\]
pour toute $f\in L^1(a_P)$. 

On munit le groupe $G(F_v)$ d'une mesure de Haar. On exige, comme il est loisible, que $\vol(K_v; G(F_v))=1$, pour presque tous les $v$. Le groupe $G(\A_F)^1$ est muni de l'unique mesure dont la multiplication avec la mesure sur $A_{G,\infty}$ vaut la mesure sur $G(\A_F)$.

La mesure sur un espace quotient est la mesure quotient. 

Pour tout $P \in \F^{G}(M_0)$, on munit $N_P(F_v)$ d'une mesure de Haar. On exige, comme il est loisible, que pour tout $s\in W_0^G$, la mesure sur $N_P(F_v)$ est obtenue par conjugaison par $w_s$ (un représentant fixé de $s$ dans $G(F)$) de la mesure sur $N_{(\Ad w_s)P}(F_v)$, et que $\vol(K_v\cap N(F_v); N(F_v))=1$ pour presque tous les $v$, enfin $\vol(N_P(F)\backslash N_P(\A_F))=1$, où $N_P(F)$ est muni de la mesure de comptage. Pour tout $M\in \L^G$ et $P,Q\in \P^G(M)$, on munit $N_Q(F_v)\cap N_P(F_v)$ d'une mesure de Haar. On exige, comme il est loisible, que $\vol(K_v\cap N_Q(F_v)\cap N_P(F_v); N_Q(F_v)\cap N_P(F_v))=1$ pour presque tous les $v$, et que $\vol((N_Q\cap N_P)(F)\backslash (N_Q\cap N_P)(\A_F))=1$ où $(N_Q\cap N_P)(F)$ est muni de la mesure de comptage.


\subsection{\texorpdfstring{$(G,M)$}{(G,M)}-familles}\label{subsec:GMfamilles}
On résume ici la théorie des $(G,M)$-familles d'Arthur. Pour tout espace vectoriel $V$ sur $\R$, on note $V_{\C}\eqdef V\otimes_{\R}\C$ sa complexification, et $iV\subseteq V_{\C}$ le sous-$\R$-espace évident où $i\in\C$ est l'unité imaginaire.

Soit $M\in\L^G$. Une $(G, M)$-famille est un ensemble d'applications lisses $\{c_P(\lambda)\mid P\in\P^G(M)\}$ de $\lambda\in ia_M^\ast$ vérifiant la propriété que si $P$ et $Q$ sont adjacents (i.e. $\Sigma(\mathfrak{p}; A_{M_0})\cap(-\Sigma(\mathfrak{q}; A_{M_0}))$ est un singleton) et si $\lambda$ appartient à l'hyperplan engendré par le mur en commun des chambres de $P$ et $Q$ dans $ia_M^\ast$ alors $c_P(\lambda) = c_{Q}(\lambda)$. \`{A} toute $(G,M)$-famille $\{c_P(\lambda)\mid P\in\P^G(M)\}$ et tout $Q\in \F^G(M)$, on associe naturellement une fonction $c_M^Q(\lambda)$, lisse en $\lambda\in ia_M^\ast$. On note $c_M^Q$ la valeur $c_M^Q(0)$. Si $Q=G$ on écrit simplement $c_M^G=c_M$.

\subsubsection{Formule de descente}
Moyennant certains choix, on peux montrer qu’il existe des applications :
\begin{equation}\label{eq:formuledescentGMsec}
\begin{split}
    d_M^G :\L^G(M)\times \L^G(M)&\rightarrow [0,+\infty[ \\
    s:\L^G(M)\times \L^G(M)&\rightarrow \F^G(M)\times \F^G(M)    
\end{split} 
\end{equation}
de sorte que, pour tout $(L_1,L_2)\in \L^G(M)\times \L^G(M)$, on ait
\begin{enumerate}
    \item si $s(L_1,L_2)=(Q_1,Q_2)$ alors $(Q_1,Q_2)\in \P^G(L_1)\times\P^G(L_2)$ ;
    \item si $s(L_1,L_2)=(Q_1,Q_2)$ alors $s(L_2,L_1)=(\overline{Q_2},\overline{Q_1})$ ;
    \item $d_M^G(L_1,L_2)\not =0$ si et seulement si l'une des flèches naturelles
    \[a_M^{L_1}\oplus a_M^{L_2}\longrightarrow a_M^G\]
    et
    \[a_{L_1}^G\oplus a_{L_2}^G\longrightarrow a_M^G\]
    est un isomorphisme, auquel cas les deux sont isomorphismes et $d_M^G(L_1,L_2)$ est le volume dans $a_M^G$ du parallélotope formé par les bases orthonormées de $a_M^{L_1}$ et de $a_M^{L_2}$ ;
    \item formule de descente : si $(c_P)_{P\in\P^G(M)}$ est une $(G,M)$-famille et $L\in \L^G(M)$, alors
    \[c_L=\sum_{L'\in \L^G(M)}d_M^G(L,L')c_M^{Q'}\]
    où l'on note $Q'$ la deuxième composante de $s(L,L')$.
\end{enumerate}

\subsubsection{Formule de scindage}
Soient $F$ un corps de nombres et $S$ un sous-ensemble non-vide de $\V_F$. Le symbole $v$ désigne ici un élément général de $S$. Soit $M\in \L^G$. Notons $M_S=\prod_{v\in S}M_v$ et le voit comme sous-groupe de Levi de $G_S=\prod_{v\in S}G_v$ défini sur $F_S$. On note de plus $\L^{G_S}(M_S)$ (resp. $\P^{G_S}(M_S)$ ; $\F^{G_S}(M_S)$) l'ensemble des produits $\prod_{v\in S}L_v$ (resp. $\prod_{v\in S}P_v$ ; $\prod_{v\in S}Q_v$) où $L_v$ est un sous-groupe de Levi défini sur $F_v$ contenant $M_v$ (resp. $P_v$ est un sous-groupe parabolique défini sur $F_v$ et de facteur de Levi $M_v$ ; $Q_v$ est un sous-groupe parabolique défini sur $F_v$ contenant $M_v$). Moyennant certains choix, on peux montrer qu’il existe des applications :
\begin{align*}
    d_M^G :\L^{G_S}(M_S)&\rightarrow [0,+\infty[ \\
    s:\L^{G_S}(M_S)&\rightarrow \F^{G_S}(M_S)
\end{align*}
de sorte que, pour tout $(L_v)_{v\in S}\in \L^{G_S}(M_S)$, on ait
\begin{enumerate}
    \item si $s\left((L_v)_{v\in S}\right)=\left((Q_v)_{v\in S}\right)$ alors $(Q_v)_{v\in S}\in \P^{G_S}(L_S)$ ;
    \item $d_M^G\left((L_v)_{v\in S}\right)\not =0$  si et seulement si l'une des flèches naturelles
    \[\bigoplus_{v\in S}a_{M_v}^{L_v}\longrightarrow a_M^G\]
    et
    \[\bigoplus_{v\in S}a_{L_v}^{G_v}\longrightarrow a_M^G\]
    est un isomorphisme, auquel cas les deux sont isomorphismes et $d_M^G\left((L_v)_{v\in S}\right)$ est le volume dans $a_M^G$ du parallélotope formé par les bases orthonormées des $a_{M_v}^{L_v}$ ;
    \item formule de scindage : supposons que nous disposons, pour tout $v\in S$, de $(c_{P_v})_{P_v\in\P^{G_v}(M_v)}$ une $(G_v,M_v)$-famille sur $F_v$. Définissons $c_P=\prod_{v\in S}c_{P_v}$. Alors $(c_P)_{P\in \P^G(M)}$ est une $(G,M)$-famille, et
    \begin{equation}\label{YDLgeomeq:formuledescindage}
    c_{M}=\sum_{(L_v)_{v\in S}\in \L^{G_S}(M_S)}d_M^G\left((L_v)_{v\in S}\right)\prod_{v\in S} c_{M_v}^{Q_v}.    
    \end{equation}
\end{enumerate}

\subsection{Notations}\label{YDLspecsubsec:prelinotations}

Toute représentation d'un groupe est supposé définie sur un espace de Hilbert le long de l'article. Soit $H$ un groupe topologique de type I. On note par $\Pi(H)$ (resp. $\Pi_\unit(H)$) l'ensemble des classes d'isomorphismes des représentations irréductibles (resp. représentations irréductibles unitaires) de $H$. L'élément $1=1_H\in \Pi_\unit(H)$ est la représentation triviale. Soit $H\curvearrowright V$ une représentation de $H$, on note $\langle-,-\rangle$ le produit sclaire sur $V$, $V_\pi$ la partie $\pi$-isotypique de $V$ pour $\pi\in \Pi(H)$, et $V=V_{\disc}\oplus V_{\cont}$ la décomposition (orthogonale) en la partie discrète et celle continue, enfin $m_{V,\disc}(\pi)$ est la multiplicité de $\pi$ dans $V_\disc$. Un entier $n$ est dit positif si $n\geq 0$, on note $\N$ l'ensemble des entiers positifs.

\subsubsection{Notations locales}
Supposons que $F$ est local.

Soit $(\pi,V)$ une représentation unitaire admissible de $G(F)$. On écrit $(\pi^\vee,V^\vee)$ la représentation contragrédiente de $(\pi,V)$. Puisque $(\pi,V)$ est unitaire, on peut identifier $(\pi^\vee,V^\vee)$ à $(\pi,V)$ via le produit scalaire $G(F)$-invariant duquel $V$ est muni. L'intégrale de Bochner 
\[\pi(f)\eqdef \int_{G(F)}f(g)\pi(g)\,dg\in \End_{\C}(V)\]
est bien définie pour toute fonction $f\in L^1(G(F))$, et donne un opérateur traçable lorsque $f$ est lisse à support compact. Le cas échéant, la distribution $f\mapsto\tr(\pi(f))$ est représentée par le caractère de Harish-Chandra $\Theta_\pi$. C'est une fonction sur $G_\rss(F)$ invariante par $G(F)$-conjugaison, avec $G_\rss$ l'ouvert dense de $G$ consiste en les éléments réguliers semi-simples. Deux représentations admissibles de longueur finie possèdent les mêmes facteurs de Jordan-H\"{o}lder si et seulement si leurs caractères de Harish-Chandra coïncident. Notons l'accouplement \[\langle \pi, f\rangle \eqdef\tr(\pi(f)).\]

On fixe $K$ un sous-groupe ouvert compact maximal de $G(F)$ si $F$ est non-archimédien (resp. $K$ un sous-groupe compact maximal de $G(F)$ si $F$ est archimédien), en bonne position par rapport à $M_0$. Pour tout $M\in\L^G(M_0)$ posons $K_M=M(F)\cap K$ qui est un sous-groupe ouvert compact maximal de $M(F)$ si $F$ est non-archimédien (resp. un sous-groupe compact maximal de $G(F)$ si $F$ est archimédien), en bonne position par rapport à $M_0$.

Soit $M\in \L^G(M_0)$. Soit $\pi$ une représentation admissible de $M(F)$. On écrit $\pi_{\lambda}=\pi\otimes \lambda$ pour la représentation tordue par $\lambda\in a_{M,\C}^\ast$, i.e.
\[\pi_{\lambda}(m)=\pi(m)e^{\lambda(H_M(m))},\,\,\,\,\forall m\in M(F).\]
Il est clair que $\langle\pi\otimes\lambda,f\rangle=\langle\pi,\lambda\cdot f\rangle$ pour tout $f$ fonction lisse à support compact et $\lambda\in a_{M,\C}^\ast$, ici à droite on identifie $\lambda$ à $1_\lambda$ la représentation triviale tordue par $\lambda$, qui est considérée comme une fonction à valeurs complexes. Soit $P\in\P^G(M)$, on note $\I_P^G(\pi_{\lambda})=\I_P^G(\pi,\lambda)$ la représentation induite de $M(F)$ à $G(F)$ par $\pi_{\lambda}$. On entend par l'induite  d'une représentation l'induite normalisée, c'est-à-dire telle que l'induite d'une représentation unitaire est une représentation unitaire. Soit $\lambda\in a_{G,\C}^\ast\subseteq a_{M,\C}^\ast$, alors $\I_P^G(\pi)\otimes\lambda=\I_P^G(\pi\otimes \lambda)$. Soit maintenant $Q\in \P^G(M)$ un autre sous-groupe parabolique, l'application qui à $\phi\in \I_P^G(\pi_\lambda)$, ce dernier vu comme espace de fonctions sur $G(F)$, associe 
\[x\in G(F)\mapsto \int_{N_P(F)\cap N_Q(F)\backslash N_Q(F)}\phi(nx)e^{\langle\lambda, H_P(nx)-H_Q(x)\rangle}\,dn,\]
définit un opérateur d'entrelacement $J_{Q|P}(\pi_\lambda)=J_{Q|P}(\pi,\lambda):\I_P^G(\pi_\lambda)\rightarrow \I_Q^G(\pi_\lambda)$ pourvu que $\Re(\lambda)$ soit assez régulier. Les opérateurs $J_{Q|P}(\pi_\lambda)$ vérifient les conditions données aux pages 26-27 de \cite{Art89I}, et admettent un prolongement méromorphe, vu comme fonction de $\lambda$, sur tout $a_{M,\C}^\ast$. Ils sont analytiques et unitaires sur $ia_M^\ast$.

Il existe des fonctions de $\lambda$ méromorphes à valeurs complexes $r_{Q|P}(\pi_\lambda)=r_{Q|P}(\pi,\lambda)$ 
telles que les opérateurs normalisés $N_{Q|P}(\pi_{\lambda})=N_{Q|P}(\pi,\lambda)\eqdef r_{Q|P}(\pi_\lambda)^{-1}J_{Q|P}(\pi_\lambda)$ vérifient les conditions données dans le théorème 2.1 de \cite{Art89I}, et sont analytiques et unitaires sur $ia_M^\ast$ (on précisera ces fonctions pour $G$ une forme intérieure d'un groupe général linéaire dans la sous-section \ref{subsec:normalisationr}). Définissons 
\begin{equation}\label{eq:opr}r_Q(P,\pi,\lambda,\Lambda)=r_{Q|P}(\pi_\lambda)
^{-1}r_{Q|P}(\pi_{\lambda+\Lambda}),\,\,\,\,\Lambda\in ia_M^\ast
\end{equation}
et
\begin{equation}\label{eq:opN}
\mathcal{N}_Q(P,\pi,\lambda,\Lambda)=N_{Q|P}(\pi_\lambda)
^{-1}N_{Q|P}(\pi_{\lambda+\Lambda}),\,\,\,\,\Lambda\in ia_M^\ast.    
\end{equation}
Alors $(r_Q(P,\pi,\lambda,\Lambda))_{Q\in \P^G(M)}$ et $(\mathcal{N}_Q(P,\pi,\lambda,\Lambda))_{Q\in \P^G(M)}$ considérés comme fonctions de $\Lambda$ sont des $(G,M)$-familles, ils donnent donc les poids $r_L(P,\pi,\lambda)$ et $\mathcal{N}_L(P,\pi,\lambda)$ pour tout $L\in \L^G(M)$. On écrit $r_L(P,\pi,\lambda)=r_L(\pi,\lambda)$ comme il est indépendant de $P$.

Supposons que $F=\R$ ou $\C$. Pour $\mathfrak{h}$ une algèbre de Lie réelle, on note $\mathcal{U}(\mathfrak{h})$ son algèbre enveloppante munie de la filtration usuelle, et $\mathcal{Z}(\mathfrak{h})$ le centre de $\mathcal{U}(\mathfrak{h})$.

\subsubsection{Notations globales}

Supposons que $F$ est global.

Dans ce texte on appelle une représentation automorphe au sens $L^2$ de $G(\A_F)^1$ tout sous-quotient irréductible de $R_G$, la représentation régulière à droite de $G(\A_F)^1$ sur $ L^2(G(F)\backslash G(\A_F)^1)$. Soit $M\in \L^G(M_0)$. Un élément de $\Pi_\unit(M(\A_F)^1)$ est parfois considéré comme un élément de $\Pi_\unit(M(\A_F))$ trivial sur $A_{M,\infty}$. On note par $\Pi_\disc(M(\A_F)^1)$ l'ensemble des classes d'isomorphismes des représentations irréductibles unitaires apparaîssant dans $L_\disc^2(M(F)\backslash M(\A_F)^1)$.  

Soient $\pi=\otimes_{v\in\mathcal{V}_F}\pi_v\in \Pi_\unit(M(\A_F))$ et $\lambda\in a_{M,\C}^\ast$. On pose $N_{Q|P}(\pi_{\lambda})=N_{Q|P}(\pi,\lambda)\eqdef \otimes_v N_{Q_v|P_v}(\pi_{v,\lambda})$, de même pour $r_{Q|P}(\pi_{\lambda})=r_{Q|P}(\pi,\lambda)$ et $J_{Q|P}(\pi_{\lambda})=J_{Q|P}(\pi,\lambda)$. Puis on définit les $(G,M)$-familles $(r_Q(P,\pi,\lambda,\Lambda))_{Q\in \P^G(M)}$ et $(\mathcal{N}_Q(P,\pi,\lambda,\Lambda))_{Q\in \P^G(M)}$ en suivant les mêmes équations \eqref{eq:opr} et \eqref{eq:opN}.

Pour $P\in \P^G(M)$ posons
$\overline{\mathcal{A}}^2(P)$ le complété hilbertien de 
\[\{\phi\in C^\infty(M(F)N_P(\A_F)\backslash G(\A_F)):\delta_P(\cdot)^{-1/2}\phi(\cdot x)\in L_\disc^2(M(F)\backslash M(\A_F)^1)\,\,\forall x\in G(\A_F)\}\]
par rapport au produit scalaire $\langle\phi_1,\phi_2\rangle=\int_{A_{M,\infty}M(F)N_P(\A_F)\backslash G(\A_F)}\phi_1(g)\overline{\phi_2(g)}\,dg$. Pour $\lambda\in a_{M,\C}^\ast$ on a la représentation $\I_{P(\A_F)}^{G(\A_F)}(R_{M,\disc,\lambda})$ du groupe $G(\A_F)$ sur $\overline{\mathcal{A}}^2(P)$. Notons $K_\infty=\prod_{v\in \V_\infty}K_v$. \'{E}crivons $\mathcal{A}^2(P)$ le sous-espace dense formé par les formes automorphes le long de $P$, i.e. les fonctions qui sont de plus $K_\infty$-finies et $\mathcal{Z}(\g({F_\infty})_\C)$-finies. Puis on écrit $\mathcal{A}_\pi^2(P)=\I_{P(\A_F)}^{G(\A_F)}(\pi)$ pour $\pi\in\Pi_{\disc}(M_P(\A_F)^1)$. 
Pour $H_1$ et $H_2$ deux sous-groups de Levi ou sous-groupes paraboliques semi-standards, on note $W^G(a_{H_1},a_{H_2})$ l'ensemble des isomorphisms de $a_{H_1}$ à $a_{H_2}$ obtenus à partir d'éléments de $W_0^G$ restreints sur $a_{H_1}$. Lorsque $H_1=H_2$ on note $W^G(a_{H_1},a_{H_2})=W^G(a_{H_1})$. Soient maintenant $P\in \P^G(M)$, $Q\in \P^G(M)$, $s\in W^G(a_P,a_Q)$, et $w_s\in G(F)$ un représentant fixé de $s$. Posons $M_{Q|P}(s,\lambda):\I_{P(\A_F)}^{G(\A_F)}(R_{M,\disc,\lambda})\rightarrow \I_{Q(\A_F)}^{G(\A_F)}(R_{M,\disc,\lambda})$ l'opérateur d'entrelacement défini par
\[(M_{Q|P}(s,\lambda)\phi)(x)=\int_{w_sN_P(\A_F)w_s^{-1}\cap N_Q(\A_F)\backslash N_Q(\A_F)}\phi(w_s^{-1}nx)e^{\langle\lambda, H_P(w_s^{-1}nx)\rangle-\langle s\lambda, H_Q(x)\rangle}\,dn\]
pourvu que $\Re(\lambda)$ soit assez régulier. Les opérateurs $M_{Q|P}(s,\lambda)$ vérifient les conditions données à la page 177 de \cite{FiLaMu11}, et admettent un prolongement méromorphe, vu comme fonction de $\lambda$, sur tout $a_{M,\C}^\ast$, analytique et unitaire sur $ia_M^\ast$. Lorsque $s=1$ on écrit $M_{Q|P}(s,\lambda)=M_{Q|P}(\lambda)$. L'opérateur $M_{Q|P}(\lambda)$ est bien sûr l'analogue global des $J_{Q|P}(\pi_{v,\lambda})$. Plus présisément, pour $\pi\in\Pi_{\disc}(M(\A_F)^1)$ on a un isomorphisme canonique de $G(\A_{F,\fin})\times (\g(F_\infty\otimes_{\R}\C),K_{\infty})$-modules 
\[j_P: \Hom(\pi, L^2(M(F)\backslash M(\A_F)^1))\otimes \I_{P(\A_F)}^{G(\A_F)}(\pi)\rightarrow \I_{P(\A_F)}^{G(\A_F)}(R_{M,\disc,\pi,\lambda}),\]
ce qui nous permet d'exprimer $M_{Q|P}(\lambda)\simeq \bigoplus_{\pi\in \Pi_{\disc}(M(\A_F)^1)} m_{\disc}(\pi)J_{Q|P}(\pi,\lambda)$. Définissons 
\[\mathcal{M}_Q(P,\lambda,\Lambda)=M_{Q|P}(\lambda)^{-1}M_{Q|P}(\lambda+\Lambda),\,\,\,\,\Lambda\in ia_M^\ast.\]
Alors $(\mathcal{M}_Q(P,\lambda,\Lambda))_{Q\in \P^G(M)}$ considéré comme fonction de $\Lambda$ est une $(G,M)$-famille, il donne donc $\mathcal{M}_L(P,\lambda)$ pour tout $L\in \L^G(M)$.

\subsection{Fomule des traces non-invariante}

Supposons que $F$ est un corps de nombres. 
Nous rappelons ici la définition du côté spectral de la formule des traces non-invariante $J_\spec^G$, et nous définissons un espace de fonctions tests. 

Soit $J$ un sous-groupe ouvert compact de $G(\A_{F,\fin})$, l'espace $G(\A_F)/J$ est une variété différentielle (une réunion dénombrable de copies de $G(F_\infty)$). Tout élément $X\in\mathcal{U}(\mathfrak{g}(F_\infty)_\C)$ de l'algèbre enveloppante définit $f\mapsto f\ast X$ un opérateur différentiel invariant à gauche sur $G(\A_F)/J$.

\begin{definition}[M\"{u}ller-Finis-Lapid]\label{def:espaceMuLaFi} Soit $\mathcal{C}(G(\A_F),J)$ l'espace des fonctions lisses invariantes à droite par $J$ sur $G(\A_F)$, telles que les fonctions elle-mêmes ainsi que toutes leurs dérivées appartiennent à $L^1(G(\A_F))$. C'est un espace de Fréchet muni des semi-normes $(\|-\|_X)_{X\in\mathcal{U}(\mathfrak{g}(F_\infty)_\C)}$ avec
\begin{equation}\label{YDLspeceq:defCMFLseminorm}
\|f\|_X\eqdef\| f\ast X\|_{L^1(G(\A_F))},\,\,\,\,f\in \mathcal{C}(G(\A_F),J).    
\end{equation}

Enfin on note $\mathcal{C}(G(\A_F))$ la limite inductive des $\mathcal{C}(G(\A_F),J)$ avec $J$ parcourt l'ensemble des sous-groupes ouverts compacts de $G(\A_{F,\fin})$. L'espace $\mathcal{C}(G(\A_F))$ est muni de la topologie limite inductive. 
\end{definition}

L'espace $\CMFL$ est une algèbre associative topologique sans unité pour le produit de convolution, et il 
contient $C_c^\infty(G(\A_F))$. 

On décrit ensuite le développement spectral de la formule des traces. Soit $M\in \L^G(M_0)$, pour tout $s\in W^G(M)$ posons $L_s$ le plus petit sous-groupe de Levi de $G$ dans $\L^G(M)$ contenant $w_s$, il est également caractérisé par la condition $a_{L_s} = \{H \in a_M\mid  sH = H\}$.

L'application de Harish-Chandra $H_G:A_{G,\infty}\rightarrow a_G$ étant un isomorphisme, elle nous permet d'interpréter $a_G$ comme l'algèbre de Lie de $A_{G,\infty}$. Pour toute fonction $f\in\CMFL$, on note $f^1$ la fonction sur $G(\A_F)^1$ définie par
\[g\in G(\A_F)^1\mapsto\int_{A_{G,\infty}}f(ag)\,da.\]
On peut également définir l'espace $\mathcal{C}(G(\A_F)^1)$, analogue à $\mathcal{C}(G(\A_F))$, en remplaçant, dans la définition de $\mathcal{C}(G(\A_F))$, le groupe $G(\A_F)$ partout par $G(\A_F)^1$ et $\mathfrak{g}(F_\infty)_\C$ par la complexification de l'algèbre de Lie de $G(F_\infty)\cap G(\A_F)^1$ (cf. \cite[page 426]{FiLa16}). Toute fonction de $\mathcal{C}(G(\A_F)^1)$ est de la forme $ f^1$ pour une fonction $f\in \mathcal{C}(G(\A_F))$. Pour le voir, on peut fixer $\phi$ une fonction positive lisse à support compact sur $A_{G,\infty}$ d'intégrale 1, alors $(h\phi)^1=h$ pour tout $h\in \mathcal{C}(G(\A_F)^1)$.

\begin{theorem}[M\"{u}ller-Finis-Lapid]\label{YDLspecthm:non-invtraceformula}~{}
\begin{enumerate}
    \item Pour toute fonction test $f\in\CMFL$, l'expression 
    \begin{align*}
    \sum_{[P]}\frac{1}{|W^G(M_P)|}\sum_{s\in W^G(M_P)}|&\det(s-1)_{a_{M_P}^{L_s}}|^{-1}\\
    &\cdot\int_{i(a_{L_s}^G)^\ast}\tr\left(\mathcal{M}_{L_s}(P,\lambda)M_{P|P}(s,0)\I_P^G(R_{M,\disc,\lambda},f^1)\right)\,d\lambda
    \end{align*}
    où la somme sur $[P]$ est portée sur l'ensemble des représentants dans $\F^G(M_0)$ des classes de sous-groupes paraboliques associés est absolument convergente (pour la norme trace), et donne une distribution $J_{\text{spec}}^G(-)$ sur $\CMFL$ (\cite[corollaire 1]{FiLaMu11}). C'est le côté spectral de la formule des traces non-invariante. 
    \item On note $J_\geom^G$ le côté géométrique de la formule des traces non-invariante (\cite[corollaire 7.2]{FiLa16}).
    Pour toute fonction test $h\in\mathcal{C}(G(\A_F)^1)$, l'expression $J_{\geom}^G(h)$ est absolument convergente et donne une distribution sur $\mathcal{C}(G(\A_F)^1)$. 
    \item Pour toute fonction test $f\in\CMFL$ telle que $f^1\in \mathcal{C}(G(\A_F)^1)$ il y a l'égalité $J_{\geom}^G(f^1)=J_{\spec}^G(f)$.
\end{enumerate}
\end{theorem}

\subsection{Développement fin du côté spectral}

Pour tout $L\in \L^G(M)$ on écrit $W^L(M)_\reg=\{s\in W^L(M):L_s=L\}$. Ainsi $s\in W^{L_s}(M)_\reg$ pour tout $s\in W^G(M)$.
Considérons l'expression 
\begin{align}
J_{\disc}^G(f)\eqdef\sum_{M\in\L^G}|W_0^M||W_0^G|^{-1}  \sum_{s\in W^G(M)_\reg}|\det(s-1)_{a_M^G}|^{-1}\tr\left(M_{P|P}(s,0)\I_P^G(R_{M,\disc,0},f^1)\right), 
\end{align}
avec $P$ un élément quelconque de $\P^G(M)$. Cela est la partie discrète de la formule des traces, elle contient la trace de $f^1$ sur le spectre automorphe discret $L_\disc^2(G(F)\backslash G(\A_F)^1)$, qui est le terme correspondant à $M=G$ et $s=1$ dans la somme. Au vu de l'admissibilité de $\I_P^G(R_{M,\disc,0})$ et le lemme de Schur, on a 
\[J_{\disc}^G(f)=\sum_{\pi\in\Pi_\unit(G(\A_F)^1)}a_\disc^G(\pi)\tr(\pi(f^1))\]
avec $a_\disc^G(\pi)$ un nombre complexe uniquement déterminés par $G$ et $\pi$. On aimerait travailler sur un sous-ensemble de $\Pi_\unit(G(\A_F)^1)$ qui contient le support de $a_\disc^G(-)$. Considérons $\Pi_\disc(G)$ le sous-ensemble constitué des sous-quotients irréductibles des représentations induites
\[\I_P^G(\sigma),\,\,\,\,M\in\L^G,P\in\P^G(M),\sigma\in\Pi_\disc(M(\A_F)^1)\]
de $G(\A_F)^1$, où $\sigma$ vérifie la condition suivante : il existe $s\in W^G(M)_\reg$ tel que $s\cdot\sigma=\sigma$. On voit $\Pi_\disc(G)$ comme un espace topologique discret dénombrable. On a alors  
\[J_{\disc}^G(f)=\sum_{\pi\in\Pi_\disc(G)}a_\disc^G(\pi)\tr(\pi(f^1)).\]
Prenons un moment pour clarifier un point qui pourrait être trompeur pour les lecteurs moins familiers avec le domaine. L'ensemble $\Pi_\disc(G(\A_F)^1)$ est l'ensemble des représentations apparaissant dans $L_\disc^2(G(F)\backslash G(\A_F)^1)$, tandis que $\Pi_\disc(G)$ est l'ensemble des représentations contribuant de façon discrète à la formule des traces. On a $ \Pi_\disc(G(\A_F)^1) \subseteq \Pi_\disc(G)$. 

On regarde ensuite l'ensemble
\[\coprod_{M\in\L^G}\coprod_{P\in\P^G(M)}\{\I_{P}^G(\pi_\lambda)\mid \pi\in\Pi_\disc(M),\lambda\in i(a_M^G)^\ast\}.\]
On définit une relation d'équivalence $\sim$ sur cet ensemble par $\I_{P}^G(\pi_\lambda)\sim \I_{P'}^G(\pi_{\lambda'}')$ si et seulement si $M_P=M_{P'}$, $\pi=\pi'$ et $\lambda=\lambda'$. Introduisons le quotient
\begin{equation}\label{eq:defPi(G)triplets}
\Pi(G)=\coprod_{M\in\L^G}\coprod_{P\in\P^G(M)}\{\I_{P}^G(\pi_\lambda)\mid \pi\in\Pi_\disc(M),\lambda\in i(a_M^G)^\ast\}/\sim.    
\end{equation}
On tient à préciser que les éléments de $\Pi(G)$ sont les classes d'équivalence des symboles formels $\I_{P}^G(\pi_\lambda)$. \`{A} chaque élément de $\Pi(G)$ est naturellement associé une classe d'isomorphisme de representations. L'ensemble $\Pi(G)$ est ainsi vu comme un espace topologique homéomorphe à une réunion dénombrable de $\mathbb{R}$-espaces vectoriels de dimension finie, et muni de la mesure $d\pi$ telle que
\[\int_{\Pi(G)}h(\pi)\,d\pi=\sum_{M\in\L^G}|W_0^M||W_0^G|^{-1}\sum_{\pi_1\in\Pi_\disc(M)}\int_{i(a_M^G)^\ast}h(\I_{P}^G(\pi_{1,\lambda}))\,d\lambda\]
pour toute fonction $h\in L^1(\Pi(G))$, avec $P$ un élément quelconque de $\P^G(M)$. Les ingrédients principaux du côté spectral de la formule des traces sont les coefficients 
\begin{equation}\label{eq:defcoeffspeca^Gtriplets}
a^G(\pi)=a_\disc^{M}(\pi_1)r_M^G(\pi_{1,\lambda}),\,\,\,\,M\in\L^G,\pi_1\in\Pi_\disc(M),\lambda\in i(a_M^G)^\ast    
\end{equation}
définis en tout point $\I_{P}^G(\pi_{1,\lambda})=\pi\in \Pi(G)$, et les caractères pondérés 
\begin{equation}\label{eqdef:weightedchara}
J_M^G(\pi,f)=\tr(\mathcal{N}_M(P,\pi)\I_P^G(\pi,f^1))
\end{equation}
puis
\begin{equation}
\J_M^G(\pi,f)=\int_{ i(a_M^G)^\ast}J_M^G(\pi_\lambda,f)\,d\lambda
\end{equation}
définis en tout point $\pi\in \Pi_\unit(M(\A_F)^1)$, avec $P$ un élément quelconque de $\P^G(M)$. L'expression \eqref{eqdef:weightedchara} ne dépend pas de $P\in\P^G(M)$ choisi.

\begin{theorem}[développement fin]\label{thm:devlopfinspec}~{}
\begin{enumerate}
    \item Pour tout $M\in\L^G$ et $\pi\in\Pi(M)$, les expressions $J_M^G(\pi,-)$ et $\J_M^G(\pi,-)$ sont des distributions (i.e. fonctionnelles continues) sur $\CMFL$.
    \item Soit $f\in \CMFL$, alors
    \[J_\spec^G(f)=\sum_{M\in\L^G}|W_0^M||W_0^G|^{-1}\int_{\Pi(M)}a^M(\pi)\J_M^G(\pi,f)\,d\pi.\]
    Toutes les expressions intervenant dans la définition de l'intégrale de droite sont absolument convergentes, et cette intégrale elles-même l'est.
\end{enumerate}
\end{theorem}
\begin{proof}
On prouve le point 2, et on verra au cours de la démonstration que $J_M^G(\pi,-)$ et $\J_M^G(\pi,-)$ convergent absolument et sont des distributions sur $\CMFL$. Amorçons par la vérification formelle de la formule en question. Soient $M\in \L^G$ et $P,Q\in \P^G(M)$, on pose $\mathcal{N}_Q(P,\lambda)\in \End(\I_P^G(R_{M,\disc,\lambda}))$ l'opérateur dont la restriction sur $\I_P^G(\pi_{\lambda})$ vaut $\mathcal{N}_Q(P,\pi,\lambda)$, pour tout $\pi\in \Pi_{\disc}(M(\A_F)^1)$. Nous avons
\begin{align*}
J_\spec^G(f)=\sum_{M\in\L^G}\sum_{L\in\L^G(M)}&|W_0^M||W_0^G|^{-1}\sum_{s\in W^L(M)_\reg}|\det(s-1)_{a_{M}^{L}}|^{-1}\\
&\cdot\int_{i(a_{L}^G)^\ast}\tr\left(\mathcal{M}_{L}(P,\lambda)M_{P|P}(s,0)\I_P^G(R_{M,\disc,\lambda},f^1)\right)\,d\lambda, \end{align*}
soit, grâce à $\mathcal{M}_L(P,\lambda)|_{\I_P^G(\pi_{\lambda})}=\sum_{S\in\L^G(L)}r_L^S(\pi_\lambda)\mathcal{N}_{S}(P,\lambda)|_{\I_P^G(\pi_{\lambda})}$,
\begin{align*}
J_\spec^G(f)=\sum_{M\in\L^G}\sum_{\pi\in\Pi_{\disc}(M(\A_F)^1)}\sum_{L\in\L^G(M)}&\sum_{S\in\L^G(L)}|W_0^M||W_0^G|^{-1}\sum_{s\in W^L(M)_\reg}|\det(s-1)_{a_{M}^{L}}|^{-1}\\
&\cdot\int_{i(a_{L}^G)^\ast}r_L^S(\pi_\lambda)\tr\left(\mathcal{N}_{S}(P,\pi,\lambda)M_{P|P}(s,0)\I_P^G(\pi_{\lambda},f^1)\right)\,d\lambda.  
\end{align*}
Par définition de $a_\disc^M$ on peut réécrire l'égalité en
\begin{equation}\label{eq:devlopfinabCVcheck}
\begin{split}
J_\spec^G(f)=\sum_{L\in\L^G}\sum_{\pi_1\in\Pi_{\disc}(L)}\sum_{S\in\L^G(L)}&|W_0^L||W_0^G|^{-1}\\
&\cdot\int_{i(a_{L}^G)^\ast}a_\disc^L(\pi_1)r_L^S(\pi_{1,\lambda})\tr\left(\mathcal{N}_{S}(Q,\pi_1,\lambda)\I_{Q}^G(\pi_{1,\lambda},f^1)\right)\,d\lambda  
\end{split}    
\end{equation}
avec $Q\in\P^G(L)$ un élément contenant $P\in \P^G(M)$. Observons que $r_L^S(\pi_{1,\lambda})$ ne dépend que de la projection de $\lambda\in i (a_L^G)^\ast$ sur $i (a_L^S)^\ast$, et que pour tout $\lambda_{L}^S\in i (a_L^S)^\ast$ fixé, on a
\begin{equation}\label{eq:weightedchar}
\begin{split}
\int_{i(a_S^G)^\ast}\tr(\mathcal{N}_S(Q,\pi_{1,\lambda_{L}^S+\lambda})&\I_Q^G(\pi_{1,\lambda_{L}^S+\lambda},f^1))\,d\lambda\\
&=\int_{i(a_S^G)^\ast}\tr(\mathcal{N}_S(R,\I_{Q\cap S}^S(\pi_{1,\lambda_{L}^S})_{\lambda})\I_R^G(\I_{Q\cap S}^{S}(\pi_{1,\lambda_{L}^S})_{\lambda},f^1))\,d\lambda\\
&=\J_S^G(\I_{Q\cap S}^{S}(\pi_{1,\lambda_{L}^S}),f)    \end{split}
\end{equation}
avec $R\in\P^G(S)$ un élément contient $Q\in \P^G(L)$. On obtient enfin
\begin{align*}
J_\spec^G(f)&=\sum_{L\in\L^G}\sum_{S\in\L^G(L)}|W_0^L||W_0^G|^{-1}\sum_{\pi_1\in\Pi_{\disc}(L)}\int_{i(a_{L}^S)^\ast}a_\disc^L(\pi_1)r_L^S(\pi_{1,\lambda})\J_S^G(\I_{Q\cap S}^{S}(\pi_{1,\lambda}),f)\,d\lambda   \\
&=\sum_{S\in\L^G}|W_0^S||W_0^G|^{-1}\int_{\Pi(S)}a^S(\pi)\J_S^G(\pi,f)\,d\pi.
\end{align*}

La convergence des formules ci-dessus s'ensuit à travers une investigation des travaux de Müller-Finis-Lapid. Soit $M\in \L^G$. Posons $\F_m^G(M)=\{P\in \F^G(M):\dim a_M^P=m\}$ pour tout entier positif $m$. Soit $H$ un sous-groupe de $G$ stable par conjugaison par $M$. On note $\Sigma^{\text{nd}}(H; A_M)$ l’ensemble des racines non-divisibles de $A_M$ sur $\mathfrak{h}$. Pour $\alpha\in \Sigma^{\text{nd}}(M;A_M)$, on dira que deux sous-groupes paraboliques $P,Q\in \P^G(M)$ sont adjacents le long de $\alpha$, et écrira $P|^{\alpha}Q$, si $\Sigma^{\text{nd}}(P;A_M)\cap (-\Sigma^{\text{nd}}(Q;A_M))=\{\alpha\}$. Pour $H$ une partie de $G$ on notera $\overline{H}$ sa clôture dans $G$. Fixons maintenant $P\in \P^G(M)$. Soit $L\in \L^G(M)$, note par $m=\dim a_L^G$ son co-rang dans $G$. Posons $\mathfrak{B}_{P,L}$ l'ensemble des $m$-uplets $\underline{\beta}=(\beta_1^\vee,\dots,\beta_m^\vee)$ d'éléments de $\Sigma^{\text{nd},\vee}(P;A_M)$ dont les projections sur $a_L$ forment une base de $a_L^G$, et écrivons $\vol(\underline{\beta})$ le covolume dans $a_L^G$ du réseau engendré par cette base. Désignons, pour tout   $\underline{\beta}=(\beta_1^\vee,\dots,\beta_m^\vee)\in\mathfrak{B}_{P,L}$, l'ensemble
\begin{align*}
\Xi_L(\underline{\beta})&=\{(Q_1,\dots,Q_m)\in \F_1^G(M): \beta_i^\vee\in a_M^{Q_i},\forall i=1,\dots,m\}\\
&=\{(\overline{P_1P_1'},\dots,\overline{P_mP_m'}): P_i|^{\beta_i}P_i',\forall i=1,\dots,m\}.
\end{align*}
Pour $P_1|^{\alpha}P_2$ une paire de sous-groupes paraboliques adjacents dans $\P^G(M)$ et $\pi\in \Pi_{\disc}(M(\A_F)^1)$, on pose
\[\delta_{P_1|P_2}^N(\pi,\lambda)=D_{\varpi}N_{P_1|P_2}(\pi,\lambda):\mathcal{A}_{\pi}^2(P_2)\rightarrow \mathcal{A}_{\pi}^2(P_1),\]
avec $\varpi\in a_M^\ast$ un vecteur tel que $\langle \varpi,\alpha^\vee\rangle=1$, et $D_{\varpi}$ la dérivée directionnelle suivant $\varpi$. On sait que $N_{P_1|P_2}(\pi,\lambda)=\Phi_\pi(\langle \lambda,\alpha^\vee\rangle)$ pour $\Phi_\pi$ une fonction méromorphe d'une variable complexe, alors $\delta_{P_1|P_2}^N(\pi,\lambda)=\Phi_\pi'(\langle \lambda,\alpha^\vee\rangle)$. Pour $\mathcal{X}=(\overline{P_1P_1'},\dots,\overline{P_mP_m'})\in \Xi_L(\underline{\beta})$, on note par $\Delta_{\chi}^N(P,\pi,\lambda)$ l'expression
\begin{align*}
\frac{\vol(\underline{\beta})}{m!}N_{P|P_1}(\pi,\lambda) \delta_{P_1|P_1'}^N(\pi,\lambda)&N_{P_1'|P_2}(\pi,\lambda)\\
&\cdots\delta_{P_{m-1}|P_{m-1}'}^N(\pi,\lambda)N_{P_{m-1}'|P_m}(\pi,\lambda)\delta_{P_m|P_m'}^N(\pi,\lambda)N_{P_m'|P}(\pi,\lambda).
\end{align*}
Soit $\underline{\mu}=(\mu_1,\dots,\mu_m)\in (a_M^\ast)^m$ en position générale, autrement dit en dehors d'une réunion finie de hyperplans, on obtient une application $\mathcal{X}_{L,\underline{\mu}}:\mathfrak{B}_{P,L}\rightarrow \F_1^G(M)$ via le procédé expliqué aux pages 220-221 de \cite{FiLa11} (ou aux pages 179-180 de \cite{FiLaMu11}), telle que $\mathcal{X}_{L,\underline{\mu}}(\underline{\beta})\in \Xi_L(\underline{\beta})$ pour tout $\underline{\beta}\in \mathfrak{B}_{P,L}$. Maintenant, les trois propriétés formelles suivantes
\begin{enumerate}
    \item $N_{Q|Q}(\pi,\lambda)=\Id$ ;
    \item $N_{Q_1|Q_2}(\pi,\lambda)\circ N_{Q_2|Q_3}(\pi,\lambda)=N_{Q_1|Q_3}(\pi,\lambda)$ ; et
    \item si $Q_1|^{\alpha}Q_2$, alors $N_{Q_1|Q_2}(\pi,\lambda)=\Phi_\pi(\langle \lambda,\alpha^\vee\rangle)$ pour $\Phi_\pi$ une fonction méromorphe d'une variable complexe ;  
\end{enumerate}
font de $(\mathcal{A}_Q;\I_Q^G(\pi,\lambda)^{J,\tau})_{Q\in \P^G(M)}$ une « famille compatible (\textit{anglais} : compatible family) » dans le langage de \cite{FiLa11}, avec $\mathcal{A}_Q$ la série de Taylor en $\Lambda=0$ de $\mathcal{N}_Q(P,\pi,\lambda,\Lambda)$ restreint à $\I_Q^G(\pi,\lambda)^{J,\tau}$ où $J$ un sous-groupe ouvert compact de $G(\A_{F,\fin})$ et $\tau \in \Pi_{\unit}(K_\infty)$. Ceci étant valable pour tout $J$ et $\tau$, le théorème 8.2 de \cite{FiLa11} dit alors 
\begin{equation}\label{eq:new(G,M)-formulaat0forN}
\mathcal{N}_L(P,\pi,\lambda)=\sum_{\underline{\beta}\in \mathfrak{B}_{P,L}} \Delta_{\mathcal{X}_{L,\underline{\mu}}(\underline{\beta})}^N(P,\pi,\lambda).
\end{equation}

Pour un endomorphisme borné d'un espace de Hilbert, on note $\|-\|_{\text{op}}$ sa norme d'opérateur et $\|-\|_1$ sa norme trace. Posons $\Omega$ (resp. $\Omega_{K_\infty}$) l'opérateur de Casimir de $G(F_\infty)$ (resp. $K_\infty$). Pour tout $\pi \in \Pi_{\disc}(M(\A_F))$ soit $\lambda_\pi$ la valeur propre de Casimir de $\pi_\infty$ ; pour tout $\tau \in \Pi_{\unit}(K_\infty)$ soit $\lambda_\tau$ la valeur propre de Casimir de $\tau$. Fixons pour la suite $J=\prod_{v\in\V_{\fin}} J_v$ un sous-groupe ouvert compact de $G(\A_{F,\fin})$, avec $J_v=K_v$ pour presque tout $v$, tel que $f\in \mathcal{C}(G(\A_F),J)$. L'image de $\I_P^G(\pi_{\lambda},f^1)$ est alors inclus dans $\mathcal{A}_{\pi}^2(P)^{J}$, donc $\I_P^G(\pi_{\lambda},f^1)=0$ sauf si $\pi\in \Pi_{\disc}(M(\A)^1)^{J\cap M(\A_{F,\fin})}$. 

Les ingrédients sont mis en place pour aborder le problème de convergence. Les constantes apparaissent plus tard dans des estimations seront notées $c_i(-)$, indexées par $i$ entier positif, et les objets dont elles dépendent sont placés entre les parenthèses. On commence par vérifier la convergence (pour la norme trace) du caractère pondéré (\ref{eq:weightedchar}) : posons l'opérateur différentiel 
\[\Delta=\Id-\Omega+2\Omega_{K_\infty}.\]
Nous savons, d'après la formule (\ref{eq:new(G,M)-formulaat0forN}) et le fait que $\I_Q^G(\pi_{1,\lambda},-)$ est unitaire sur $\lambda\in i(a_L^G)^\ast$, que pour $k\in \mathbb{N}$,
\begin{equation}\label{eq:majorationderiveJspec}
\begin{split}
\|\tr(\mathcal{N}_S&(Q,\pi_{1,\lambda_L^S+\lambda})\I_Q^G(\pi_{1,\lambda_L^S+\lambda},f^1))\|_1
\\&\leq c_1(G)\|\Delta_{\mathcal{X}}^N(Q,\pi_1,\lambda_L^S+\lambda)\|_{\text{op}}\|\I_Q^G(\pi_{1,\lambda_L^S+\lambda},\Delta^{2k})^{-1}\|_1\|\I_Q^G(\pi_{1,\lambda_L^S+\lambda},(f\ast\Delta^{2k})^1)\|_{\text{op}}\\
&\leq c_1(G)\|\Delta_{\mathcal{X}}^N(Q,\pi_1,\lambda_L^S+\lambda)\|_{\text{op}}\|\I_Q^G(\pi_{1,\lambda_L^S+\lambda},\Delta^{2k})^{-1}\|_1|f\ast \Delta^{2k}|_{L^1(G(\A_F))}
\end{split}    
\end{equation}
pour certains $\underline{\beta}\in\mathfrak{B}_{S,Q}$ et $\mathcal{X}\in \Xi_S(\underline{\beta})$. Pour tout $\tau \in \Pi_{\unit}(K_\infty)$, $\I_Q^G(\pi_{1,\lambda_L^S+\lambda},\Delta)$ agit sur $\mathcal{A}_{\pi_1}^2(Q)^{J,\tau}$ par un scalaire $\mu(\pi_1,\lambda_L^S+\lambda,\tau)$ vérifiant $|\mu(\pi_1,\lambda_L^S+\lambda,\tau)|^2\geq \frac{1}{4}(1+\|\lambda_L^S\|^2+\|\lambda\|^2+\lambda_{\pi_1}^2+\lambda_\tau^2)$ (\cite[équation (5.2)]{FiLaMu11}). 
Supposons à présent $\mathcal{X}=(\overline{Q_1Q_1'},\dots,\overline{Q_mQ_m'})$ avec $ Q_i|^{\alpha_i}Q_i'$, utilisons l'unitarité des opérateurs $N_{Q'|Q''}(\pi_1,-)$ sur $i(a_L^G)^\ast$ et l'inégalité $\|A\|_1\leq (\dim V)\|A\|_{\text{op}}$ pour tout endomorphisme $A$ d'un espace de Hilbert de dimension finie $V$, on est conduit à
\begin{equation*}
\begin{split}
|\J_S^G(&\I_{Q\cap S}^{S}(\pi_{1,\lambda_L^S}),f)|\leq c_2(G,k,f)\times\\
&\sum_{\tau\in\Pi_{\unit}(K_\infty)}\dim(\mathcal{A}_{\pi_1}^2(Q)^{J,\tau})\int_{i(a_S^G)^\ast}|\mu(\pi_1,\lambda_L^S+\lambda,\tau)|^{-2k}\prod_{i=1}^m\|\delta_{Q_i|Q_i'}^N(\pi_1,\lambda_L^S+\lambda)|_{\mathcal{A}_{\pi_1}^2(Q)^{J,\tau}}\|_{\text{op}}\,d\lambda     
\end{split}
\end{equation*}
Par construction des facteurs de normalisation (cf. \cite{Art89I}), on sait qu'il existe un ensemble fini $S(J)\subset \V_{F}$, contenant les places archimédiennes et dépendant uniquement de $J$, tel que $N_{Q'|Q''}(\pi_{1,v},\lambda_L^S+\lambda)|_{\I_{Q''}^G(\pi_{1,v},\lambda_L^S+\lambda)^{J_v}}=\Id$ sauf si $v\in S(J)$. On obtient alors $\delta_{Q_i|Q_i'}^N(\pi_1,\lambda_L^S+\lambda)=\sum_{v\in S(J)}D_{\varpi_i}N_{Q_i|Q_i'}(\pi_v,\lambda_L^S+\lambda)$ avec  $\varpi_i\in a_L^\ast$ tel que $\langle \varpi_i,\alpha_i^\vee\rangle=1$, on est en conséquence amené au cas adjacent. Selon les équations (5.3), (5.4), le lemme 2 puis le lemme 1 de \cite{FiLaMu11}, on voit que  
\begin{equation*}
\begin{split}
\int_{i(a_S^G)^\ast}|\mu(\pi_1,\lambda_L^S+\lambda,\tau)|^{-2k}\prod_{i=1}^m&\|\delta_{Q_i|Q_i'}^N(\pi_1,\lambda_L^S+\lambda)|_{\mathcal{A}_{\pi_1}^2(Q)^{J,\tau}}\|_{\text{op}}\,d\lambda \\
&\leq c_3(G,k,f,J)(1+\|\lambda_L^S\|^2+\lambda_{\pi_1}^2+\lambda_\tau^2)^{-k-c_4(G,J)} 
\end{split}
\end{equation*}
pour $k$ assez grand. Il vient
\begin{equation}\label{eq:weightedcharMAJ}
\begin{split}
|\J_S^G&(\I_{Q\cap S}^{S}(\pi_{1,\lambda_L^S}),f)|\\
&\leq c_5(G,k,f,J)\sum_{\tau\in\Pi_{\unit}(K_\infty)}\dim(\mathcal{A}_{\pi_1}^2(Q)^{J,\tau})(1+\|\lambda_L^S\|^2+\lambda_{\pi_1}^2+\lambda_\tau^2)^{-k-c_4(G,J)}
\end{split}
\end{equation}
pour $k$ est assez grand, ceci converge grâce à \cite[corollaire 0.3]{Mu98}. Il est en outre clair que les caractères pondérés et les $\J_S^G(\pi,-)$ convergent absolument et sont des distributions sur l'espace $\CMFL$ comme $ c_5(G,k,f,J)=c_6(G,k,J)|f\ast\Delta^{2k}|_{L^1(G(\A_F))}$ pour un certain facteur $c_6(G,k,J)$.

On enchaîne sur la convergence absolue de l'équation (\ref{eq:devlopfinabCVcheck}), i.e. $J_{\spec}^G(f)$. \'{E}crivons $\pi_1$ comme l'unique sous-quotient irréductible de $\I_R^L(\pi_{2,\lambda})$ où $R\in\F^L, \pi_2\in\Pi_\disc(M_R(\A_F)^1),\lambda\in i(a_R^L)^\ast$. Nous avons, selon la majoration (\ref{eq:weightedcharMAJ}) et \cite[théorème 5.4]{Mu02}
\begin{equation}
\begin{split}
\int_{i(a_{L}^G)^\ast}|r_L^S(\pi_{1,\lambda})|\times&\|\mathcal{N}_{S}(Q,\pi_1,\lambda)\I_{Q}^G(\pi_{1,\lambda},f^1)\|_1\,d\lambda\\ &\leq    c_7(G,k,f,J)\sum_{\tau\in\Pi_{\unit}(K_\infty)}\dim(\mathcal{A}_{\pi_2}^2(Q)^{J,\tau})(1+\lambda_{\pi_2}^2+\lambda_\tau^2)^{-k-c_8(G,J)}   
\end{split}
\end{equation} 
pour $k$ est assez grand. Or en partant de la définition on voit que $|a_{\disc}^L(\pi_1)|\leq c_9(L)m_{\disc}^{M_R}(\pi_2)$. L'expression de droite de l'équation (\ref{eq:devlopfinabCVcheck}) est alors bornée par
\begin{align*}
c_{10}(G,k,f,J)\sum_{M\in \L^G}\sum_{L\in \L^G(M)}\sum_{\pi\in \Pi_{\disc}(M(\A_F)^1)}&\sum_{\tau\in\Pi_{\unit}(K_\infty)}\\
&m_{\disc}^{M}(\pi)\dim(\mathcal{A}_{\pi}^2(Q_L)^{J,\tau})(1+\lambda_{\pi}^2+\lambda_\tau^2)^{-k-c_8(G,J)},    
\end{align*}
avec $Q_L\in \P^G(L)$ pour tout $L$. On voit que cette expression converge, lorsque $k$ est assez grand, à la faveur de \cite[corollaire 0.3]{Mu98} (cf. aussi \cite[pp. 710-711]{Mu02} pour le passage du contexte adélique à celui classique). Le problème de convergence est achevé. On remarque de surcroît que cela justifie toute interversion somme-intégrale dans le déroulement de la déduction du développement fin.
\end{proof}

\section{Conventions}\label{YDLspecsec:convention}

\`{A} partir du moment présent, jusqu'à la clôture de l'article, on suit les conventions décrites ci-dessous, sauf mention explicite du contraire.

Soit $F$ un corps local ou global de caractéristique 0. Soit $D$ une algèbre à division centrale sur $F$ de degré $d$. Soient $G=\GL_{m,D}$ et $G^\ast=\GL_{n,F}$ avec $n=md$. On fixe $M_0\eqdef\GL_{1,D}^m$ le sous-groupe diagonal. C'est un sous-groupe de Levi minimal de $G$. Soit $\overline{F}$ une clôture algébrique de $F$. Pour un groupe $H$ défini sur $F$ on note $H_{\overline{F}}$ le changement de base $H\times_F\overline{F}$. On fixe $\eta:\g_{\overline{F}}\rightarrow \g_{\overline{F}}^\ast$ une forme intérieure. Par restriction on a que $\eta|_{G_{\overline{F}}}:G_{\overline{F}}\rightarrow G_{\overline{F}}^\ast$ est une forme intérieure. Pour alléger l'écriture, on omettra $\overline{F}$  dans l’indice de $\g_{\overline{F}}\rightarrow \g_{\overline{F}}^\ast$, $G_{\overline{F}}\rightarrow G_{\overline{F}}^\ast$ etc lorsqu’on parle d'un torseur intérieur.

On note $\nu_\g:\g\to \gl_{1,F}$ la norme réduite sur $\g$. Le groupe des $F$-caractères sur $G$ est engendré par $\nu_\g|_G$. On identifie $a_{G}^\ast$ à $\R$ via l'isomorphisme d'espaces linéaires réels $r\nu_{\g}|_G\mapsto r$, et on identifie $a_{G,\C}^\ast$ à $\C$ par l'isomorphisme d'espaces linéaires complexes
\begin{equation}\label{eq:idenaGCtoCI}
\nu_{\g}|_G\otimes s
\mapsto s.      
\end{equation}

On note $e(G)=(-1)^{m(d-1)}$ le signe de Kottwitz de $G$ (\cite[pp. 291-292]{Kott83}). Quand $F$ est global, on se donne un caractère additif non-trivial $\psi:F\backslash\A_F\to \C$ ; quand $F$ est local, on se donne un caractère additif non-trivial $\psi:F\to \C$.

Si $F$ est un corps local archimédien, alors $D=\R$, $\C$, ou $\mathbb{H}$, l'algèbre des quaternions de Hamilton. On note $X\mapsto\overline{X}$ l'involution principale sur $D$. Plus précisément, si $F=\R$ et $D=F$ alors $\overline{a}=a$ pour tout $a\in \R$ ; si $F=\C$ et $D=F$ alors $\overline{a+bi}=a-bi$ pour tout $a,b\in \R$ ; si $F=\R$ et $D=\mathbb{H}$ alors $\overline{a+bi+cj+dk}=a-bi-cj-dk$ pour tout $a,b,c,d\in \R$. On note ensuite $X\mapsto \overline{X}$ l'involution principale et $X\mapsto {}^tX$ la transposée  sur $\gl_m(D)$. Donc si $X=(X_{\alpha,\beta})_{1\leq \alpha,\beta\leq m}\in \gl_m(D)$ alors $\overline{X}=(\overline{X_{\alpha,\beta}})_{1\leq \alpha,\beta\leq m}$ et ${}^tX=(X_{\beta,\alpha})_{1\leq \alpha,\beta\leq m}$.

Si $F$ est local, on fixe $K$ le sous-groupe de $G(F)$ suivant. 
Primo, si $F$ est non-archimédien, alors $K\eqdef\GL_m(\O_D)$ où $\O_D$ est l'anneau des entiers de $D$. C'est un sous-groupe ouvert compact maximal de $G(F)$. Secundo, si $F$ est archimédien, alors $K\eqdef\{g\in \GL_m(D)\mid g\cdot{}^t\overline{g}=\text{Id}\}$. C'est un sous-groupe compact maximal de $G(F)$. Soit $M\in\L^G(M_0)
$, on note $K_M\eqdef K\cap M(F)$. Si $F$ est un corps local non-archimédien, posons $\mathfrak{k}\eqdef\gl_m(\O_D)\subseteq \g(F)$.


\section{Analyse spectrale autour des espaces \texorpdfstring{$\cZ(G(F),s)$}{Z(G(F),s)} et \texorpdfstring{$\cZ(G(\A_{F}),s)$}{Z(G(AF),s)}}\label{sec:analyseespaceGJ}

Soit $F$ un corps local ou global de caractéristique 0. Lorsque $F$ est un corps $p$-adique on note $q$ le cardinal du corps résiduel.

Pour un groupe topologique localement compact et séparé $H$, on désigne par $\S(H)$ l'espace de Schwartz-Bruhat de $H$, muni de sa topologie usuelle (cf. \cite[page 61]{Bruhat1961}).

Nous fixons les mesures sur les groupes impliqués dans les discussions suivantes. Lorsqu'une mesure n'est pas explicitement précisée, cela indique que le choix de celle-ci ne revêt pas d'importance pour notre objectif.

\subsection{Rappel des fonctions \texorpdfstring{$L$}{L}, \texorpdfstring{$\epsilon$}{e} et \texorpdfstring{$\gamma$}{r} de Godement-Jacquet}\label{subsec:reviewGJtheory}
Cette partie a trait à un rappel de la théorie de \cite{GoJa72}, dont on modifiera légèrement la présentation. On ne s'intéresse qu'aux représentations unitaires. On va toujours identifier l'espace vectoriel sous-jacent d'une représentations unitaire à son contragrédient via son produit scalaire.

\subsubsection{Théorie locale}\label{subsubsec:GJtheorielocale}
Soit $F$ un corps $p$-adique. On note $\langle-,-\rangle$ la forme bilinéaire canonique de $\g(F)$. Autrement dit $\langle X,Y\rangle$ est la trace réduite de $XY$. On prend sur $\g(F)$ la mesure auto-duale par rapport à la transformée de Fourier induite par $\psi(\langle-,-\rangle)$. Pour une fonction $\Phi\in \S(\g(F))$ on désigne par $\widehat{\Phi}\in \S(\g(F))$ sa transformée de Fourier, à savoir
\[\widehat{\Phi}(Y)=\int_{\g(F)}\Phi(X)\psi(\langle X,Y\rangle)\,dX,\,\,\,\,\forall Y\in \g(F),\] 

Soient $(\pi,V)$ une représentation admissible unitaire de $G(F)$, non nécessairement irréductible, et $v,v^\vee\in V$. Pour toute fonction $\Phi\in \S(\g(F))$ définissons la fonction zêta locale 
\[Z(v,v^\vee,\Phi,s)=Z((v,v^\vee)_\pi,\Phi,s)\eqdef\int_{G(F)}\Phi(g)\langle \pi(g)v,v^\vee\rangle|\nu_\g(g)|^{s+\frac{n-1}{2}}\,dg.\]

\subsubsection{Cas non-archimédien}

Soit $F$ un corps $p$-adique.

\begin{theorem}[{{\cite[proposition 1.1, théorème 3.3]{GoJa72}}}]\label{thm:GJzetalocalnon-archi} Soit $(\pi,V)$ une représentation irréductible unitaire (donc admissible) de $G(F)$.
\begin{enumerate}
    \item L'intégrale $Z(v,v^\vee,\Phi,s)$ converge absolument sur le demi-plan $\{s\in\C\mid\Re(s)>\frac{n+1}{2}-d\}$ pour tout $v,v^\vee$ et $\Phi$. 
    \item Pour tout $v,v^\vee\in V$ on a $Z(v,v^\vee,\Phi,s)\in \C(q^{-s})$, et la famille $(Z(v,v^\vee,\Phi,s))_{v,v^\vee,\Phi}$ engendrent un idéal fractionnaire de $\C
[q^{\pm s}]$ de la forme $\frac{1}{P_\pi(q^{-s})}\C[q^{\pm s}]$ pour $P_\pi$ un polynôme avec $P_\pi(0)=1$. On appelle $\frac{1}{P_\pi(q^{-s})}$ la fonction $L$ locale de $\pi$ et on la note $L(\pi,s)$.
    \item Il y a une équation fonctionnelle pour tout $v,v^\vee\in V$ et $\Phi\in \S(\g(F))$ : 
    \[Z(v^\vee,v,\widehat{\Phi},1-s)=e(G)\gamma(\pi,\psi,s)Z(v,v^\vee,\Phi,s)\]
    pour une fonction rationnelle $\gamma(\pi,\psi,s)$, appelée le facteur gamma local de $\pi$.
    \item Soit 
    \[\epsilon(\pi,\psi,s)\eqdef \gamma(\pi,\psi,s)\frac{L(\pi,s)}{L(\pi^\vee,1-s)}\]
    le facteur epsilon local de $\pi$, qui est un monôme de $q^{-s}$. Alors pour tout $v,v^\vee\in V$ et $\Phi\in \S(\g(F))$,
    \[\frac{Z(v^\vee,v,\widehat{\Phi},1-s)}{L(\pi^\vee,1-s)}=e(G)\epsilon(\pi,\psi,s)\frac{Z(v,v^\vee,\Phi,s)}{L(\pi,s)}.\]
    Pour rappel $\pi^\vee$ est la représentation contragrédiente de $\pi$. On a aussi $\epsilon(\pi,\psi,s)\epsilon(\pi^\vee,\psi,1-s)=z_\pi(-1)$ avec $z_\pi$ le caractère central de $\pi$.
\end{enumerate}
\end{theorem}

On voit, selon les points 1 et 2, que toute fonction zêta locale $Z(v,v^\vee,\Phi,s)$ admet un prolongement méromorphe au plan complexe.

Lors même que $\pi$ ne serait pas irréductible, l'énoncé du théorème garde un sens quand elle admet un caractère central.

\begin{theorem}[{{\cite[théorème 3.4, proposition 3.5, corollaire 3.6]{GoJa72}}}]~{}
\begin{enumerate}
    \item Soient $L=\prod_i L_i$ un sous-groupe de Levi de $G$ avec $L_i$ de la forme $\GL_{m_i,D}$, $\sigma_i$ une représentation admissible unitaire (non nécessairement irréductible) de $L_i(F)$ admettant un caractère central, tels que le théorème précédent soit valable pour $\sigma_i$. Soit $P\in\P^G(L)$, alors le théorème précédent est valable pour $\I_P^G(\otimes_i \sigma_i)$. De surcroît, la formation de la fonction $L$, du facteur epsilon et du facteur gamma est compatible à l'induction parabolique, autrement dit 
    \begin{align*}
    L(\I_P^G(\otimes_i \sigma_i),s)=\prod_{i}L(\sigma_i,s),\,\,\epsilon(\I_P^G(\otimes_i \sigma_i),\psi,s)&=\prod_{i}\epsilon(\sigma_i,\psi,s)\\
    \,\,\text{et}\,\,\gamma(\I_P^G(\otimes_i \sigma_i),\psi,s)&=\prod_{i}\gamma(\sigma_i,\psi,s).    
    \end{align*}
    \item Soit $\pi$ une représentation unitaire admettant un caractère central et pour laquelle le théorème précédent reste valable. Soit $\lambda\in ia_{G,\C}^\ast$. Alors le théorème est également valable pour $\pi_{\lambda}$. La formation de la fonction $L$, du facteur epsilon et du facteur gamma est compatible à la torsion au sens où
    \[L(\pi_\lambda,s)=L(\pi,s+\lambda),\epsilon(\pi_\lambda,s)=\epsilon(\pi,s+\lambda),\,\,\text{et}\,\,\gamma(\pi_\lambda,s)=\gamma(\pi,s+\lambda).\]
    \item Soit $\pi$ une représentation unitaire admettant un caractère central et pour laquelle le théorème précédent reste valable. Alors le théorème est également valable pour tout sous-quotient irréductible $\pi'$. De plus 
    \[L(\pi',s)L(\pi,s)^{-1}\in \C[q^{-s}],\,\,\text{et}\,\,\gamma(\pi', \psi, s)=\gamma(\pi, \psi, s).\]
\end{enumerate}  
\end{theorem}

\subsubsection{Cas archimédien}
Soit $F=\R$ ou $\C$.

Introduisons $\Gamma_\R(s)=\pi^{-\frac{1}{2}s}\Gamma\left(\frac{1}{2}s\right)$ et $\Gamma_\C(s)=2(2\pi)^{-s}\Gamma(s)$ où $\Gamma$ est la fonction gamma d'Euler et $\pi$ la constante d’Archimède. Si $F=\R$, un facteur eulérien est une fonction de la forme 
\[P(s)\prod_{i\text{ fini}}\Gamma_\R(s+s_i)\prod_{j\text{ fini}}\Gamma_\C(s+s_j)\]
où $P$ est un polynôme et $s_i$ et $s_j$ sont certaines constantes ; si $F=\C$, un facteur eulérien est une fonction de la forme 
\[P(s)\prod_{j\text{ fini}}\Gamma_\C(s+s_j)\]
où $P$ est un polynôme et $s_j$ sont certaines constantes. Introduisons ensuite le sous-espace 
\begin{equation}\label{YDLspeceq:defS-Bspacearchiconstraint}
\S(\g(F),\psi)\subseteq \S(\g(F))    
\end{equation} 
définit comme suit : si $\psi(\cdot)=\exp(-2i\pi\tr_{F/\R}(a_\psi\cdot))$ pour $a_\psi\in F\setminus\{0\}$, alors $\S(\g(F),\psi)$ est l'espace des fonctions de la forme $X\mapsto P(X)\exp(-[F:\R]d\pi|a_\psi|_F^{\frac{1}{[F:\R]}}\tau_\g(X\cdot{}^t\overline{X}))$ où $P\in \text{Sym}(\Res_{F/\R}\g^\vee)(\R)$ (i.e. une fonction polynômiale sur $\g(F)$ vu comme espace vectoriel réel). L'espace $\S(\g(F),\psi)$ est stable par la transformée de Fourier et l'action de l'algèbre enveloppante $\mathcal{U}(\g(F)_\C)$. 

\begin{theorem}[{{\cite[proposition 1.1, théorème 8.7]{GoJa72}, \cite[équation (4.3.1), (4.5) proposition]{Ja79}}}]\label{thm:GJzetalocalarchi} Soient $(\pi,V)$ une représentation irréductible unitaire (donc admissible) de $G(F)$.
\begin{enumerate}
    \item L'intégrale $Z(v,v^\vee,\Phi,s)$ converge absolument sur le demi-plan $\{s\in\C\mid\Re(s)>\frac{n+1}{2}-d\}$ pour tout $v,v^\vee$ et $\Phi\in \S(\g(F))$. 
    Pour $s\in\C$ avec $\Re(s)>\frac{n+1}{2}-d$, l'application $\Phi\mapsto Z(v,v^\vee,\Phi,s)$ est une distribution tempérée, dépendant holomorphiquement de $s$. 
    \item Il existe un facteur eulérien $L(\pi,s)$ appelé la fonction $L$ locale de $\pi$ vérifiant que si $\Xi(v,v^\vee,\Phi,s)\eqdef Z(v,v^\vee,\Phi,s)L(\pi,s)^{-1}$, alors pour toute fonction $\Phi\in \S(\g(F),\psi)$ l'expression $\Xi(v,v^\vee,\Phi,s)$ est de la forme
    \[P(s)|a_\psi|_F^{-\frac{1}{2}sn}\] 
    avec $P$ un polynôme. La fonction $L(\pi,s)$ s'annule en aucun point de $\C$.
    \item Il y a une équation fonctionnelle pour tout $v,v^\vee$ et $\Phi\in \S(\g(F),\psi)$ :
    \[Z(v^\vee,v,\widehat{\Phi},1-s)=e(G)\gamma(\pi,\psi,s)Z(v,v^\vee,\Phi,s)\]
    pour une fonction rationnelle $\gamma(\pi,\psi,s)$ appelée le facteur gamma local de $\pi$.
    \item Soit 
    \[\epsilon(\pi,\psi,s)\eqdef \gamma(\pi,\psi,s)\frac{L(\pi,s)}{L(\pi^\vee,1-s)}\]
    le facteur epsilon local de $\pi$, qui est une fonction exponentielle de $s$. Alors
    \[\frac{Z(v^\vee,v,\widehat{\Phi},1-s)}{L(\pi^\vee,1-s)}=e(G)\epsilon(\pi,\psi,s)\frac{Z(v,v^\vee,\Phi,s)}{L(\pi,s)}.\]
    On a aussi $\epsilon(\pi,\psi,s)\epsilon(\pi^\vee,\psi,1-s)=z_\pi(-1)$ avec $z_\pi$ le caractère central de $\pi$.
\end{enumerate}
\end{theorem}


On voit, selon les points 1 et 2, que toute fonction zêta locale $Z(v,v^\vee,\Phi,s)$ admet un prolongement méromorphe au plan complexe.

\subsubsection{Théorie globale}\label{subsubsec:globalGJtheory}
Soit $F$ un corps de nombres. Soit $K_\infty\eqdef \prod_{v\in \V_{\infty}} K_v$. Soit $\pi$ un $G(\A_{F,\fin})\times (\g(F_\infty)_\C,K_\infty)$-module admissible dont le caractère central est trivial sur $F^\times$. On a $\pi=\bigotimes_{v\in\V_F}\pi_v$. Le facteur epsilon global $\epsilon(\pi,s)\eqdef\prod_{v\in\V_F}\epsilon_v(\pi_v,\psi_v,s)$ est bien défini (\cite[proposition 6.12 (2)]{GoJa72}), et ne dépend pas du choix de $\psi$. On sait aussi qu'il existe un nombre réel $s(G)$ indépendant de $\pi$ tel que  la fonction $L$ globale $L(\pi,s)\eqdef\prod_{v\in \V_F}L(\pi_v,s)$ soit bien définies sur le demi-plan $\{s\in\C\mid \Re(s)>s(G)\}$ (\cite[lemme 10.20]{GoJa72}).

On note $\langle-,-\rangle=\prod_{v\in\V_F}\langle-,-\rangle_v$ la forme bilinéaire canonique de $\g(\A_F)$, où $\langle-,-\rangle_v$ est la forme bilinéaire canonique de $\g(F_v)$. On prend sur $\g(\A_F)$ la mesure auto-duale par rapport à la transformée de Fourier induite par $\psi(\langle-,-\rangle)$. Pour une fonction $\Phi\in \S(\g(\A_F))$ on désigne par $\widehat{\Phi}\in \S(\g(\A_F))$ sa transformée de Fourier, à savoir
\[\widehat{\Phi}(Y)=\int_{\g(F)}\Phi(X)\psi(\langle X,Y\rangle)\,dX,\,\,\,\,\forall Y\in \g(\A_F),\] 

\begin{theorem}[{{\cite[théorème 12.4, théorème 13.8, corollaire 13.8, proposition 13.9]{GoJa72} pour le cas cuspidal ; \cite[(6.2) théorème]{Ja79} pour le cas général}}] Soit $(\pi,V)$ une représentation automorphe de $G(\A_F)^1$. 
\begin{enumerate}
    \item Définissons la fonction zêta globale, pour $v,v^\vee\in V$ et $\Phi\in \S(\g(\A_F))$, 
    \[Z(v,v^\vee,\Phi,s)\eqdef\int_{G(\A_F)}\Phi(g)\langle \pi(g)v,v^\vee\rangle|\nu_\g(g)|_{\A_F}^{s+\frac{n-1}{2}}\,dg.\]
    Alors elle converge absolument sur le demi-plan $\{s\in\C\mid\Re(s)>\frac{n+1}{2}\}$, et admet un prolongement méromorphe au plan complexe, entière si $\pi$ est cuspidal et $m>1$. Il y a une équation fonctionnelle 
    \[Z(v,v^\vee,\Phi,s)=Z(v^\vee,v,\widehat{\Phi},1-s)\]
    \item La fonction $L(\pi,s)$ admet un prolongement méromorphe au plan complexe, entière si $\pi$ est cuspidal et $m>1$. Elle a un nombre fini de pôles et est bornée à l'infini dans toute bande verticale de largeur finie. Il y a une équation fonctionnelle 
    \[L(\pi,s)=\epsilon(\pi,s)L(\pi^\vee,1-s).\]
\end{enumerate}
\end{theorem}

\subsection{Rappel de l'isomorphisme de Satake}\label{subsec:rappelSatake}

Soit dans cette partie $F$ un corps $p$-adique. Supposons que $G$ est déployé, i.e. $D=F$. Il convient d'observer que $M_0$ est un sous-tore déployé maximal de $G$. 

Fixons les normalisations des mesures de cette sous-section. On munit $K$ de la mesure de Haar telle que $\vol(K)=1$. On munit, pour tout $P\in \P^G$, le groupe $N_P(F)$ de la mesure de Haar telle que $\vol(K\cap N_P(F);N_P(F))=1$, et le groupe $M_P(F)$ de la mesure de Haar telle que $\vol(K\cap M_P(F);M_P(F))=1$.

Soit $J$ un sous-groupe ouvert compact de $G(F)$. Posons $\H(G(F))$ l'algèbre des fonctions continues à support compact munie du produit de convolution comme loi multiplicative, et $\H(G(F)\ssla J)$ la sous-algèbre des fonctions bi-invariantes par $J$. Soient $M\in \L^G(M_0)$ et $P\in \P^G(M)$. On définit la fonction module pour tout $p \in P(F)$ par $\delta_P (p)=|\det(\Ad(p);\mathfrak{n}_P (F))|$. Le terme constant 
\[f\mapsto (f_P: m\mapsto \delta_P(m)^{1/2}\int_K\int_{N_P(F)}f(k^{-1}mnk)\,dn\,dk,\,\,\,\,\forall m\in M(F))\] 
fournit une application d'algèbres 
\begin{equation}\label{eq:constterm-Satake}
 \H(G(F)\ssla K)\rightarrow \H(M(F)\ssla K_M).
\end{equation}
On a la formule d'adjonction : pour toute $\pi$ représentation admissible de $M(F)$,
\[\langle\I_P^G(\pi)\otimes \lambda,f\rangle=\langle \pi\otimes \lambda,f_P\rangle,\,\,\forall \lambda\in a_{G,\C}^\ast.\]

On affuble un groupe réductif sur $F$ d'un chapeau lorsqu'on parle de son dual de Langlands (\cite[page 29]{Bor79}). On abrège $H_{M_0}$ en $H_0$ et $a_{M_0}$ en $a_0$. Les représentations irréductibles ($K_{M_0}$-)non-ramifiées de $M_0(F)$ sont données par $m\mapsto \exp(\langle \lambda, H_0(m)\rangle)$ avec $\lambda\in a_{0,\C}^\ast$ (on tient à préciser que l'exponentielle $\exp=\exp_F$ est de base $q^{-1}$), 
et évidemment cette représentation ne dépend que de l'image de $\lambda$ par 
\begin{equation}\label{eq:Satakespaces}
\lambda\in a_{0,\C}^\ast=X^\ast(M_0)\otimes_\Z\C=\Lie(\widehat{M_0})\xrightarrow{\exp}\widehat{M_0},    
\end{equation} 
On identifie dans la suite l'ensemble des classes d’isomorphismes des représentations irréductibles non-ramifiées de $M_0(F)$ à $a_{0,\C}^\ast/\ker(\exp)$.

Notons $W(\widehat{G},\widehat{M_0})$ le groupe de Weyl absolue de $\widehat{M_0}$ dans $\widehat{G}$. Il y a une bijection entre $  \widehat{M_0}/W(\widehat{G},\widehat{M_0})$ et l'ensemble des classes des représentations irréductible non-ramifiées de $G(F)$, donnée par l'association à la $W(\widehat{G},\widehat{M_0})$-orbite de $\exp(\lambda)\in \widehat{M_0}$ la classe de l'unique sous-quotient non-ramifié de $\I_{P_0}^G(\lambda)$. Soit $\pi$ une représentation non-ramifiée de $G(F)$. Soit $\lambda\in a_{0,\C}^\ast$ tel que $\pi$ est l'unique sous-quotient non-ramifié de $\I_{P_0}^G(\lambda)$. On note 
\[z_\pi\in \mathcal{X}^{G}\eqdef\widehat{M_0}/W(\widehat{G},\widehat{M_0})\]
l'image de $\lambda$ dans $\mathcal{X}^G$ et on l'appellera le paramètre de Hecke de $\pi$. 
En parallèle, il y a un morphisme d'algèbres $\H(M_0(F)\ssla K_{M_0})\to \C[X^{\ast}(\widehat{M_0})]$. En le composant avec \eqref{eq:constterm-Satake} où $M=M_0$, on obtient l'isomorphisme de Satake 
\[\mathcal{S}:\H(G(F)\ssla K)\xrightarrow{\sim} \C[X^{\ast}(\widehat{M_0})]^{W(\widehat{G},\widehat{M_0})}=\Gamma(\mathcal{X}^G,\O_{\mathcal{X}^G}),\]
où $\Gamma(\mathcal{X}^G,\O_{\mathcal{X}^G})$ est l'anneau des fonctions régulières de la variété complexe $\mathcal{X}
^G$. On a la formule d'adjonction
\[\langle \pi,f\rangle=\mathcal{S}f(z_\pi)\]
pour tout $f\in\H(G(F)\ssla K)$. 

\begin{proposition}\label
{prop:Satakeunrcalcul} 
~{}
\begin{enumerate}
    \item Posons $1_{\mathfrak{k}}$ la fonction indicatrice de $\mathfrak{k}$. Soit $(\pi,V)$ une représentation non-ramifiée de $G(F)$. Alors pour tout $v,v^\vee\in V$ vérifiant $\langle v,v^\vee\rangle\not=0$ nous avons
    \[Z(v,v^\vee,1_{\mathfrak{k}},s)\in \C^\times \cdot L(s,\pi)\]
    \item Soit $\pi$ une représentation non-ramifiée de $G(F)$. Soit $\lambda\in a_{0,\C}^\ast$ tel que $\pi$ est l'unique sous-quotient non-ramifié de $\I_{P_0}^G(\lambda)$. Identifions $a_{0,\C}^\ast=X^\ast(M_0)\otimes_\Z\C$ à $\C^n$ en envoyant $\nu_{\gl_{1,F}}|_{\GL_{1,F}}\otimes 1$ sur $1$, et écrivons $\lambda=(\lambda_1,\dots,\lambda_n)\in\C^n$, alors
    \[L(s,\pi)=\prod_{i=1}^n(1-q^{-s-\lambda_i})^{-1}.\]

    \item L'adhérence $\mathcal{X}_{\unit,\text{ur}}^G$ du sous-ensemble des paramètres de Satake des représentations unitaires non-ramifiées est une partie compacte de $\mathcal{X}^G$.  Quantitativement, si l'on identifie $a_{0,\C}^\ast=X^\ast(M_0)\otimes_\Z\C$ à $\C^n$ en envoyant $\nu_{\gl_{1,F}}|_{\GL_{1,F}}\otimes 1$ sur $1$ alors 
    \[\mathcal{X}_{\unit,\text{ur}}^G\subseteq\exp\left(\{s\in\C\mid\Re(s)\in \left[-\frac{n-1}{2},\frac{n-1}{2}\right]\}^n\right)/W(\widehat{G},\widehat{M_0}).\]
\end{enumerate}
\end{proposition}
\begin{proof}~{}
\begin{enumerate}
    \item Cela est \cite[lemme 6.10]{GoJa72}.
    \item Notons, pour $i=1,\dots,n$, $\lambda_i:x\mapsto |\nu_{\gl_{1,F}}(x)|^{\lambda_i}$ une représentation de $\GL_1(F)$. Alors
    \[L(s,\pi)=\prod_{i=1}^mL(s,\lambda_i)=\prod_{i=1}^m(1-q^{-s-\lambda_i})^{-1}.\]
    Les égalités dans l'ordre sont respectivement \cite[proposition 6.9]{GoJa72} et \cite[proposition 4.7]{GoJa72}.
    \item Cela est \cite[proposition 6.12 (4)]{GoJa72}.\qedhere
\end{enumerate}   
\end{proof}

\subsection{Espaces \texorpdfstring{$\cZ(G(F),s)$}{Z(G(F),s)} et \texorpdfstring{$\cZ(G(\A_{F,S}),s)$}{Z(G(A(F,S),s))}}\label{subsec:GJspace}
Soit dans ce numéro $F$ un corps local ou global de caractéristique 0.

On rappelle la topologie sur l'espace de Schwartz-Bruhat $\S(\g(F))$ lorsque $F$ est local : quand $F$ est non-archimédien, il s'agit de la topologie discrète ; et quand $F$ est archimédien, il s'agit de la topologie engendrée par les semi-normes $\|-\|_{P,Q}$ pour $P$ et $Q$ polynômes à $\dim_\R \g(F)=[F:\R]^{n^2}$ variables réelles et à coefficients complexes, avec 
\begin{equation}\label{YDLspeceq:defSchseminormFarchi}
\|f\|_{P,Q}=\sup_{x\in \g(F)}\left|P(x)Q\left(\frac{\partial}{\partial x}\right)f(x)\right|.    
\end{equation}

\begin{definition}[Domaines de convergence $\D(G)$, $\D_\temp(G)$, $\D_\unit(G)$, et $\D_\adel(G)$] ~{}
\begin{enumerate}
    \item Définissons $\D(G)\eqdef\{-\frac{n-1}{2}\}\cup\{s\in\C\mid \Re(s)>-\frac{n-1}{2}\}$. 
    \item Soit $F$ local. Définissons  $\D_\temp(G)\eqdef\{s\in\C \mid \Re(s)>\frac{-d+1}{2}\}$ et $\D_\unit(G)\eqdef \{s\in\C\mid \Re(s)>\frac{n+1}{2}-d\}$.
    \item Soit $F$ global. Définissons $\D_\adel(G)\eqdef\{s\in \C \mid \Re(s)>\frac{n+1}{2}\}$. 
\end{enumerate}
\end{definition}

Pour $F$ local on a $\D_\unit(G)\subseteq \D_\temp(G)\subseteq\D(G)$ ; pour $F$ global on a $\D_\adel(G)\subseteq \D(G)$. On a identifié $a_{G,\C}^\ast$ à $\C$ par l'isomorphisme \eqref{eq:idenaGCtoCI}. Les ensembles $\D(G)$, $\D_\temp(G)$, $\D_\unit(G)$, et $\D_\adel(G)$ sont donc des domaines de la variété complexe $a_{G,\C}^\ast$. 

\subsubsection{Espaces $\cZ(G(F),s)$ et $\cZ(G(\A_{F,S}),s)$}

\begin{definition}\label{def:espacedeGJ} Par convention, $0^0=1$.  
\begin{enumerate}
    \item Supposons que $F$ est local. Soient $s\in \D(G)$ et $f\in \S(\g(F))$. Posons $f_s$ la fonction $X\in \g(F)\mapsto f(X)|\nu_{\g}(X)|_F^{s+\frac{n-1}{2}}\in \C$. Définissons $\cZ(\g(F),s)\eqdef \{f_s \mid f\in \S(\g(F))\}$. C'est un espace vectoriel complexe.
    \item Supposons que $F$ est global. Soit $S\subseteq \V_F$. On note $\A_{F,S}\eqdef \prod_{v\in S}'F_v$ le produit restreint des $F_v$ ($v\in S$) par rapport à la famille $(\O_v)_{v\in S\cap \V_{\fin}}$. Soient $s\in \D(G)$, $f\in \S(\g(\A_{F,S}))$. Posons $f_s$ la fonction $X\in \g(\A_{F,S})\mapsto f(X)|\nu_{\g}(X)|_{\A_{F,S}}^{s+\frac{n-1}{2}}\in \C$. Définissons $\cZ(\g(\A_{F,S}),s)\eqdef \{f_s \mid f\in \S(\g(\A_{F,S}))\}$. C'est un espace vectoriel complexe.
\end{enumerate}
\end{definition}

On munit $\D(G)$ de la topologie induite par celle de $\C$. 

\begin{definition}[Topologie]\label{def:GJuniontopo} 
Soit $I$ une sous-partie de $\D(G)$. 
\begin{enumerate}
    \item Supposons que $F$ est local. On définit la topologie sur 
    \[\bigcup_{s\in I}\cZ(\g(F),s)\,\,\,\,\text{(réunion dans l'espace des fonctions sur $\g(F)$)}\] 
    en poussant en avant la topologie produit sur $I\times \S(\g(F))$. Autrement dit un sous-ensemble $E$ est un ouvert de $\cup_{s\in I}\cZ(\g(F),s)$ si et seulement si $\iota^{-1}(E)$ est un ouvert dans $I\times \S(\g(F))$ avec $\iota :(s,f)\in I\times \S(\g(F))\mapsto f_s\in \cup_{s\in I}\cZ(\g(F),s)$.
    \item Supposons que $F$ est global. On définit la topologie sur 
    \[\bigcup_{s\in I}\cZ(\g(\A_F),s)\,\,\,\,\text{(réunion dans l'espace des fonctions sur $\g(\A_F))$}\]
    en poussant en avant la topologie produit sur $I\times \S(\g(\A_F))$. Autrement dit un sous-ensemble $E$ est un ouvert de $\cup_{s\in I}\cZ(\g(\A_F),s)$ si et seulement si $\iota^{-1}(E)$ est un ouvert dans $I\times \S(\g(\A_F))$ avec $\iota :(s,f)\in I\times \S(\g(\A_F))\mapsto f_s\in \cup_{s\in I}\cZ(\g(\A_F),s)$.
\end{enumerate}    
\end{definition}

    

Supposons que $F$ est local (resp. global). On observe que $G(F)$ (resp. $G(\A_F)$) est un sous-ensemble dense de $\g(F)$ (resp. $\g(\A_F)$). 
Toute fonction continue sur $\g(F)$ (resp. $\g(\A_F)$) est donc complètement déterminée par sa restriction sur $G(F)$ (resp. $G(\A_F)$). Pour tout $I\subseteq \D(G)$, l'application $\cup_{s\in I}\cZ(\g(F),s)\to\cup_{s\in I}\cZ(\g(F),s)|_{G(F)}$ (resp. $\cup_{s\in I}\cZ(\g(\A_F),s)\to\cup_{s\in I}\cZ(\g(\A_F),s)|_{G(\A_F)}$) est alors bijective. On obtient de ce fait une topologie sur $\cup_{s\in I}\cZ(\g(F),s)|_{G(F)}$ (resp. $\cup_{s\in I}\cZ(\g(\A_F),s)|_{G(\A_F)}$) par cette bijection.

Supposons que $F$ est local (resp. global). Le groupe $G(F)$ (resp. $G(\A_F)$) agit sur l'espace de Schwartz-Bruhat $\S(\g(F))$ (resp. $\S(\g(\A_F))$) à gauche et à droite via
\[(g_1\cdot f\cdot g_2):X\in\g(F)\mapsto f(g_1Xg_2^{-1})\in \C,\,\,\,\,\text{où $g_1,g_2\in G(F)$, $f\in \S(\g(F))$, $X\in \g(F)$}.\]
(resp. via
\[(g_1\cdot f\cdot g_2):X\in\g(\A_F)\mapsto f(g_1Xg_2^{-1})\in \C,\,\,\,\,\text{où $g_1,g_2\in G(\A_F)$, $f\in \S(\g(\A_F))$, $X\in \g(\A_F)$)}.\]

Supposons que $F$ est $p$-adique. Soit $J\subseteq G(F)$ un sous-groupe ouvert compact, on note $\S(\g
(F)\ssla J)$ l'espace des fonctions Schwartz-Bruhat bi-invariantes par $J$. De même $\cZ(\g(F)\ssla J,s)\eqdef\{f_s\mid f\in \S(\g(F)\ssla J)\}$. On munit $\H(G(F))$ de la topologie discrète. Il importe d'observer qu'il existe des flèches naturelles entre les espaces vectoriels topologiques 
\[\H(G(F))\hookrightarrow \S(\g(F))\xrightarrow[ \times |\nu_{\g}|_F^{s+\frac{n-1}{2}}]{\sim} \cZ(\g(F),s),\]
ainsi que 
\[\H(G(F)\ssla J)\hookrightarrow \S(\g(F)\ssla J)\xrightarrow[ \times |\nu_{\g}|_F^{s+\frac{n-1}{2}}]{\sim} \cZ(\g(F)\ssla J,s).\]

Supposons que $F$ est global. Soit $J\subseteq G(\A_{F,\fin})$ un sous-groupe ouvert compact, on écrit $\S(\g
(\A_F)\ssla J)$ l'espace des fonctions Schwartz-Bruhat bi-invariantes par $J$. De même $\cZ(\g(\A_F)\ssla J,s)\eqdef\{f_s\mid f\in \S(\g(\A_F)\ssla J)\}$. Il importe d'observer qu'il existe des isomorphismes naturels d'espaces vectoriels topologiques 
\[\S(\g(\A_F))\xrightarrow[ \times |\nu_{\g}|_{\A_F}^{s+\frac{n-1}{2}}]{\sim} \cZ(\g(\A_F),s),\,\,\,\,\S(\g(\A_F)\ssla J)\xrightarrow[ \times |\nu_{\g}|_{\A_F}^{s+\frac{n-1}{2}}]{\sim} \cZ(\g(\A_F)\ssla J,s).\]

Suppsons que $F$ est $p$-adique. Soit $s\in\D(G)$. Le groupe $G(F)$ agit sur l'espace $\cZ(\g(F),s)$ à droite via
\[(f_s\cdot g):X\in\g(F)\mapsto f_s(Xg^{-1})\in \C,\,\,\,\,\text{où $g\in G(F)$, $f\in \S(\g(F))$, $X\in \g(F)$}.\]
Ce qui procure une structure de $\H(G(F))$-module à l'espace $\cZ(\g(F),s)$. On résume ce fait dans la proposition suivante.

\begin{proposition}\label{prop:GJspaceoverHeckealg}
Soient $s\in\D(G)$, $h\in \H(G(F))$ et $f\in\S(\g(F))$, le produit de convolution
\[h\ast f_s:X\in\g(F)\longmapsto \int_{G(F)}h(g)f_s(Xg^{-1})\,dg\]
converge et définit une fonction dans $\cZ(\g(F),s)$. En d'autres mots l'espace $\cZ(\g(F),s)$ acquiert une structure de module sur l'algèbre de Hecke via la convolution
\begin{equation}\label{eq:GJspaceoverHeckealgconv}
\H(G(F))\times \cZ
(\g(F),s)\overset{\ast}{\longrightarrow}\cZ(\g(F),s).    
\end{equation}
Pour tout $J\subseteq G(F)$ sous-groupe ouvert compact, nous avons
\[\H(G(F)\ssla J)\times \cZ(\g(F)\ssla J,s)\overset{\ast}{\longrightarrow}\cZ(\g(F)\ssla J,s).\]
\end{proposition}

\subsubsection{Espace de Schwartz-Harish-Chandra}

Soit dans cette partie $F$ local. Fixons dans la suite $P_0\in\P^G(M_0)$ un sous-groupe parabolique minimal de $G$. Soit $M_0(F)^+\eqdef\{m\in M_0\mid \alpha(H_{M_0}(m))\geq 0\,\,\forall \alpha\in \Delta_{P_0}\},$ avec $\Delta_{P_0}$ le sous-ensemble des racines simples dans $G$ relativement à $M_0$ associées à $P_0$. On a la décomposition de Cartan $G(F)=KM_0(F)^+K$. 
On rappelle que $\rho_{P_0}\in a_{0}^\ast$ est la demi-somme des éléments de $\Delta_{P_0}$. Posons 
\[\Xi(g)\eqdef \int_{K}e^{-\rho_{P_0}(H_{M_0}(kg))}\,dk\]
pour tout $g\in G(F)$ (cf. \cite[p.153]{Sil79} pour $F$ non-archimédien ; \cite[p.329]{Vara77} pour $F$ archimédien).

\begin{lemma}[{{\cite[théorème 4.2.1, lemme 4.2.3]{Sil79} pour $F$ non-archimédien ; \cite[théorème II.8.30, proposition II.8.16 (v)]{Vara77} pour $F$ archimédien}}]\label{YDLspeclem:estimateHCfunctions}~{}
\begin{enumerate}
    \item Il existe une constante $c>0$ telle que 
    \[\Xi(l)^{-1}<ce^{\rho_{P_0}(H_{M_0}(l))}\,\,\,\,\forall l\in M_0(F)^+.\] 
    \item Soit $I$ une partie compacte de $G(F)$. Alors il existe des constantes $c_1,c_2>0$ telles que 
    \[c_1\Xi(g)<\sup_{x_1,x_2\in I}\Xi(x_1gx_2)<c_2\Xi(g)\,\,\,\,\forall g\in G(F).\]
\end{enumerate}    
\end{lemma}

Munissons le $F$-espace vectoriel $\g(F)$ d'une norme $\|-\|_{\g}$. Définissons $\|g\|_G\eqdef \max\{\|g\|_\g,\|g^{-1}\|_\g\}$ pour tout $g\in G(F)$. Enfin, on pose $\sigma(g)\eqdef \log \|g\|_G$ pour tout $g\in G(F)$ (cf. \cite[p.150]{Sil79} pour $F$ non-archimédien ;  \cite[p.320]{Vara77} pour $F$ archimédien). Deux fonctions sur $G(F)$ à valeurs dans $\R_{\geq 1}$ sont dites équivalentes si leur quotient est à valeurs dans un compact de $\R_{>0}$. Quand on change la norme vectorielle $\|-\|_{\g}$, la classe d'équivalence de $1+\sigma$ reste la même.

\begin{definition}[Espace de Schwartz-Harish-Chandra $\mathcal{HC}(G(F))$]\label{YDLspecdef:SCHspace}~{}
\begin{enumerate}
    \item (\cite[p.174]{Sil79}) Supposons que $F$ est non-archimédien. Pour tout sous-groupe ouvert-compact $J$ de $G(F)$ on note $\mathcal{HC}(G(F)\ssla J)$ l'espace des fonctions $f$ sur $G(F)$ bi-invariantes par $J$ telles que $\|f\|_{r}\eqdef\sup_{g\in G(F)}\frac{(1+\sigma(g))^r|f(g)|}{\Xi(g)}<+\infty$ pour tout $r\in \R_{\geq 0}$. C'est un espace de Fréchet muni des semi-normes $(\|-\|_{r})_{r\in \R_{\geq 0}}$. Enfin on note $\mathcal{HC}(G(F))$  la limite inductive des $\mathcal{HC}(G(F)\ssla J)$  avec $J$ parcourt l’ensemble des sous-groupes ouverts compacts de $G(F)$, il est muni de la topologie limite inductive. 
    \item (\cite[pp.341-342]{Vara77}) Supposons que $F$ est archimédien. On rappelle que l'algèbre enveloppante $\mathcal{U}(\g(F)_\C)$ agit à gauche et à droite sur l'espace des fonctions lisses sur $G(F)$ par dérivation. On note $\mathcal{HC}(G(F))$ l'espace des fonctions $f$ lisses sur $G(F)$ telles que $\|f\|_{A,B,r}\eqdef\sup_{g\in G(F)}\frac{(1+\sigma(g))^r|AfB(g)|}{\Xi(g)}<+\infty$ pour tout $A,B\in \mathcal{U}(\g(F)_\C)$ et $r\in \R_{\geq 0}$. C'est un espace de Fréchet muni des semi-normes $(\|-\|_{A,B,r})_{A,B\in \mathcal{U}(\g(F)_\C),r\in \R_{\geq 0}}$.
\end{enumerate}    
L'espace $\mathcal{HC}(G(F))$ est appelé l'espace de Schwartz-Harish-Chandra.
\end{definition}

On a $C_c^\infty(G(F))\subseteq \mathcal{HC}(G(F))$. 
On dit qu'une fonctionnelle sur $G(F)$ est une distribution tempérée si elle est continue sur $\mathcal{HC}(G(F))$.



\begin{proposition}[{{\cite[lemme 4.5.2]{Sil79} pour $F$ non-archimédien ; \cite[théorème 12.23]{Knap86} pour $F$ archimédien}}]\label{prop:temprepHSC}
Soit $\pi$ une représentation irréductible de $G(F)$. Alors $\tr(\pi(-))$ est une distribution tempérée si et seulement si $\pi$ est  une représentation tempérée si et seulement si tout coefficient matriciel de $\pi$ est de classe Schwartz-Harish-Chandra. Le cas échéant, la distribution $\tr(\pi(-))$ sur $\mathcal{HC}(G(F))$ est représentée par le caractère de Harish-Chandra $\Theta_\pi$, c’est une fonction de classe Schwartz-Harish-Chandra invariante par $G(F)$-conjugaison.    
\end{proposition}

\begin{proposition}\label{prop:GJinSHC}~{}
\begin{enumerate}
    \item Il y a
    \begin{equation}\label{eq:GJinHSC}
    \bigcup_{s\in \D_\temp(G)}\cZ(\g(F),s)|_{G(F)}\subseteq \mathcal{HC}(G(F)).    
    \end{equation}
    \item L'inclusion \eqref{eq:GJinHSC} est continue.
\end{enumerate}   
\end{proposition}
\begin{proof}~{}
\begin{enumerate}
    \item Soient $f\in \S(\g(F))$ et $s\in \D_\temp(G)$. 

    Supposons d'abord que $F$ est non-archimédien. Soit $r\in \R_{\geq 0}$. On veut prouver que $\|f_s\|_r<+\infty$.

    On définit $\|X\|_{\g}=\max_{1\leq i,j\leq m}|\nu_{\gl_{1,D}}(X_{ij})|$ pour $X=(X_{ij})_{1\leq i,j\leq m}\in \g(F)$. C'est une norme vectorielle bi-invariante par $K$. 

    Le support de $f$ est contenu dans un ouvert compact dans $\g(F)$ bi-$K$-invariant. 
    Il suffit donc de prouver que la fonction
    \[g\mapsto(1+\sigma(g))^r|\nu_\g(g)|^{s+\frac{n-1}{2}}\Xi(g)^{-1}\]
    est bornée sur toute intersection d'un ouvert compact de $\g(F)$ avec $G(F)$. Puis d'après le point 2 du lemme \ref{YDLspeclem:estimateHCfunctions} et la décomposition de Cartan, il suffit de prouver que la fonction
    \begin{equation}\label{eq:interstepGJinSHC}
    l\mapsto(1+\sigma(l))^r|\nu_\g(l)|^{s+\frac{n-1}{2}}\Xi(l)^{-1}    
    \end{equation}
    est bornée sur toute intersection d'un ouvert compact de $\g(F)$ avec $M_0(F)^+$. Ou encore, d'après le point 1 du lemme \ref{YDLspeclem:estimateHCfunctions}, il suffit de prouver que la fonction
    \begin{equation}\label{eq:interstepGJinSHCpadic}
    l\mapsto(1+\sigma(l))^r|\nu_\g(l)|^{s+\frac{n-1}{2}}e^{\rho_{P_0}(H_{M_0}(l))}   
    \end{equation}
    est bornée sur toute intersection d'un ouvert compact de $\g(F)$ avec $M_0(F)^+$.

    Soit 
    \[l=\begin{bmatrix}
     a_{1} &  &  & \\ 
     & a_{2} &  & \\ 
     &  &  \ddots & \\ 
     &  &   & a_{m} 
    \end{bmatrix}\in M_0(F).\]
    Nous avons 
    \[1+\sigma(l)=1+\log\max\left\{\max_{1\leq i\leq m}|\nu_{\gl_{1,D}}(a_{i})|,\max_{1\leq i\leq m}|\nu_{\gl_{1,D}}(a_{i})|^{-1}\right\},\]
    et
    \[\nu_\g(l)=\prod_{1\leq i\leq m}|\nu_{\gl_{1,D}}(a_{i})|,\]
    enfin
    \[e^{\rho_{P_0}(H_{M_0}(l))}=\prod_{1\leq i\leq m}|\nu_{\gl_{1,D}}(a_{i})|^{(\frac{m+1}{2}-i)d}.\] 
    On est alors conduit à montrer que la fonction qui à $(a_1,\dots,a_m)$ associe
    \begin{equation}\label{eq:GJinHSCpadicestimate}
    \left(1+\log\max\left\{\max_{1\leq i\leq m}|\nu_{\gl_{1,D}}(a_{i})|,\max_{1\leq i\leq m}|\nu_{\gl_{1,D}}(a_{i})|^{-1}\right\}\right)^r\prod_{1\leq i\leq m}|\nu_{\gl_{1,D}}(a_{i})|^{s+\frac{d-1}{2}+(n-id)}.    
    \end{equation}
    est bornée sur toute intersection d'un ouvert compact de $\gl_1(D)^m$ avec $\GL_1(D)^m$. Cela est vrai.

    Supposons ensuite que $F$ est archimédien. Soient $A,B\in \mathcal{U}(\g(F)_\C)$ et $r\in \R_{\geq 0}$. On veut prouver que $\|f_s\|_{A,B,r}<+\infty$. Par linéarité on peut et on va supposer que $A$ et $B$ sont des monômes. 
    
    Soit $C\in \mathcal{U}(\mathfrak{g}(F)_\C)$ un élément de degré 1, on a 
    \begin{equation}\label{eq:GJinHSCarchiestimate1}
    \begin{cases}
     &C\ast f_s=(C\ast f)_s+f_s\left(s+\frac{n-1}{2}\right)\tau_{\g,\C}(C) \\  
     &f_s\ast C=(f\ast C)_s+f_s\left(s+\frac{n-1}{2}\right)\tau_{\g,\C}(C),
    \end{cases}
    \end{equation}
    où $\tau_\g:\g(F)\to F$ est la trace réduite et $\tau_{\g,\C}: U\otimes t\in \g(F)\otimes_{\R}\C\mapsto t\tau_{\g}(U)\in \C$ sa complexification. Il est clair que toutes les dérivées d'une fonction de classe Schwartz-Bruhat restent de classe Schwartz-Bruhat. Ainsi, dans la suite, pour simplifier la démonstration, on peut et on va supposer que $A$ et $B$ sont 1.

    On a $D=\R$, $\C$, ou $\mathbb{H}$, l'algèbre des quaternions de Hamilton. On note $Y\in D\mapsto \overline{Y}\in D$ l'involution principale sur $D$. On note $P_0$ la fonction polynômiale à $[F:\R]^{n^2}$ variables réelles $P_0(X)=\sum_{1\leq i,j\leq m}X_{ij}\overline{X_{ij}}$ pour $X=(X_{ij})_{1\leq i,j\leq m}\in \g(F)$. On définit $\|X\|_{\g}=\sqrt{P_0(X)}$ sur $\g(F)$. C'est une norme vectorielle. Clairement il existe des constantes strictement positives $c_1$ et $c_2$ dépendant uniquement de $G$ tel que
    \[c_1P_0(X)<P_0(k_1Xk_2)<c_2P_0(X),\,\,\,\,\forall k_1,k_2\in K, X\in \g(F).\]

    Soit $B_1$ l'intersection de $G(F)$ avec la boule fermée centrée en 0 de rayon 1 dans $\g(F)$ pour la norme $\|-\|_\g$, et $B_1^c$ son complémentaire dans $G(F)$. Le lemme \ref{YDLspeclem:estimateHCfunctions} nous assure, afin de prouver la majoration $\|f_s\|_{A=1,B=1,r}<+\infty$, qu'il suffit de prouver que 
    \begin{enumerate}
        \item la fonction
    \begin{equation}\label{eq:GJinHSCarchiestimate2}
    l\mapsto(1+\sigma(l))^r|\nu_\g(l)|^{s+\frac{n-1}{2}}\|f\|_{P=1,Q=1}e^{\rho_{P_0}(H_{M_0}(l))}   
    \end{equation}
    est bornée sur $M_0(F)\cap B_1$, pour rappel $\|-\|_{P,Q}$ est la semi-norme \eqref{YDLspeceq:defSchseminormFarchi} sur $\S(\g(F))$;
        \item la fonction
    \begin{equation}\label{eq:GJinHSCarchiestimate3}
    l\mapsto(1+\sigma(l))^r|\nu_\g(l)|^{s+\frac{n-1}{2}}P_0(l)^{-(s+\frac{2n-1-d}{2})-1}\|f\|_{P=P_0,Q=1}^{(s+\frac{2n-1-d}{2})+1}e^{\rho_{P_0}(H_{M_0}(l))}  
    \end{equation}
    est bornée sur $M_0(F)\cap B_1^c$. 
    \end{enumerate}

    Soit 
    \[l=\begin{bmatrix}
     a_{1} &  &  & \\ 
     & a_{2} &  & \\ 
     &  &  \ddots & \\ 
     &  &   & a_{m} 
    \end{bmatrix}\in M_0(F).\]
    Nous avons 
    \[1+\sigma(l)=1+\frac{1}{2}\log\max\left\{\sum_{1\leq i\leq m}a_i\overline{a_i},\sum_{1\leq i\leq m}(a_i\overline{a_i})^{-1}\right\},\]
    et
    \[\nu_\g(l)=\prod_{1\leq i\leq m}|\nu_{\gl_{1,D}}(a_{i})|,\]
    enfin
    \[e^{\rho_{P_0}(H_{M_0}(l))}=\prod_{1\leq i\leq m}|\nu_{\gl_{1,D}}(a_{i})|^{(\frac{m+1}{2}-i)d}.\]

    Primo, dans le cas (a), les deux quantités $(1+\sigma(l))^r$, $|\nu_\g(l)|^{s+\frac{n-1}{2}}e^{\rho_{P_0}(H_{M_0}(l))}$ sont bornées par des constantes. Secundo, dans le cas (b), $(1+\sigma(l))^r|\nu_\g(l)|^{s+\frac{n-1}{2}}P_0(l)^{-(s+\frac{2n-1-d}{2})-1}e^{\rho_{P_0}(H_{M_0}(l))}$ est bornée par une constante. En conclusion on a la majoration $\|f_s\|_{A=1,B=1,r}<+\infty$, ce qu'il fallait.

    \item Supposons d'abord que $F$ est non-archimédien. Rappelons que la topologie sur $\S(\g(F))$ est la topologie discrète. On déduit la continuité de l'inclusion \eqref{eq:GJinHSC} de l'équation \eqref{eq:GJinHSCpadicestimate}. Supposons ensuite que $F$ est archimédien. La topologie de $\D_\temp(G)$ est induite par la valeur absolue $|-|$ sur $\C$. D'après la théorie générale, l'inclusion \eqref{eq:GJinHSC} est continue si et seulement si pour tout $A,B\in \mathcal{U}(\g(F)_\C)$, $r\in \R_{\geq 0}$, $f_0\in \S(\g(F))$ et $s_0\in \D_\temp(G)$, on peut majorer $\|f_s-f_{0,s_0}\|_{A,B,r}$ par une combinaison linéaire positive d'un nombre fini de semi-normes $|s-s_0|+\|f-f_0\|_{P,Q}$ pour tout $s$ dans un voisinage de $s_0$ et $f$ dans un voisinage de $f_0$. 
    Or $f_s-f_{0,s_0}=(f-f_0)_{s_0}+f_{s_0}|\nu_\g|^{s-s_0}$. La majoration voulue est de ce fait une conséquence des équations \eqref{eq:GJinHSCarchiestimate1}, \eqref{eq:GJinHSCarchiestimate2} et \eqref{eq:GJinHSCarchiestimate3}. \qedhere
\end{enumerate}
\end{proof}

Soient $L\in \L^G$, $P\in \P^G(L)$, $Q\in \F^G(L)$, $\lambda\in a_{L,\C}^\ast$ et $\pi\in \Pi_\temp(L(F))$. On déduit de l'équation \eqref{eq:opN} le poids $\mathcal{N}_L^Q(P,\pi,\lambda)$, et cette fonction méromorphe de $\lambda$ est analytique sur $ia_L^\ast$. On écrit $\mathcal{N}_L^Q(P,\pi)=\mathcal{N}_L^Q(P,\pi,0)$. Introduisons le caractère pondéré local, pour $h$ une fonction définie sur $G(F)$,
\begin{equation}\label{eq:defcarapondlocaltemp}
J_L^Q(\pi,h)\eqdef\tr(\mathcal{N}_L^Q(P,\pi)\I_P^G(\pi,h)),    
\end{equation}
cette valeur ne dépend pas de $P\in \P^G(L)$.

\begin{proposition}\label{prop:caraunitponddecomp}
Soient $L\in \L^G$, $Q\in \F^G(L)$, et $\pi\in \Pi_\temp(L(F))$. Le caractère pondéré local $J_L^Q(\pi,-)$ est une distribution sur $\mathcal{HC}(G(F))$.     
\end{proposition}
\begin{proof}
Cela vient de la proposition \ref{prop:temprepHSC} et du fait que $\mathcal{N}_L^Q(P,\pi)$ est unitaire.
\end{proof}


\begin{lemma}\label{YDLspeclem:holomorphieint}
Soit $X$ un espace mesurable (dont on néglige la tribu et la mesure dans la notation), et $\mathcal{D}\subseteq\C$ un domaine de $\C$. Soit $h:\mathcal{D}\times X\rightarrow\C$ une fonction assujettie aux conditions suivantes.
\begin{enumerate}
    \item Pour tout $s\in \mathcal{D}$, $h(s,\cdot)$ est mesurable.
    \item Pour tout $x\in X$, $h(\cdot,x)$ est holomorphe.
    \item La fonction $\int_{X}|h(\cdot,x)|\,dx$ est localement bornée.
\end{enumerate}
Alors la fonction $\int_{X}h(\cdot,x)\,dx$ est holomorphe sur $\mathcal{D}$ et toutes ses dérivées s'obtiennent par dérivation sous le signe intégrale.
\end{lemma}
\begin{proof}
Le lemme est connu et on omet la preuve.    
\end{proof}

\begin{proposition}\label{prop:holotempcarainDt}
Soient $L\in \L^G$, $Q\in \F^G(L)$, $\pi\in \Pi_\temp(L(F))$.
\begin{enumerate}
    \item Supposons que $F$ est non-archimédien. Soit $f\in\S(\g(F))$. Alors l'application $s\mapsto J_{L}^Q(\pi,f_s)$ définit une fonction holomorphe sur le domaine $s \in \D_\temp(G)$.
    \item Supposons que $F$ est archimédien. Soit $f\in\S(\g(F),\psi)$. Alors l'application $s\mapsto J_{L}^Q(\pi,f_s)$ définit une fonction holomorphe sur le domaine $s \in \D_\temp(G)$.
\end{enumerate}    
\end{proposition}

\begin{proof}Nous allons réduire la situation au lemme précédent. Supposons d'abord que  $F$ est non-archimédien. Soit $J$ un sous-groupe ouvert compact de $G(F)$ tel que $f$ soit bi-invariante par $J$. On pose $(e_{i})_{i\in I}$ une base orthonormée de $\text{Im}(\I_P^G(\pi)(f))\subseteq \I_P^G(\pi)^{J}$, avec $I$ un ensemble fini. Supposons ensuite que $F$ est archimédien. Comme tout élément de $\S(\g(F),\psi)$ est $K\times K$-fini, on en déduit l'existence de $\xi:K\to \C$ qui est une somme finie de fonctions $k\in K\mapsto (\dim \sigma_{j})\tr (\sigma_{j}(k))$ où $\sigma_{j}$ est une représentation irréductible de $K$, telle que $\xi\ast f \ast \xi=f$. On pose $(e_{i})_{i\in I}$ une base orthonormée de $\text{Im}(\I_P^G(\pi)(\xi))$, avec $I$ un ensemble fini.

Nous avons 
\begin{align*}
\tr\left(\mathcal{N}_{L}^{Q}(P,\pi)\I_{P}^{G}(\pi, f_{s}))\right)&=\sum_{i\in I}\int_{G(F)}f(g)\langle\mathcal{N}_{L}^{Q}(P,\pi)\I_{P}^{G}(\pi)(g)e_{i},e_{i} \rangle |\nu_\g(g)|_v^{s+\frac{n-1}{2}}\,dg    \\
&=\int_{G(F)}\sum_{i\in I}f(g)\langle\mathcal{N}_{L}^{Q}(P,\pi)\I_{P}^{G}(\pi)(g)e_{i},e_{i} \rangle |\nu_\g(g)|_v^{s+\frac{n-1}{2}}\,dg.    
\end{align*}
L'interversion série-intégrale se justifie par l'inégalité de Cauchy-Schwartz et par le fait que $\mathcal{N}_{L}^{Q}(P,\pi)$ est unitaire et $I$ est fini. 
On prend $\mathcal{D}=\D_{\temp}(G)$, $X=G(F)\times I$, et $h(s,g,i)=f(g)\langle\mathcal{N}_{L}^{Q}(P,\pi)\I_{P}^{G}(\pi)(g)e_{i},e_{i} \rangle |\nu_\g(g)|_v^{s+\frac{n-1}{2}}$ dans le lemme plus haut. Les conditions 1 et 2 sont clairement vérifiées. Puis la condition 3 provient du point 2 de la proposition \ref{prop:GJinSHC} et  de la proposition \ref{prop:caraunitponddecomp}.
\end{proof}

\subsubsection{Caractère pondéré unitaire local}

Soit dans cette partie $F$ local. Soient $L\in \L^G$, $P\in \F^G(L)$, $Q\in \F^G(L)$, et $\pi\in \Pi_\unit(L(F))$. On définit le caractère pondéré unitaire local $J_L^Q(\pi,-)$ par la même équation \eqref{eq:defcarapondlocaltemp} (cf. \cite[section 2.7]{Walds12}). C'est une distribution, autrement dit une fonctionnelle continue sur $C_c^\infty(G(F))$.

\begin{proposition}\label{prop:holounitcarainDu} Soient $L\in \L^G$, $Q\in \F^G(L)$, $\pi\in \Pi_\unit(L(F))$.
\begin{enumerate}
    \item Supposons que $F$ est non-archimédien. Soit $f\in\S(\g(F))$. Alors l'application $s\mapsto J_{L}^Q(\pi,f_s)$ définit une fonction holomorphe sur le domain $s \in \D_\unit(G)$.
    \item Supposons que $F$ est archimédien. Soit $f\in\S(\g(F),\psi)$. Alors l'application $s\mapsto J_{L}^Q(\pi,f_s)$ définit une fonction holomorphe sur le domain $s \in \D_\unit(G)$.
\end{enumerate}    
\end{proposition}

\begin{proof}
C'est presque le même argument que celui de la proposition \ref{prop:holotempcarainDt}. On définit comme dans la preuve de la proposition \ref{prop:holotempcarainDt} un ensemble fini $I$. Nous avons 
\begin{align*}
\tr\left(\mathcal{N}_{L}^{Q}(P,\pi)\I_{P}^{G}(\pi, f_{s}))\right)&=\sum_{i\in I}\int_{G(F)}f(g)\langle\mathcal{N}_{L}^{Q}(P,\pi)\I_{P}^{G}(\pi)(g)e_{i},e_{i} \rangle |\nu_\g(g)|_v^{s+\frac{n-1}{2}}\,dg    \\
&=\int_{G(F)}\sum_{i\in I}f(g)\langle\mathcal{N}_{L}^{Q}(P,\pi)\I_{P}^{G}(\pi)(g)e_{i},e_{i} \rangle |\nu_\g(g)|_v^{s+\frac{n-1}{2}}\,dg.    
\end{align*}
L'interversion série-intégrale se justifie par l'inégalité de Cauchy-Schwartz et par le fait que $\mathcal{N}_{L}^{Q}(P,\pi)$ est unitaire et $I$ est fini. 
On prend $\mathcal{D}=\D_{\unit}(G)$, $X=G(F)\times I$, et $h(s,g,i)=f(g)\langle\mathcal{N}_{L}^{Q}(P,\pi)\I_{P}^{G}(\pi)(g)e_{i},e_{i} \rangle |\nu_\g(g)|_v^{s+\frac{n-1}{2}}$ dans le lemme plus haut. Les conditions 1 et 2 sont clairement vérifiées. Puis la condition 3 provient du fait que
\[\int_{G(F)}\sum_{i\in I}|f(g)\langle\mathcal{N}_{L}^{Q}(P,\pi)\I_{P}^{G}(\pi)(g)e_{i},e_{i} \rangle ||\nu_\g(g)|_v^{s+\frac{n-1}{2}}\,dg\leq |I|\int_{G(F)}|f(g)||\nu_\g(g)|^{s+\frac{n-1}{2}}\,dg\]
converge et définit une fonction holomorphe sur le domaine $s\in \D_\unit(G)$ selon la preuve de \cite[proposition]{GoJa72} (la preuve ne demande que la finitude des semi-normes $(\|-\|_{P,Q=1})_P$ de l'équation \eqref{YDLspeceq:defSchseminormFarchi}, donc elle fonctionne bien que $|f|$ ne soit pas de classe Schwartz-Bruhat). 
\end{proof}

\subsubsection{Formule des traces non-invariante}\label{subsubsec:FTnon-inv}

On impose les normalisations des mesures consignées dans le numéro \ref{subsec:normalisationmesureofspec1}.

\begin{proposition}\label{prop:espacedeGJsubseteqMuFiLa}~{}
\begin{enumerate}
    \item Il y a 
    \begin{equation}\label{YDLspeceq:GJinCMFL}
    \bigcup_{s\in \D_\adel(G)}\cZ(\g(\A_F),s)|_{G(\A_F)}\subseteq \CMFL.    
    \end{equation}
    \item L'inclusion \eqref{YDLspeceq:GJinCMFL} est continue.
\end{enumerate}  
\end{proposition}



\begin{proof}~
\begin{enumerate}
    \item Soit $f_s\in\cZ(\g(\A_F),s)$. Il est clair que toutes les dérivées d'une fonction de classe Schwartz-Bruhat restent de classe Schwartz-Bruhat. Soit $C\in \mathcal{U}(\mathfrak{g}(F_\infty)_\C)$ un élément de degré 1, on a 
\begin{equation}\label{YDLspeceq:GJinCMLFinterstep1}
f_s\ast C=(f\ast C)_s+f_s\times \left(s+\frac{n-1}{2}\right)\left(\sum_{v\in\V_{\infty}}\tau_{\g_v,\C}(C_v)\right),    
\end{equation}
où $\tau_{\g_v}:\g(F_v)\to F_v$ est la trace réduite et $\tau_{\g_v,\C}: U\otimes t\in \g(F_v)\otimes_{\R}\C\mapsto t\tau_{\g}(U)\in \C$ sa complexification.
Il suffit alors de montrer que $f_s\in L^1(G(\A_F))$. La proposition s'ensuit maintenant du fait que 
\[\int_{G(\A_F)} |f(g)||\nu_\g(g)|^{s+\frac{n-1}{2}}dg\]
est convergente dans le domaine $s\in\D_\adel(G)$ 
: en effet la mesure $dg$ sur $G(\A_F)$ est donnée par $d^+g/|\nu_{\g}(g)|^n$ avec $d^+g$ la mesure sur $\g(\A_F)$, et pour presque tout $v$ on a $\int_{G(F_v)} |1_{\mathfrak{k}_v}(g)||\nu_\g(g)|_v^{s+\frac{n-1}{2}}dg=\prod_{i=0}
^{n-1}L_v(1_{\GL_{1}(F_v)},s-i+\frac{n-1}{2})$, le produit est donc convergente dans le domaine indiqué (point 1 du théorème \ref{thm:GJzetalocalarchi}).
    
    \item Il nous suffit de prouver que l'inclusion \begin{equation}\label{YDLspeceq:GJinCMLFinterstep2}
    \bigcup_{s\in \D_\adel(G)}\cZ(\g(\A_F)\ssla J,s)|_{G(\A_F)}\subseteq \mathcal{C}(G(\A_F),J).    
    \end{equation}
    est continue, pour tout $J$ sous-groupe ouvert compact de $G(\A_{F,\fin})$. Soit $s\in \D(G)$. Il est évident que $\cZ(\g(\A_F),s)=\cZ(\g(\A_{F,\fin}),s)\widehat{\otimes}\cZ(\g(F_{\infty}),s)$, où  $\otimes$ est le produit tensoriel d'espaces vectoriels topologiques, muni de sa topologie usuelle (cf. \cite[définition 43.1 ou 43.2]{BookTreve}, les espaces en questions étant nucléaires, la $\epsilon$-topologie et la $\pi$-topologie sur le produit tensoriel sont donc la même), puis $\widehat{\otimes}$ le complété pour cette topologie. Or $\cZ(\g(\A_{F,\fin}),s)$ est muni de la topologie discrète. Il vient que l'inclusion \eqref{YDLspeceq:GJinCMLFinterstep2} est continue si et seulement si pour tout $X\in \mathcal{U}({\g(F_\infty)_\C})$, $f_{\fin}\in \S(\g(\A_{F,\fin}))$, $f_0\in \S(\g(F_\infty))$ et $s\in \D_\adel(G)$, on peut majorer $\|(f_\fin\otimes f)_s-(f_\fin\otimes f_0)_{s_0}\|_{X}$ ($\|-\|_X$ est la semi-norme \eqref{YDLspeceq:defCMFLseminorm}) par une combinaison linéaire positive d'un nombre fini de semi-normes $|s-s_0|+\|f-f_{0}\|_{P,Q}$ pour tout $s$ dans un voisinage de $s_0$ et $f$ dans un voisinage de $f_0$. Or $(f_\fin\otimes f)_s-(f_\fin\otimes f_0)_{s_0}=f_\fin\otimes (f-f_0)_{s_0}+f_\fin\otimes f_{s_0}|\nu_\g|_{\infty}^{s-s_0}$. 
    La majoration voulue est de ce fait une conséquence de l'équation \eqref{YDLspeceq:GJinCMLFinterstep1}.\qedhere
\end{enumerate}
\end{proof}

On notera $\S(\g(\A_F),\psi)$ le sous-espace de $\S(\g(\A_F))$ engendré par les tenseurs purs $f=\bigotimes_{v\in \V_F}f_v$ avec $f_v\in \S(\g(F_v),\psi_v)$ pour tout $v\in\V_{\infty}$.

\begin{proposition}Soit $f\in \S(\g(\A_F),\psi)$. L'équation \eqref{eq:devlopfinabCVcheck}, i.e. (on change les symboles $L$ en $M$, $S$ en $L$, $\pi_1$ en $\pi$ et $Q$ en $P$)
\begin{equation*}
\begin{split}
J_\spec^G(f_s)=\sum_{M\in\L^G}&\sum_{\pi\in\Pi_{\disc}(M)}\sum_{L\in\L^G(M)}\int_{i(a_{M}^G)^\ast}\\
&|W_0^M||W_0^G|^{-1}a_\disc^M(\pi_1)r_M^L(\pi_{\lambda})\tr\left(\mathcal{N}_{L}(P,\pi_{\lambda})\I_{P}^G(\pi_{\lambda},(f_s)^1)\right)\,d\lambda 
\end{split}    
\end{equation*}    
avec $P\in\P^G(M)$, définit une fonction holomorphe sur le domaine $s\in\D_{\adel}(G)$, toutes ses dérivées s'obtiennent par dérivation sous les signes sommes et intégrales.
\end{proposition}
\begin{proof}
Nous allons réduire la situation au lemme \ref{YDLspeclem:holomorphieint}. On reprend les notations de la démonstration du théorème \ref{thm:devlopfinspec}. Afin d'appliquer le lemme, il suffit de prouver que 
\begin{enumerate}
    \item chaque terme $\tr\left(\mathcal{N}_{L}(P,\pi_{\lambda})\I_P^G(\pi_{\lambda},(f_s)^1)\right)$ est holomorphe sur le demi-plan $\Re(s)>\frac{n+1}{2}$ ; et
    \item $|f_s\ast\Delta^{2k}|_{L^1(G(\A_F))}$ est localement borné dans le demi-plan $\Re(s)>\frac{n+1}{2}$, pour $k$ un entier positif fixé assez grand. Ici on utilise implicitement la dernière ligne de l'équation \eqref{eq:majorationderiveJspec}.
\end{enumerate} 

L'assertion 2 est limpide. Quant à la première, on fixe $M\in \L^G$, $\pi\in\prod_{\disc}(M)$, $L\in\L^G(M)$, $\lambda\in i(a_M^G)^\ast$, $P\in\P^G(M)$. On écrit $M\simeq \prod_{i=1}^r\GL_{m_i,D}$, et $\pi=\pi_1\otimes\cdots\otimes\pi_r$ pour $\pi_i$ une représentation irréducible unitaire de $\GL_{m_i,D}(\A_F)^1$. Puis on identifie $a_{M,\C}^\ast=X^\ast(M)\otimes_\Z\C$ à $\C^r$ en envoyant $\nu_{\gl_{m_i,D}}\otimes 1$ sur le vecteur $(0,\dots,0,1,0,\dots,0)$ où $1$ est la $i$-ème coordonnée, et on écrit $\lambda=(\lambda_1,\dots,\lambda_r)\in\C^r$. On prend ensuite $S$ un sous-ensemble fini de $\V_F$, contenant les places archimédiennes, tel que $f^S=1_{\mathfrak{k}^S}$. Ainsi l'opérateur $\I_{P}^G(\pi_{\lambda},(f_s)^1)$ est non-nul uniquement si $\pi^{S}$ est non-ramifié. Le cas échéant d'après la construction d'Arthur (cf. \cite[théorème 2.1 (R8)]{Art89I}) et la formule de scindage pour une $(G,L)$-famille il est clair que, pour tout $\lambda\in i(a_M^G)^\ast$, 
\begin{align*}
\tr(\mathcal{N}_L&(P,\pi_\lambda)\I_P^G(\pi_\lambda,(f_s)^1) )\\
&=\left(\sum_{(A_v)_{v\in S}\in \L^{G_{S}}(L_{S})}d_L^G\left((A_v)_{v\in S}\right)\prod_{v\in S} \tr\left(\mathcal{N}_{L_v}^{Q_{A_v}}(P_v,\pi_{v,\lambda})\I_{P_v}^{G_v}(\pi_{v,\lambda}, f_{s,v}))\right)\right)\\
&\hspace{10cm}\cdot \prod_{i=1}^rL^{S}(\pi_i,s+\lambda_i).
\end{align*}   
L'holomorphie du facteur $\prod_{i=1}^rL^{S}(\pi_i,s+\lambda_i)$ dans le demi-plan $\Re(s)>\frac{n+1}{2}$ résulte de la proposition \ref{prop:Satakeunrcalcul}, et celle de $\tr\left(\mathcal{N}_{L_v}^{Q_{A_v}}(P_v,\pi_{v,\lambda})\I_{P_v}^{G_v}(\pi_{v,\lambda}, f_{s,v}))\right)$ résulte de la proposition \ref{prop:holounitcarainDu}.
\end{proof}

Comme corollaire, toute décomposition du côté spectral compatible à l'équation \eqref{eq:devlopfinabCVcheck} nous gratifie, pour $f\in\S(\g(\A_F),\psi)$, une fonction holomorphe sur le domaine $s\in \D_{\adel}(G)$, dont les dérivées s’obtiennent par dérivation sous les signes sommes et intégrales.

\`{A} titre d'exemple, soit $\mathfrak{X}^G$ l'ensemble des données cuspidales (i.e. les classes de $G(F)$-conjugaison des couples $(L,\sigma)$ consistant en un sous-groupe de Levi semi-standard $L$ et une représentation irréductible cuspidale $\sigma$ de $L(\A_F)^1$). On a la décomposition de Langlands (cf. \cite[III.2.6. théorème]{MW95})
\begin{equation}\label{eq:Langlandsdecomp}
L^2(G(F)\backslash G(\A_F)^1)=\bigoplus_{\chi\in\mathfrak{X}^G}L_{\chi}^2(G(F)\backslash G(\A_F)^1).    
\end{equation}
Cela donne décomposition du côté spectral de la formule des traces $J_{\spec}^G(-)=\sum_{\chi\in\mathfrak{X}^G}J_{\chi}^G(-)$. 

\begin{corollary}\label{coro:CVAholocotespec}
Pour tout $\mathcal{X}\subseteq \mathfrak{X}^G$, la somme
\[\sum_{\chi\in \mathcal{X}}J_{\chi}^G(f_s)\]
définit une fonction holomorphe sur le domaine $s\in \D_\adel(G)$, toutes ses dérivées s'obtiennent par dérivation sous le signe somme.    
\end{corollary}

\subsection{Fonctions méromorphes à plusieurs variables}\label{YDLspecsubsec:fonctionmeropluvar}  

Soit $F$ un corps $p$-adique. On dit qu'un groupe réductif est du type GL (sur $F$) s’il est le groupe des unités d’une algèbre séparable sur $F$ (\cite[définition 2.1]{YDL23b}). Jusqu'à présent, pour des raisons de lisibilité nous avons choisi de travailler avec $G=\GL_{m,D}$ plutôt qu'avec un groupe du type GL général. Pourtant, pour certaines applications, telles que la descente parabolique, il semble inévitable de généraliser les démarches précédentes à tous les groupes du type GL. Nous expliquons alors brièvement ces démarches.

Soit $G$ un groupe du type GL sur $F$. On l'écrit en
\begin{equation}\label{eq:defGLtypegrpstdform}
G=\prod_{\alpha\in A}G_\alpha,\,\,\,\,\,\,\,\,G_\alpha=\Res_{E_\alpha/F}\GL_{m_\alpha,D_\alpha},  
\end{equation}
avec $A$ un ensemble fini, $E_\alpha$ une extension finie de $F$, $D_\alpha$ une algèbre à division de centre $E_\alpha$, et $m_\alpha$ un entier strictement positif. Une représentation irréductible admissible $(\pi,V)$ de $G(F)$ est alors $\pi=\bigotimes_\alpha \pi_\alpha$ avec $\pi_\alpha$ une représentation irréductible unitaire (donc admissible) de $G(F)$ sur un espace de Hilbert $V_\alpha$, et $V=\widehat{\bigotimes}_{\alpha}V_\alpha$. Soit $v=\bigotimes_{\alpha} v_\alpha\in V$ et $v^\vee=\bigotimes_{\alpha}v_\alpha^\vee\in V$. Soit $\Phi=\bigotimes_{\alpha}\Phi_\alpha$ une fonction décomposable avec $\Phi_{\alpha}\in \S(\g_\alpha(F))$. Soit $s=(s_\alpha)_{\alpha\in A}\in \C^A$. On définit 
\[Z(v,v^\vee,\Phi,s)\eqdef \prod_{\alpha\in A}Z(v_\alpha,v_\alpha^\vee,\Phi_\alpha,s_\alpha),\]
c'est une fonction à $|A|=\dim_{\C} a_{G,\C}^\ast$ variables complexes, holomorphe sur un certain domaine de $\C^A$, et elle admet un prolongement méromorphe à $\C^A$. On définit aussi 
\[L(\pi,s)\eqdef \prod_{\alpha\in A}L(\pi_\alpha,s_\alpha),\,\,\,\,\gamma(\pi,\psi,s)\eqdef\prod_{\alpha\in A}\gamma(\pi_\alpha,\psi,s_\alpha),\,\,\,\,\epsilon(\pi,\psi,s)\eqdef \prod_{\alpha\in A}\epsilon(\pi_\alpha,\psi,s_\alpha).\]
Enfin on a $e(G)=\prod_{\alpha}e(G_\alpha)=\prod_{\alpha} (-1)^{m_\alpha(d_\alpha-1)}$ pour le signe de Kottwitz. De même on définit les objets globaux du numéro \ref{subsubsec:globalGJtheory} par produit. Les résultats des sous-sections \ref{subsec:reviewGJtheory} et \ref{subsec:rappelSatake} s'étendent à $G$ de façon évidente.

On définit ensuite les domaines $\D(G)$, $\D_\temp(G)$, $\D_\unit(G)$, et $\D_\adel(G)$ :
\begin{definition}[Domaines de convergence $\D(G)$, $\D_\temp(G)$, $\D_\unit(G)$, et $\D_\adel(G)$] ~{}
\begin{enumerate}
    \item Définissons $\D(G)\eqdef\{s=(s_\alpha)_{\alpha\in A} \mid \forall \alpha,\,\,s_\alpha\in \{-\frac{n_\alpha-1}{2}\}\sqcup\{\Re(s_\alpha)>-\frac{n_\alpha-1}{2}\}\}\subseteq \C^A$. 
    \item Soit $F$ local. Définissons $\D_\temp(G)\eqdef\{s=(s_\alpha)_{\alpha\in A} \mid \forall \alpha,\,\,\Re(s_\alpha)>\frac{-d_\alpha+1}{2}\}\subseteq \C^A$ et $\D_\unit(G)\eqdef\{s=(s_\alpha)_{\alpha\in A} \mid \forall \alpha,\,\,\Re(s_\alpha)>\frac{n_\alpha+1}{2}-d_\alpha\}\subseteq \C^A$. 
    \item Soit $F$ global. Définissons $\D_\adel(G)\eqdef\{s=(s_\alpha)_{\alpha\in A} \mid \forall \alpha,\,\,\Re(s_\alpha)>\frac{n_\alpha+1}{2}\}\subseteq \C^A$. 
\end{enumerate}
\end{definition}

Soit $(G_i)_{i \text{ fini}}$ une famille finie de groupes du type GL. On a bien sûr $\D(\prod_i G)=\prod_i \D(G_i)$
, $\D_\temp(\prod_i G_i)=\prod_i \D_\temp(G_i)$, $\D_\unit\prod_i G_i)=\prod_i \D_\unit(G_i)$ et $\D_\adel(\prod_i G_i)=\prod_i \D_\adel(G_i)$. 

On identifie $a_{G,\C}^\ast=\bigoplus_{\alpha\in A}a_{G_\alpha,\C}^\ast$ à $\C^A$ par l'isomorphism d'espaces linéaires complexes
\begin{equation}\label{eq:idenaGCtoCIgeneral}
(\nu_{\g_{\alpha}}|_{G_\alpha}\otimes s_\alpha)_{\alpha\in A}\mapsto(s_\alpha)_{\alpha\in A},      
\end{equation}
et considérer $\D(G)$, $\D_\temp(G)$, $\D_\unit(G)$, et $\D_\adel(G)$ comme domaines de la variété complexe $a_{G,\C}^\ast$.

On définit les suites $n\eqdef (n_\alpha)_{\alpha \in A}$, $m\eqdef (m_{\alpha})_{\alpha\in A}$, et $d\eqdef (d_{\alpha})_{\alpha\in A}$. Soient $a=(a_\alpha)_{\alpha\in A}\in \C^A$, $b=(b_\alpha)_{\alpha\in A}\in \C^A$, et $c\in \C$. On définit $a+b\eqdef(a_\alpha+{b_\alpha})_{\alpha\in \A}\in \C^A$, $ab\eqdef(a_\alpha b_\alpha)_{\alpha\in \A}\in \C^A$, $a+c=c+a\eqdef (c+a_\alpha)_{\alpha\in \A}\in \C^A$ et $ac=ca\eqdef (ca_\alpha)_{\alpha\in \A}\in \C^A$. Soit $s=(s_\alpha)_{\alpha\in A}\in \C^A$. Quand $F$ est local, par abus de notation on note $|\nu_{\g}(X)|_F^s\eqdef \prod_{\alpha\in A}|\nu_{\g_\alpha}(X_\alpha)|_F^{s_\alpha}$ pour $X\in \g(F)$ ; quand $F$ est global, par abus de notation on note $|\nu_{\g}(X)|_{\A_F}^s\eqdef \prod_{\alpha\in A}|\nu_{\g_\alpha}(X_\alpha)|_{\A_F}^{s_\alpha}$ pour $X\in \g(\A_F)$. Avec ces conventions, on peut définir l'espace $\cZ(G(F),s)$ suivant exactement les mêmes lignes de la définition \ref{def:espacedeGJ}. 

On impose les normalisations des mesures consignées dans le numéro \ref{subsec:rappelSatake}. On munit ensuite pour tout $P\in\F^G$ le groupe $\n_P(F)$ de la mesure correspondant à celle de $N_P(F)$ au sens suivant : prenant une filtration par des sous-groupes normaux $N_P = N_{P,0} \supseteq N_{P,1} \supseteq\cdots\supseteq N_{P,r} = \{1\}$ telle que $N_{P,k}/N_{P,k+1}$ est abélien et $[N_{P,0},N_{P,k}] \subseteq N_{P,k+1}$ pour tout $k$, le groupe $N_{P,k}/N_{P,k+1}$ s'identifie canoniquement avec $\n_{P,k}/\n_{P,k+1}$. On exige que les mesures sur $N_P(F)$ et $\n_P(F)$ sont normalisées de sorte qu'il existe des mesures sur les $N_{P,k}(F)$ (resp. $\n_{P,k}(F)$) telles que la mesure de $N_{P,k}(F)/N_{P,k+1}(F)$ coïncide avec celle de $\n_{P,k}(F)/\n_{P,k+1}(F)$.

Soient $M\in \L^G$ et $P\in \P^G(M)$. On écrit $M^\ast=\prod_{\alpha\in A}M_\alpha^\ast$ où $M_\alpha^\ast\in\L^{G_\alpha^\ast}$. Supposons que  $M_\alpha^\ast=\prod_{j}\GL_{n_{\alpha,j},F}$ avec $\sum_j n_{\alpha,j}=n_\alpha$. La descente parabolique $f_s\mapsto (f_{s,P}: M\mapsto \int_K\int_{\n_P(F)}f_s(k^{-1}(M+N)k)\,dN\,dk)$ fournit une application 
\begin{equation}\label{eq:GJspaceparadescent}
\cZ(G(F)\ssla K,s=(s_\alpha)_\alpha)\rightarrow \cZ(M(F)\ssla K_M,(s_{\alpha,j})_{\alpha,j}),\,\,\,\, \text{avec } s_{\alpha,j}\eqdef s_\alpha+\frac{n_\alpha-n_{\alpha,j}}{2}    
\end{equation}
qui rend le diagramme
\[\begin{tikzcd}
	{\H(G(F)\ssla K)\times\cZ(\g(F)\ssla K,(s_\alpha)_\alpha)} && {\cZ(\g(F)\ssla K,(s_\alpha)_\alpha)} \\
	\\
	{\H(M(F)\ssla K_M)\times\cZ(\m(F)\ssla K_M,(s_{\alpha,j})_{\alpha,j})} && {\cZ(\m(F)\ssla K_M,(s_{\alpha,j})_{\alpha,j})}
	\arrow[from=1-1, to=1-3]
	\arrow[from=3-1, to=3-3]
	\arrow["{((-)_P,(-)_P)}"', from=1-1, to=3-1]
	\arrow["{(-)_P}", from=1-3, to=3-3]
\end{tikzcd}\]
commutatif. On souligne la propriété que l'application $(s_\alpha)_\alpha\mapsto (s_{\alpha,j})_{\alpha,j}$ définie par l'équation \eqref{eq:GJspaceparadescent}, définit des applications $\mathbb{D}(G)\to \mathbb{D}(M)$, $\D_\temp(G)\to \D_\temp(M)$ et $\D_\unit(G)\to \D_\unit(M)$. 

Compte tenu de ces considérations, il découle que le diagramme suivant est commutatif, pour $s\in \D_{\unit}(G)$, 
\begin{equation}
\begin{tikzcd}
	{\H(G(F)\ssla K)} && {\Gamma(\mathcal{X}^G,\O_{\mathcal{X}^G})} \\
	\\
	{\cZ(\g(F)\ssla K,s)} && {C^0(\mathcal{X}_{\unit,\text{ur}}^G;
\C)}
	\arrow["{\text{res}}", from=1-3, to=3-3]
	\arrow[hook, from=1-1, to=3-1]
	\arrow["{\mathcal{S}}", from=3-1, to=3-3]
	\arrow["{\mathcal{S}}", from=1-1, to=1-3]
	\arrow["\sim"', from=1-1, to=1-3]
\end{tikzcd}
\end{equation}
où $C^0(\mathcal{X}_{\unit,\text{ur}}^G;\C)$ est l'algèbre de Banach des fonctions continues à valeurs complexes sur $\mathcal{X}_{\unit,\text{ur}}^G$, 
et la structure de $\cZ(G(F)\ssla K,s)$ en tant que $\H(G(F)\ssla K)$-module est donné par la multiplication des fonctions dans $C^0(\mathcal{X}_{\unit,\text{ur}}^G;\C)$. Par abus de notation on continue à noter la flèche de bas par $\S$. Le théorème de Stone-Weierstrass affirme la densité de $\S\H(G(F)\ssla K)$ dans $C^0(\mathcal{X}_{\unit,\text{ur}}^G;\C)$.

Il serait plus judicieux de travailler avec un groupe de la forme \eqref{eq:defGLtypegrpstdform} dans certains contextes spécifiques. Cependant, pour les besoins de notre article, notamment en ce qui concerne les résultats de comparaison (théorèmes \ref{thm:idenendoglobal}, \ref{thm:transfertspecmeaning}, et \ref{thm:transfergeom'specequiv}), les conclusions dans le cas général peuvent être déduites de celles dans le cas simplifié expliqué dans la section \ref{YDLspecsec:convention}, à travers des réductions claires et immédiates telles que la formule de produit pour les $(G,M)$-familles.

\subsection{Speculation au-delà de l'endoscopie}

On aimerait comprendre la contribution du spectre discret lorsque l'on plonge $f_s$, avec $s\in \mathbb{D}_\adel(G)$ et $f=\bigotimes_{v\in \V_F}f_v\in  \S(\g(\A_F),\psi)$ un tenseur pur, dans la formule des traces non-invariante. Soit $J=\prod_{v\in \V_\fin} J_v\subseteq G(\A_{F,\fin})$ un sous-groupe ouvert compact tel que $f$ soit bi-invariante par $J$. Soit $\pi\in\Pi_\disc(G(\A_F)^1)$. Soit $v\in \V_{\fin}$, on pose $(e_{\pi,v,i})_{i\in I_{\pi,f,v}}$ une base orthonormée de $\text{Im}(\pi_v(f_v))\subseteq V_{\pi_v}^{J_v}$ avec $I_{\pi,f,v}$ un ensemble fini. Alors $\langle \pi_v,f_{s,v}\rangle=\sum_{i\in I_{\pi,f,v}}Z_v(e_{\pi,v,i},e_{\pi,v,i},f_v,s)$ ; soit $v\in \V_{\infty}$. Comme tout élément de $\S(g(F_v),\psi_v)$ est $K_v\times K_v$-fini, on en déduit l'existence de $\xi_v:K_v\to \C$ qui est une somme finie de fonctions $k_v\in K_v\mapsto(\dim \sigma_{v,j})\tr (\sigma_{v,j}(k_v))$ où $\sigma_{v,j}$ est une représentation irréductible de $K_v$, telle que $\xi_v\ast f_v \ast \xi_v=f_v$. On pose $(e_{\pi,v,i})_{i\in I_{\pi,f,v}}$ une base orthonormée de $\text{Im}(\pi_v(\xi_v))$ avec $I_{\pi,f,v}$ un ensemble fini. Alors $\langle \pi_v,f_{s,v}\rangle=\sum_{i\in I_
{\pi,f,v}}Z_v(e_{\pi,v,i},e_{\pi,v,i},f_v,s)$. 

En somme, lorsque l'on plonge $f_s\in \cZ(\g(\A_F),\psi,s)$, avec $f=\bigotimes_{v\in \V_F}f_v$ un tenseur pur, dans la formule des traces non-invariante, la contribution du spectre discret vaut
\begin{equation}\label{eq:discpartfstestfunction}
\sum_{\pi\in\Pi_\disc(G(\A_F)^1)}\prod_{v\in \V_F} \left(\sum_{i\in I_{\pi,f,v}}Z_v(e_{\pi,v,i},e_{\pi,v,i},f_v,s)\right).    
\end{equation}
Pour presque tous les $v$, la somme des fonctions zêta locales vaut la fonction $L$ locale $L_v(\pi_v,s)$. En parallèle, la contribution du spectre continu fait intervenir des opérateurs d'entrelacement et les spectres discrets des sous-groupes de Levi. Comme le côté spectral vaut une somme de fonctions de nature arithmétique, on espère pouvoir en tirer des informations arithmétiques intéressantes. Cette idée a été suggérée par Langlands dans son plan d'aller « au-delà de l'endoscopie » dans un cadre plus général (cf. \cite{Lang04,Lang07}, cf. aussi des analyses ultérieures dans \cite{FLN10,Lang13}).

\section{Correspondance de Jacquet-Langlands}\label{sec:corrJL}

\subsection{Correspondance de Jacquet-Langlands locale}\label{subsec:corrJLlocal} Soit dans cette partie $F$ un corps local de caractéristique 0.

\begin{theorem}[{{\cite[VII.2.6. théorème]{Re10} pour $F$ non-archimédien ; \cite[théorèmes 14.76, 14.90]{Knap86} pour $F$ archimédien}}]\label{thm:classtempbysquareint}On fixe un sous-groupe parabolique minimal de $G$ contenant $M_0$.
\begin{enumerate}
    \item Soit $\pi$ une représentation admissible irréductible et tempérée de $G(F)$. Alors il existe un sous-groupe parabolique standard $P$ de $G$ et une représentation admissible irréductible de carré intégrable $\sigma$ de $M_P(F)$ tels que $\pi$ soit sous-quotient de $\I_P^G(\sigma)$.
    \item Si $P_1$ et $P_2$ sont deux sous-groupes paraboliques standard, $\sigma_1$ et $\sigma_2$ sont respectivement représentations admissibles irréductibles de carré intégrable de $M_{P_1}(F)$ et $M_{P_2}(F)$. Alors $\I_{P_1}^G(\sigma_1)$ et $\I_{P_2}^G(\sigma_2)$ admettent un sous-quotient en commun si et seulement s'il existe $w\in W_0^G$ tel que $M_{P_1}=(\Ad w)M_{P_2}$ et $\sigma_1\simeq (\Ad w)\sigma_2$.
\end{enumerate}
\end{theorem}

Soit $\Pi_{2}(G(F))$ (resp. $\Pi_{\temp}(G(F))$ ; $\Pi_{\unit}(G(F))$) l'ensemble des classes d’isomorphismes des représentations irréductibles de carré intégrable (resp. tempérées ; unitaires).

\begin{theorem}[{{\cite[VII.2.1 remarque 2, VII.2.6. corollaire]{Re10} pour $F$ non-archimédien ; \cite[theorems 8.51, 8.53.]{Knap86} pour $F$ archimédien}}]\label{thm:propertiesofrep} On a $\Pi_{2}(G(F))\subseteq \Pi_{\temp}(G(F))\subseteq \Pi_{\unit}(G(F))$.
\end{theorem}

\begin{theorem}\label{thm:irredunitparaind} Soient $\pi_1\in \Pi_{\unit}(\GL_{m_1}(D))$ et $\pi_2\in \Pi_{\unit}(\GL_{m_2}(D))$. Soit $P$ un sous-groupe parabolique de $\GL_{m_1+m_2,D}$ dont le facteur de Levi est isomorphe à $\GL_{m_1,D}\times\GL_{m_2,D}$. Alors $\I_P^{\GL_{m_1+m_2,D}}(\pi_1\otimes \pi_2)\in \Pi_{\unit}(\GL_{m_1+m_2}(D))$.   
\end{theorem}
\begin{remark}\label{rem:R-grpistriv}
Il est facile de voir que le R-groupe d'un paramètre de Langlands de $G(F)$ est toujours trivial. Ce théorème dit que le R-groupe d'une représentation irréductible unitaire de $G(F)$ est également trivial.   
\end{remark}
\begin{proof}
Ce théorème a été établi, pour $D=F$ non-archimédien par Bernstein \cite[0.2. corollaire]{Bern84} ; pour $D$ quelconque sur $F$ un corps local non-archimédien par Sécherre \cite[proposition 3.1, section 3.3]{Sech09} ; pour $D$ quelconque sur $F$ un corps local archimédien par Vogan \cite[théorème 17.6]{Vog86}. Voir \cite[section 7.2]{Badu09} pour l'historique et aussi \cite[section 12]{Badu09} pour plus d'explication sur l'article de Vogan cité. 
\end{proof}

\`{A} titre de rappel, si $\pi$ est une représentation de $G(F)$ alors on note $\Theta_{\pi}$ son caractère de Harish-Chandra. On note $G_{\rss}$ le lieu semi-simple régulier de $G$. On définit la notion de transfert des classes de conjugaisons et des éléments des groupes des $F$-points comme dans le texte précédent (\cite[définitions 2.4]{YDL23b}). On suit aussi les conventions de \cite[sous-section 2.4]{YDL23b} : on note, en ajoutant une étoile $^\ast$, en exposant les objets relativement à $G^\ast$ qui se transfèrent à $G$ ; pour les objets généraux relativement à $G^\ast$ on les note en ajoutant un signe prime $'$ en exposant. Par exemple en écrivant $g^\ast\in G^\ast(F)$ il est sous-entendu qu'il existe $g\in G(F)$ tel que $(\Ad G^\ast(F))g^\ast$ se transfère en $(\Ad G(F))g$. Alors que l'écriture $g'\in G^\ast(F)$ n'impose pas de condition particulière sur $g'$. Le contexte devrait permettre d’enlever toute ambiguïté. Soit $g'\in G^\ast(F)$. Quand $g'$ se transfère on écrit $g'=g^\ast$ et $g'\arr g^\ast$, avec $g\in G(F)$ un transfert de $g^\ast$.

Une représentation $\pi'$ irréductible de $\GL_{md}(F)$ est dite $d$-compatible si $\Theta_{\pi'}(g^\ast)\not=0$ pour un certain $\GL_{m,\rss}(D)\ni g\arr g^\ast\in \GL_{md,\rss}^\ast(F)$. La condition dépend uniquement de l'ordre de $D$ dans le groupe de Brauer $\text{Br}(F)$. \'{E}crivons $\Pi_{\unit}^{d}(\GL_{md}(F))$ l’ensemble des classes d'isomorphismes des représentations irréductibles unitaire $d$-compatible.

On écrit $\mathcal{R}(\GL_k(D))$ pour le groupe de Grothendieck des représentations admissibles de longueur finie de $\GL_k(D)$. Posons $\mathcal{R}(D)=\bigoplus_{k\in \N}\mathcal{R}(\GL_k(D))$ avec pour convention $\GL_0(D)=\mathcal{R}(\GL_0(D))=\Z$, c'est un anneau commutatif avec l'induction parabolique comme loi multiplicative que l'on notera $\times$.


\begin{proposition}[{{\cite[proposition 15.1]{Badu09}}}]\label{prop:d-compcommutewithind}
Soient $\pi_i'\in \Pi(\GL_{n_i}(F))$ pour $i=1,\dots,k$. Alors $\pi_1'\times\cdots\times \pi_k'$ est $d$-compatible si et seulement si $d\mid n_i$ et $\pi_i'$ est $d$-compatible pour tout $i$.
\end{proposition}

\begin{theorem}\label{thm:appLJlocal}
Il existe un morphisme d'anneaux $\bLJ:\mathcal{R}(F)\to\mathcal{R}(D)$ vérifiant 
\begin{enumerate}
    \item $\bLJ=0$ sur $\bigoplus_{k\in \N,d\nmid k}\mathcal{R}(\GL_k(F))$, et $\bLJ(\mathcal{R}(\GL_{md}(F)))\subseteq\mathcal{R}(\GL_{m}(D))$ pour tout $m\in \mathbb{N}$ ;
    \item $\Theta_{\pi'}(g^\ast)=e(\GL_{m,D})\Theta_{\bLJ(\pi')}(g)$ pour tout $m\in \N$, $\GL_{m,\rss}(D)\ni g\arr g^\ast\in \GL_{md,\rss}^\ast(F)$ et $\pi'\in \mathcal{R}(\GL_{md}(F))$ ; 
    \item $\bLJ$ est surjectif ;
    \item $\bLJ$ se restreint à une bijection $\bLJ:\Pi_{2}(\GL_{md}(F))\rightarrow \Pi_{2}(\GL_{m}(D))$ pour tout $m\in \mathbb{N}$;
    \item $\bLJ(\Pi_{\unit}(\GL_{md}(F)))\subseteq \{0\}\sqcup\pm \Pi_{\unit}(\GL_m(D))$ pour tout $m\in \N\setminus\{0\}$, aussi pour $\pi'\in \Pi_{\unit}(\GL_{md}(F))$, $\bLJ(\pi')\not=0$ si et seulement si $\pi'\in \Pi_{\unit}^{d}(\GL_{md}(F))$.
    \item $\bLJ(\pi_{\lambda^\ast}')=\bLJ(\pi')_\lambda$ pout tous $m\in\N$, $\pi'\in \mathcal{R}(\GL_{md}(F))$ et $\lambda^\ast\in a_{\GL_{md,F}}^\ast$, ici $\lambda\in a_{\GL_{m,D}}^\ast$ est l'image réciproque de $\lambda^\ast$ par l'isomorphisme $a_{\GL_{m,D}}^\ast\xrightarrow{\sim}a_{\GL_{md,F}}^\ast$ induit par le torseur intérieur associé $\eta: \GL_{m,D}\to \GL_{md,F}$.
\end{enumerate}
C'est l'unique morphisme d'anneaux vérifiant les points 2, 4 et 6.
\end{theorem}
\begin{proof}
Pour $F$ non-archimédien, on prend $\bLJ$ le morphisme défini par Badulescu en \cite[proposition 3.3]{Badu07}. Voir \cite[théorème 3.1, proposition 3.3]{Badu07} pour les points 1 et 4, puis \cite[proposition 3.3]{Badu07} pour les points 2 et 3, ensuite \cite[proposition 3.9 (c), proposition 3.10 (b)]{Badu08} pour le point 5. Enfin le point 6  découle directement de la définition  \cite[section 3.1]{Badu07}.

Quant à $F$ archimédien, on note $\bLJ'$ le morphism de Badulescu-Renard en \cite[discussion après le théorème 4.1]{Badu09}, ensuite on pose $\bLJ:\mathcal{R}(F)\to\mathcal{R}(D)$ l'unique morphisme d'anneaux vérifiant le point 6 et $\bLJ=e(\GL_{m,D})\bLJ'$ sur $\Pi_{2}(\GL_{md}(F))$, pour tout $m\in \N$ (c'est un morphisme d'anneaux car $e(\GL_{m_1,D})e(\GL_{m_1,D})=e(\GL_{m_1+m_2,D})$). Voir la discussion après \cite[théorème 4.1]{Badu09} pour les points 1, 3 et 4, puis \cite[théorème 13.8]{Badu09} pour le point 5, et encore une fois le point 6 découle directement de la définition, enfin la preuve du point 2 est la même qu'au cas non-archimédien. \qedhere
\end{proof}

Bien sûr $\bLJ$ est l'application identité lorsque $D=F$. Grâce au point 5, on voit que $\bLJ(\Pi_{\unit}^{d}(\GL_{md}(F))\subseteq \pm \Pi_{\unit}(\GL_m(D))$. On note $|\bLJ|:\Pi_{\unit}^{d}(\GL_{md}(F))\rightarrow \Pi_{\unit}(\GL_m(D))$ l'application qui à $\pi'$ associe la vraie représentation parmi $\{\bLJ(\pi'),-\bLJ(\pi')\}$.

\begin{corollary}\label{coro:localJLtemp}Soit $m\in \mathbb{N}\setminus\{0\}$. Alors, par restriction, nous avons $\bLJ:\Pi_{\temp}(\GL_{md}(F))\rightarrow \{0\}\sqcup \Pi_{\temp}(\GL_{m}(D))$, et son image contient $\Pi_{\temp}(\GL_{m}(D))$. De plus, toute fibre sur $\Pi_{\temp}(\GL_{m}(D))$ de cette restriction est un singleton.
\end{corollary}
\begin{proof}
Soit $\pi'\in \Pi_{\temp}(\GL_{md}(F))$. Au moyen du point 1 du théorème \ref{thm:classtempbysquareint} et du théorème \ref{thm:irredunitparaind}, on peut écrire $\pi'=\sigma_1'\times \cdots \times\sigma_l'$ avec $\sigma_i'\in \Pi_2(\GL_{k_i}(F))$ et $\sum_i k_i=md$. On a $\bLJ(\pi')=\prod_i\bLJ(\sigma_i')$. S'il existe $k_i$ qui n'est pas un multiple de $d$ alors $\bLJ(\pi')=0$. Si tous les $k_i$ sont des multiples de $d$ alors $\bLJ(\pi')\in \Pi_{\temp}(\GL_{m}(D))$ en vertu de la proposition précédente. 

Si $\pi\in \Pi_{\temp}(\GL_{m}(D))$, on peut écrire $\pi=\sigma_1\times \cdots \times\sigma_l$ avec $\sigma_i\in \Pi_2(\GL_{m_i}(D))$. En notant $\sigma_i'$ l'image réciproque dans $\Pi_2(\GL_{m_id}(F))$ de $\sigma_i$ par $\bLJ$, nous avons $\pi=\bLJ(\sigma_1'\times\cdots\times \sigma_l')$. Avec l'appui du théorème \ref{thm:classtempbysquareint} (2), on voit que $\bLJ^{-1}(\pi)\cap \Pi_\temp(\GL_{md}(F))$ un singleton.
\end{proof}

\subsection{Lemmes supplémentaires}
Soit $F$ un corps $p$-adique.

\begin{lemma}[Décomposition de la « fonction basique »]\label{lem:functionstratification} Soit $i\in \Z$.
Définissons 
\[I^{G,[i]}=\{g\in G(F)\mid |\nu_{\g}(g)|=q^{i}\}.\]
C'est un sous-ensemble ouvert de $G(F)$. On note $1_{I^{G,[i]}}$ sa fonction indicatrice.
\begin{enumerate}
    \item Soit $J\subseteq G(F)$ un sous-groupe ouvert compact. Alors $J\cap I^{G,[0]}$ est un sous-groupe ouvert compact de $G(F)$.
    \item Soient $J\subseteq G(F)$ un sous-groupe ouvert compact, et $f\in\S(\g(F)\ssla J)$. Soit $s\in \D(G)$. On a, sur $G(F)$,
    \begin{equation}\label{eq:eqfunctionstratification}
    f_s=\sum_{i\in \Z} q^{i\left(s+\frac{n-1}{2}\right)}(f\cdot 1_{I^{G,[i]}})    
    \end{equation} 
    une décomposition en une série de fonctions de $\H(G(F) \ssla J\cap I^{G,[0]})$ pour la convergence simple. Lorsque $s\in \D_\temp(G)$, cette décomposition converge dans l'espace de Schwartz-Harish-Chandra $\mathcal{HC}(G(F))$.
    
    \item Pour tout sous-groupe ouvert compact $J\subseteq G(F)$ et toute fonction $f\in\S(\g(F)\ssla J)$, il existe $i_f\in \Z$ tel que $f\cdot 1_{I^{G,[i]}}=0$ pour tout $i$ tel que $i\geq i_{f}$.
\end{enumerate} 
\end{lemma}
\begin{remark} Cette décomposition de la fonction basique commmute au transfert géométrique (la notion est aussi rappelée dans la définition \ref{def:transfergeomgeom'spec}), autrement dit $f\cdot 1_{I^{G,[i]}}\underset{\geom}{\arr} f^\ast \cdot 1_{I^{G^{\ast},[i]}}$ pour tout $i\in \Z$.
\end{remark}

\begin{proof}
Tout sauf l'énoncé à propos de la convergence dans $\mathcal{HC}(G(F))$ est clair. On prouve cette convergence : soit $j\in \Z$, alors la fonction
\begin{equation}
\phi_j\eqdef f_s-\sum_{i\in \Z:i>j} q^{i\left(s+\frac{n-1}{2}\right)}(f\cdot 1_{I^{G,[i]}})    
\end{equation}
est à support dans l'ouvert compact $\cup_{k\in \Z:k\leq j}I^{G,[k]}$, et cet ouvert compact est stable par $K$. On peut prendre un $j$ tel que $\phi_{j'}$ est bi-invariante par $K$ pour tout $j'\in \Z$ tel que $j'<j$. Il suffit ainsi d'étudier la suite des fonctions $(\phi_{j'})_{j'\in \Z:j'\leq j}$ restreintes sur $M_0(F)^+$. Avec les expressions explicites des fonctions $1+\sigma$, $\nu_\g$ et $e^{\rho_{P_0}(H_{M_0}(\cdot))}$, obtenues dans la preuve de la proposition \ref{prop:GJinSHC}, on conclut facilement que $\phi_{j'}$ tend vers $0$ dans l'espace de Schwartz-Harish-Chandra lorsque $j'\to -\infty$, ce qu'il fallait.
\end{proof}

\begin{lemma}\label{lem:heckematrixcoeff}Soit $(\pi,V)$ une représentation unitaire admissible de $G(F)$ et $J\subseteq G(F)$ un sous-groupe ouvert compact. Pour tout $v,v^\vee\in V^J$ il existe $h_{v,v^\vee}\in \H(G(F)\ssla J)$ tel que
\[\langle \pi(g)v,v^\vee\rangle=\tr(\pi(g)\pi(h_{v,v^\vee}))\]
\end{lemma}
\begin{proof}
On sait qu'il existe $h_{v,v^\vee}\in \H(G(F)\ssla J)$ tel que 
\[\pi(h_{v,v^\vee})\cdot w= \langle w,v^\vee\rangle v,\,\,\,\,\forall w\in V^J\]
grâce à la surjectivité de $\H(G(F)\ssla J)\rightarrow\End_{\C}(V^J)$ (cf. \cite[corollaire XVII.3.3]{Lg02}). Soit $(e_i)_i$ une base orthonormée de $V^J$. Nous avons
\[\tr(\pi(g)\pi(h_{v,v^\vee}))=\sum_{i}\langle \pi(g)\pi(h_{v,v^\vee})e_i,e_i\rangle=\sum_{i}\langle \langle e_i,v^\vee\rangle \pi(g)v,e_i\rangle=\langle \pi(g)v,v^\vee\rangle,\]
ce qu'il fallait.
\end{proof}

On a définit la notion de transfert géométrique des fonctions tests $\underset{\geom}{\arr}$ dans l'article précédent (cf. \cite[théorème 4.4]{YDL23b} ; cf. aussi définition \ref{def:transfergeomgeom'spec}).

\begin{lemma}\label{YDLspeclem:tracetransferiden}
Pour tout $\pi'\in \mathcal{R}(F)$, il existe $c(\pi')\in \R$ tel que pour tout $\S(\g(F))\ni \Phi\underset{\geom}{\arr} \Phi^\ast\in \S(\g^\ast(F))$ et $s\in\C$ avec $\Re(s)>c(\pi')$ on a 
\[e(G)\langle \bLJ(\pi'),\Phi_s\rangle=\langle\pi',\Phi_s^\ast\rangle.\]
\end{lemma}
\begin{proof}
Considérons d'abord le cas où $\pi'\in \Pi_{\temp}(G^\ast(F))$. Nous allons prouver que l'on peut prendre $c(\pi')=0$. Fixons $s\in\C$ avec $\Re(s)>0$. Soit $\S(\g(F))\ni \Phi\underset{\geom}{\arr} \Phi^\ast\in \S(\g^\ast(F))$. Comme $\Phi_s\in \mathcal{HC}(G(F))$ et $\Phi_s^\ast\in \mathcal{HC}(G^\ast(F))$ (proposition \ref{prop:GJinSHC}), la formule d'intégration de Weyl nous donne
\begin{equation*}
\begin{cases}
&e(G)\langle \bLJ(\pi'),\Phi_s\rangle=\sum_{T}|W^G(T)|^{-1}\int_{T(F)}|D^G(x)|e(G)\Theta_{\bLJ(\pi')}(x)J_G^G(x,\Phi_s)\,dx,\\
&\langle\pi',\Phi_s^\ast\rangle=\sum_{T'}|W^{G^\ast}(T')|^{-1}\int_{T'(F)}|D^{G^\ast}(x')|\Theta_{\pi'}(x')J_{G^\ast}^{G^\ast}(x',\Phi_s^\ast)\,dx',   
\end{cases}    
\end{equation*}
où la somme est sur l'ensemble des représentants des classes de $G(F)$-conjugaison (resp. $G^\ast(F)$-conjugaison) des tores maximaux de $G$ (resp. $G^\ast$). Ici pour $H$ un groupe réductif et $A$ un tore maximal de $H$ on note $W^H(A)$ le groupe de Weyl relatif de $A$ dans $H$. Puis on regarde uniquement les sous-parties $T(F)\cap G_\rss(F)$ et $T'(F)\cap G_\rss^\ast(F)$ dans les domaines d'intégration. On a $|D^{G^\ast}(x')|\Theta_{\pi'}(x')J_{G^\ast}^{G^\ast}(x',\Phi_s^\ast)=|D^G(x)|e(G)\Theta_{\bLJ(\pi')}(x)J_G^G(x,\Phi_s)$ si $T'=T^\ast$ se transfère (auquel cas $x'=x^\ast$ se transfère aussi et $|W^{G^\ast}(T^\ast)|=|W^G(T)|$), et $|D^{G^\ast}(x')|\Theta_{\pi'}(x')J_{G^\ast}^{G^\ast}(x',\Phi_s^\ast)=0$ si $T'=T^\ast$ ne se transfère pas. Or il existe $\Phi$ tel que $\langle\bLJ(\pi'),\Phi_s\rangle\not=0$ (soit $g\in G_{\rss}(F)$ tel que $\Theta_{\bLJ(\pi')}(g)\not=0$, on peut prendre $\Phi$ la fonction indicatrice d'un voisinage assez petit de $g$). On en déduit donc $e(G)\langle \bLJ(\pi'),\Phi_s\rangle=\langle\pi',\Phi_s^\ast\rangle$.

Soient ensuite $\pi'\in \Pi_{\temp}(G^\ast(F))$ et $\lambda^\ast\in a_{G^\ast}^\ast$. Puisque $\bLJ$ commute à la torsion par caractère (point 6 du théorème \ref{thm:appLJlocal}), on a $\langle \bLJ(\pi_{\lambda^\ast}'),\Phi_s\rangle=\langle \bLJ(\pi'),\Phi_{s+\lambda}\rangle$ et  $\langle\pi_{\lambda^\ast}',\Phi_s^\ast\rangle=\langle\pi',\Phi_{s+\lambda^\ast}^\ast\rangle$.  On sait que ce lemme est valable pour $\pi_{\lambda^\ast}'$. Enfin comme $\bLJ$ commute à l'induite parabolique (théorème \ref{thm:appLJlocal}), on sait que ce lemme est valable pour tout module standard. Puisque les modules standards forment une base de $\mathcal{R}(F)$, la démonstration se conclut par la linéarité.
\end{proof}

\subsection{Fonctions \texorpdfstring{$L$}{L}, \texorpdfstring{$\epsilon$}{e} et  \texorpdfstring{$\gamma$}{r} sous la correspondance de Jacquet-Langlands locale}

Soit $F$ un corps $p$-adique.

\begin{proposition}\label{prop:Lepsilongammalocal}~{}
\begin{enumerate}
    \item Soit $\pi'\in \Pi_{\unit}^d(G^\ast(F))$, alors
    \begin{align*}
    L(\pi',s)\mid L(|\bLJ|(\pi'),s),\,\,\epsilon(\pi',\psi,s)&=\epsilon(|\bLJ|(\pi'),\psi,s) \\
    \,\,\text{et}\,\,\gamma(\pi',\psi,s)&=\gamma(|\bLJ|(\pi'),\psi,s).\    
    \end{align*}
    Plus précisément le quotient de fonctions $L$ est un polynôme de $q^{-s}$ de terme constant 1.
    \item Soit $\pi^\ast\in \Pi_{\unit}^{d}(G^\ast(F))\cap \Pi_\temp(G^\ast(F))$, alors
    \[L(\pi^\ast,s)=L(\bLJ(\pi^\ast),s)\]
\end{enumerate}
\end{proposition}

\begin{remark}
Certaines des conclusions peuvent être trouvées dans \cite[section 6]{Badu08} et \cite[théorème 19.1]{Badu09}. Elles ont, en revanche, été démontrées par d'autres méthodes, à savoir des calculs explicites de ces quantités.
\end{remark}

\begin{proof}
Prouvons le point 1. Amorçons par la relation de division des fonctions $L$. Soit $(\pi,V)\in\Pi_{\unit}(G(F))$ une représentation, $v,v^\vee\in V$. Soit $J$ un sous-groupe ouvert compact de $G(F)$ tel que $v, v^\vee\in V^J$. Soit $(e_i)_{i\,\,\text{fini}}$ une base orthonormée de $V^J$. Pour toute fonction test $f\in \S(\g(F))$ et tout $s\in\D(G)$ on a 
{\allowdisplaybreaks
\begin{align}
Z(v,v^\vee,f,s) &= \int_{G(F)}\langle\pi(g)v,v^\vee\rangle f(g)|\nu_\g(g)|^s\,dg \notag\\
&=\int_{g\in G(F)}\int_{x\in G(F)}\sum_{i}\langle \pi(g)\pi(x)e_i,e_i\rangle h_{v,v^\vee}(x)f(g)|\nu_\g(g)|^s\,dx\,dg \,\,\,\,\text{(lemme \ref{lem:heckematrixcoeff})}\notag\\
&=\int_{y\in G(F)}\sum_{i}\langle \pi(y)e_i,e_i\rangle\left(\int_{g\in G(F)}h_{v,v^\vee}(g^{-1}y)f(g)|\nu_\g(g)|^s\, dg \right)\,dy \,\,\,\,(y=gx) \notag\\
&=\int_{G(F)}\sum_{i}\langle \pi(y)e_i,e_i\rangle\left(\int_{g'\in G(F)}h_{v,v^\vee}(g')f(yg'{}^{-1})|\nu_\g(yg'{}^{-1})|^s\, dg' \right)\,dy \,\,\,\,(g'=g^{-1}y)\notag\\
&=\int_{G(F)}\sum_{i}\langle \pi(y)e_i,e_i\rangle (h_{v,v^\vee}\ast f_s)(y)\,dy\notag\\
&=\sum_{i\text{ fini}}\sum_{j\in \Z,\text{ fini}} q^{-j\left(s+\frac{n-1}{2}\right)}Z(e_i,e_i,h_{v,v^\vee}^{[j]}\ast f,s)\label{eq:zetadiagonalsum}
\end{align}
}
dans la dernière ligne $h_{v,v^\vee}^{[j]}\eqdef h_{v,v^\vee}\cdot 1_{I^{G,[j]}}$ (cf. lemme \ref{lem:functionstratification}).

Soient maintenant $\pi'\in \Pi_{\unit}^d(G^\ast(F))$ et $\pi =|\bLJ|(\pi')$. Soit $\S(\g(F))\ni \Phi\underset{\geom}{\arr} \Phi^\ast\in \S(\g^\ast(F))$. Posons $J$ et $J'$ sous-groupes ouverts compacts de $G(F)$ et $G^\ast(F)$, tels que $\Phi\in \S(\g(F)\ssla J)$ et $\Phi^\ast\in \S(\g^\ast(F)\ssla J')$. Soit $(e_i')_{i\,\,\text{fini}}$ une base orthonormée de $V'{}^{J'}$. On 
a $\langle \pi ,\Phi_s\rangle=\varepsilon_{\pi'}\langle \pi' ,\Phi_s^\ast\rangle$ pour un signe $\varepsilon_{\pi'}=\pm 1$ uniquement déterminé. Comme le caractère de Harish-Chandra est, au sens des distributions, une somme diagonale de coefficients matriciels, on obtient que
\begin{align*}
L(\pi',s)^{-1}\sum_{i }Z(e_i,e_i,\Phi,s)&=\varepsilon_{\pi'}L(\pi',s)^{-1}\sum_{i }Z(e_i',e_i',\Phi^\ast,s)\\
&=\varepsilon_{\pi'}\sum_{i}Z(e_i',e_i',\Phi^\ast,s)L(\pi',s)^{-1}
\end{align*}
est une somme finie de polynômes de $q^{\pm s}$ donc l'est aussi. Cela est valable pour toute fonction $\Phi\in\S(\g(F))$.

On a, d'après l'équation \eqref{eq:zetadiagonalsum}, où l'on suppose que $J$ est assez petit de sorte que $f\in \S(\g(F)\ssla J)$ et $h_{v,v^\vee}^{[j]}\in \H(G(F)\ssla J)$ pour tout $j$,
\[L(\pi
',s)^{-1}Z(v,v^\vee,f,s)=\sum_{j\in\mathbb{Z}}\sum_{i}q^{-j\left(s+\frac{n-1}{2}\right)}\left(L(\pi',s)^{-1}Z(e_i,e_i,h_{v,v^\vee}^{[j]}\ast f,s)\right)\]
donc $L(\pi',s)^{-1}Z(v,v^\vee,f,s)$ est un polynôme de $q^{\pm s}$ pour tout $f$, $v$ et $v^\vee$, cela entraîne que $L(\pi',s)^{-1}L(\pi,s)\in \C[q^{\pm s}]$. Or ce quotient est de la forme $P_{\pi'}(q^{-s})/P_{\pi}(q^{-s})$ pour $P_{\pi'}$ et $P_\pi$ deux polynômes à coefficient constant $1$. On voit, en faisant $\R\ni s\to+\infty$, que 
\[\frac{L(\pi,s)}{L(\pi',s)}\in 1+q^{-s}\C[q^{- s}].\]

Secundo, à propos de l'égalité des facteurs gamma $\gamma(\pi',\psi,s)=\gamma(\pi,\psi,s)$ on raisonne ainsi : reprenons les notations plus hautes, soient $\pi'\in \Pi_{\unit}^{d}(G^\ast(F))$ et $\pi =|\bLJ|(\pi')$. Prenons $\S(\g(F))\ni \Phi\underset{\geom}{\arr} \Phi^\ast\in \S(\g^\ast(F))$ tels que $\langle \pi ,\Phi_s\rangle\not =0$. Posons $J$ et $J'$ sous-groupes ouverts compacts de $G(F)$ et $G^\ast(F)$, tels que $\Phi,\widehat{\Phi}\in \S(\g(F)\ssla J)$ et $\Phi^\ast,\widehat{\Phi^\ast}\in \S(\g^\ast(F)\ssla J')$. Soient $(e_i)_{i\,\,\text{fini}}$ une base orthonormée de $V{}^{J}$ et $(e_i')_{i\,\,\text{fini}}$ une base orthonormée de $V'{}^{J'}$. On a 
\begin{equation}\label{eq:gammaequ1}
\langle \pi ,\Phi_s\rangle=\varepsilon_{\pi'}\langle \pi' ,\Phi_s^\ast\rangle    
\end{equation}
pour un signe $\varepsilon_{\pi'}=\pm 1$ uniquement déterminé, cette égalité est valable pour tout $s$ par prolongement méromorphe compte tenu de ce que les deux côtés sont des sommes finies de fonctions zêtas locales, on en déduit 
\begin{equation}\label{eq:gammaequ2}
e(G)\langle \pi ,\widehat{\Phi}_{1-s}\rangle=\varepsilon_{\pi'}\langle \pi' ,\widehat{\Phi^\ast}_{1-s}\rangle.    
\end{equation}
En divisant l'égalité \eqref{eq:gammaequ2} par \eqref{eq:gammaequ1} et en utilisant l'autodualité de $(e_i)_i$ et $(e_i')_i$ on obtient $\gamma(\pi,\psi,s)=\gamma(\pi',\psi,s)$, l'égalité recherchée. 

Tertio, on sait que
\[(1+q^{-s}\C[q^{- s}])\cdot\frac{1}{1+q^{-s}\C[q^{- s}]}\ni \frac{L(\pi,s)}{L(\pi',s)}\cdot\frac{L(\pi'{}^\vee,1-s)}{L(\pi^\vee,1-s)}=\frac{\epsilon(\pi,\psi,s)}{\epsilon(\pi',\psi,s)}\]
est un quotient de monômes de $q^{-s}$, la seule possibilité est $\epsilon(\pi',\psi,s)=\epsilon(\pi,\psi,s)$.

Prouvons le point 2. On applique \cite[théorème 6.2]{Badu08} 
et le corollaire \ref{coro:localJLtemp}. 
\end{proof}

\subsection{Rappel des facteurs de normalisation \texorpdfstring{$r_{Q|P}(\pi_\lambda)$}{rQPrr}}\label{subsec:normalisationr}

Afin de procéder à la comparaison des formules des traces pour démontrer la correspondance de Jacquet-Langlands globale, nous devons choisir les facteurs de normalisation $r_{Q|P}(\pi_\lambda)$ de façon compatible au torseur intérieur. Soient $F$ un corps local de caractéristique 0 et $\eta:G=\GL_{m,D}\to G^\ast=\GL_{n,F}$ une forme intérieure.

Supposons pour le moment que $G=G^\ast$. Soient $M\in \L^G$, et $P,Q\in\P^G(M)$. On écrit $M=\prod_{i=1}^r \GL_{n_i,F}$. Toute $\pi$ représentation irréductible unitaire de $M(F)$ est de la forme $\pi_1\otimes\cdots\otimes \pi_r$ avec $\pi_i$ représentation irréductible unitaire de $\GL_{n_i}(F)$. Supposons que $\pi$ est tempéré. L'ensemble des racines $\Sigma(\g;A_M)$ est en bijection avec l'ensemble des paires $(i,j)$ d'entiers différents entre $1$ et $r$ : soit $(a_1\Id_{n_1\times n_1},\dots,a_r\Id_{n_r\times n_r})\in A_M(F)$ avec $a_i\in F^\times$, alors la racine qui correspond à $(i,j)$ est la racine $\alpha_{(i,j)}(a_1\Id_{n_1\times n_1},\dots,a_r\Id_{n_r\times n_r})=a_ia_j^{-1}$. Posons 
\[r_{\alpha(i,j)}(\pi,z)\eqdef\frac{L(\pi_i\otimes\pi_j^\vee,z)}{\epsilon(\pi_i\otimes\pi_j^\vee,\psi,z)L(\pi_i\otimes\pi_j^\vee,z+1)},\,\,\,\,z\in\C,\]
où les fonctions $L$ et $\epsilon$ sont celles de Rankin-Selberg. 

Enlèvons à présent l'hypothèse que $\eta$ est trivial. On définit la notion de transfert des sous-groupes de Levi et des sous-groupes paraboliques, comme dans le texte précédent \cite[définitions 2.7]{YDL23b}. On suit aussi les conventions de \cite[sous-section 2.4]{YDL23b} : on note, en ajoutant une étoile $^\ast$, en exposant les objets relativement à $G^\ast$ qui se transfèrent à $G$ ; pour les objets généraux relativement à $G^\ast$ on les note en ajoutant un signe prime $'$ en exposant. Par exemple en écrivant $M^\ast\in \L^{G^\ast}$ il est sous-entendu qu'il existe $M\in\L^G$ tel que $M^\ast$ se transfère $M$. Alors que l'écriture $M'\in \L^{G^\ast}$ n'impose pas de condition particulière sur $M'$. Le contexte devrait permettre d’enlever toute ambiguïté. Soit $M'\in\L^{G^\ast}$. Quand $M'$ se transfère on écrit $M'=M^\ast$ et $M'\arr M^\ast$, avec $M\in\L^G$ le transfert de $M^\ast$. Les mêmes notations valent pour les sous-groupes paraboliques. Supposons que $\pi$ est tempéré. Compte tenu de l'isomorphisme induit $\eta:A_M\xrightarrow{\sim}A_{M^\ast}$, toute racine $\alpha$ de $\Sigma(\g;A_M)$ correspond à une racine $\alpha^\ast$ de $\Sigma(\g^\ast;A_{M^\ast})$, définissons $r_{\alpha}(\pi,z)=r_{\alpha^\ast}(\pi^\ast,z)$ avec $\pi^\ast$ l'unique élément de $\bLJ^{-1}(\pi)\cap \Pi_{\temp}(G^\ast(F))$. Le facteur $r_\alpha(\pi,z)$ est méromorphe et vérifie $r_\alpha(\pi_\lambda,z)=r_\alpha(\pi,\langle\lambda,\alpha^\vee\rangle +z)$ pour $\lambda\in ia_{M}^\ast$. 

Pour $\pi$ général on écrit $\pi$ comme le quotient de Langlands d'un module standard $\I_R^M(\rho_\lambda)$ et définit $r_\alpha(\pi,z)=\prod_{\beta\in\Sigma(\m;A_{M_R}):\beta|_{A_{M}}=\alpha}r_\beta(\rho,\langle\lambda,\alpha^\vee\rangle +z)$. Posons enfin
\[r_{Q|P}(\pi_\lambda)=\prod_{\alpha\in \Sigma(\mathfrak{q};A_M)\cap(-\Sigma(\mathfrak{p};A_M))}r_\alpha(\pi, \langle\lambda,\alpha^\vee\rangle).\]


\begin{proposition}[{{\cite[chapitre 2, lemme 2.1]{AC}}}]
Les fonctions méromorphes $r_{Q|P}(\pi_\lambda)$ vérifient les conditions données dans le théorème 2.1 de \cite{Art89I}.    
\end{proposition}

\subsection{Côté spectral de la formule des traces non-invariante}\label{subsec:cotespecTFnon-inv} Adressons-nous ensuite à la correspondance de Jacquet-Langlands globale. Soit dorénavant $F$ un corps de nombres. Pour $H$ un groupe algébrique sur $F$ et $v\in \V_F$, on note $H_v$ le changement de base $H\times_FF_v$, aussi $H_{\infty}=H\times_FF_\infty$. On reprend dans la suite, jusqu'à la fin de la section \ref{sec:corrJL}, les notations ainsi que les normalisations des mesures de \cite[section 9]{YDL23b} : soit donc $\eta: \g=\gl_{m,D}\to \g^\ast=\gl_{n,F}$ une forme intérieure, et on suppose que  $\eta_\infty:\g_\infty\to \g_\infty^\ast$ est un torseur trivial. Fixons $\underline{S}_0$ un sous-ensemble fini de places, contenant les places archimédiennes, tel que pour tout $v\in (\V_F\setminus \underline{S}_0)\cup\V_{\infty}$ on ait que $\eta_v:G_v\to G_v^\ast$ est un torseur trivial. On fixe, pour tout $v\in (\V_F\setminus \underline{S}_0)\cup \V_{\infty}$, un élément $x_v\in G^\ast(\overline{F_v})$, tel que $(\Ad x_v)\circ\eta_v|_{G_v}$ soit un isomorphisme défini sur $F_v$ de $G_v$ sur $G_v^\ast$. 

Comme dans le numéro \ref{YDLspecsec:convention}, on prend, pour $v\in \V_\fin$, $\mathfrak{k}_v$ (resp. $\mathfrak{k}_v^\ast$) un réseau de $\g(F_v)$ (resp. $\g^\ast(F_v)$) ; pour $v\in \V_F$, $K_v$ (resp. $K_v^\ast$)un sous-groupe compact maximal de $G(F_v)$ (resp. $G^\ast(F_v)$). On exige que ces objets sont choisis de sorte que : si $v\in \V_F\setminus \underline{S}_0$ alors $(\Ad x_v)\circ\eta_v|_{\mathfrak{k}_v}$ est un isomorphisme de $\mathfrak{k}_v$ sur $\mathfrak{k}_v^\ast$, et si $v\in (\V_F\setminus \underline{S}_0)\cup \V_\infty$ alors $(\Ad x_v)\circ\eta_v|_{K_v}$ est un isomorphisme de $K_v$ sur  $K_v^\ast$.

Fixons également $\underline{S}\supseteq \underline{S}_0$ un sous-ensemble fini de places comme en \cite[sous-section 10.4]{YDL23b}. On suit les mêmes normalisations des mesures que dans la sous-section \ref{subsec:normalisationmesureofspec1}. Puis pour tout $v\in\underline{S}$, et tout tore maximal de $G_v$ ou de $G_v^\ast$ on munit son groupe des $F_v$-points une mesure de Haar, et on demande que si deux tores maximaux sont conjugués par un élément de $G(F_v)$ ou de $G^\ast(F_v)$ alors les mesures de leur groupes des $F_v$-points se correspondent par conjugaison par ce même élément. On exige que les mesures vérifient les compatibilités suivantes.
\begin{enumerate}[label=(\roman*)]
    \item Pour tout $v\in (\V_F\setminus \underline{S})\cup \V_{\infty}$, les mesures sur $G_v(F_v)$ et $G_v^\ast(F_v)$ se correspondent par l'isomorphisme $(\Ad x_v)\circ\eta_v|_{G_v(F_v)}$ ;
    \item pour tout $P\in \F^G$ ou $\F^{G^\ast}$, $\vol(N_P(F)\backslash N_P(\A_F))=1$ ;
    \item pour tout $v\in \underline{S}$,  $M\in\L^G(M_{0,v})$, le torseur intérieur $\eta_v$ induit un isomorphisme d'espaces  vectoriels réels $a_{M}\xrightarrow{\sim} a_{M^\ast}$, on veut que les mesures sur $a_{M}$ et $a_{M^\ast}$ se correspondent par cet isomorphisme ;
    \item le torseur intérieur $\eta$ induit un isomorphisme d'espaces  vectoriels réels $a_{G}\xrightarrow{\sim} a_{G^\ast}$, on veut que les mesures sur $a_{G}$ et $a_{G^\ast}$ se correspondent par cet isomorphisme, ainsi la mesure sur $A_{G,\infty}$ correspondent à celle sur $A_{G^\ast,\infty}$ par $\eta$ ;
    \item pour tout $v\in \underline{S}\setminus\V_{\infty}$, si $T$ est un tore maximal de $G_v$ et $T^\ast$ un tore maximal de $G_v^\ast$, et $y\in G_v^\ast(\overline{F_v})$ est tel que $(\Ad y)\circ\eta|_T$ soit un isomorphisme défini sur $F_v$ de $T$ sur $T^\ast$, alors les mesures sur $T(F_v)$ et $T^\ast(F_v)$ se correspondent par cet isomorphisme.
\end{enumerate}

\begin{definition}[{{\cite[définition 9.1]{YDL23b}}}]\label{def:deftransfertgeomoffunctioninspecart}
Soit $f\in \S(\g(\A_F))$ un tenseur pur. Nous dirons qu'un tenseur pur $f^\ast\in\S(\g^\ast(\A_F))$ se transfère de façon géométrique en $f$ et nous noterons $f\underset{\geom}{\arr}f^\ast$, ou plus simplement $f\arr f^\ast$, s'il existe des décompositions  $f=\bigotimes_{v\in \V_F} f_v$ et $f^\ast=\bigotimes_{v\in \V_F} f_v^\ast$ telles que les conditions suivantes sont remplies.
    \begin{enumerate}
    \item Pour tout $v\in \underline{S}\setminus\V_{\infty}$, on a $f_v\arr f_v^\ast$.
    \item Pour tout $v\in\V_{\infty}$, on a $f_v\circ ((\Ad x_v)\circ\eta_v)^{-1}=f_v^\ast$.
    \item Pour tout $v\in \V_F\setminus \underline{S}$, on  a $f_v=1_{\mathfrak{k}_v}$ et $f_v^\ast=1_{\mathfrak{k}_v^\ast}$.
    \end{enumerate}
\end{definition}

Pour tout $S\subseteq \V_F$ ensemble fini contenant $\underline{S}$, l'espace vectoriel $\S(\g(F_S))$ s'injecte dans $\S(\g(\A_F))$ par produit tensoriel avec $\prod_{v\not\in S}1_{\mathfrak{k}_v}$. De même pour $\S(\g^\ast(F_S))\hookrightarrow\S(\g^\ast(\A_F))$

Nous voulons avoir recours au développement fin du côté spectral de la formule des traces pour la déduction de la correpondance de Jacquet-Langlands globale. Il n'est pas possible d'utiliser une fonction de classe Schwartz-Bruhat comme fonction test pour des raisons de convergence. Il nous faut faire appel à l'espace $\mathcal{C}(G(\A_F))$ (définition \ref{def:espaceMuLaFi}).

\begin{proposition}\label{prop:innerformspecsideequality}Soient $\S(\g(\A_F))\ni f\underset{\geom}{\arr} f^\ast\in\S(\g^\ast(\A_F))$ et $s\in \D_\adel(G)$, alors
\[J_\text{spec}^G(f_s)=J_\text{spec}^{G^\ast}(f_s^\ast)\]
\end{proposition}
\begin{proof}
On aura besoin des estimations pour justifier certaines interversions entre intégrales, la tâche est réalisée par Matz dans \cite{Matz16}. Le lemme suivant est une conséquence immédiate de \cite[lemme 3.7]{Matz16} (on remarque que dans son article $\g(F)$ est décomposé selon les classes de conjugaison semi-simples et non les classes de conjugaison). Les objets $\O^\g$, $K_{\o}^{\g,T}$, et $J_{\geom}^{G,T}$ suivant sont définis en \cite[sections 6.2, 8.2]{YDL23b}.
\begin{lemma}Soit $T\in a_0$ assez régulier. Pour $a\in A_{G,\infty}$ et $f\in \S(\g(\A_F))$ soit $f^{a}\in\S(\g(\A_F))$ la fonction $f^{a}(X)=f(aX)$. Il existe une semi-norme $\|\cdot\|_T$ sur $\S(\g(\A_F))$ telle que   
\begin{align*}
\int_{A_{G,\infty}G(F)\backslash G(\A_F)}\int_{A_{G,\infty}}\sum_{\o\in \O^\g,\,\o\text{ non nilpotent}}|K_{\o}^{\g,T}(x,f^{a})|\,da\,dx\leq \|f\|_T\cdot \int_{A_{G,\infty}}|a|_{\A_F}^{-n^2}\,da
\end{align*}
\end{lemma}
Il est énoncé pour $F=\Q$ dans l'article, cela dit le cas général s'ensuit par une restriction des scalaires et le simple fait $|x|_{\A_F}=|N_{F/\Q}(x)|_{\A_\Q}$ avec $x\in \A_F^\times$ et $N_{F/\Q}$ la norme de l'extension $F/\Q$ (qui donne les normes locales).

Revenons sur la preuve. On a vu en \cite[théorème 9.8]{YDL23b} que si $\S(\g(\A_F))\ni f\underset{\geom}{\arr} f^\ast\in\S(\g^\ast(\A_F))$ alors $J_{\geom}^G(f)=J_{\geom}^{G^\ast}(f^\ast)$. Or il est clair que $f^{ a}|\nu_\g(a\cdot)|^{s+\frac{n-1}{2}}\arr f^{\ast, a^\ast}|\nu_{\g^\ast}(a^\ast\cdot)|^{s+\frac{n-1}{2}}$ pour tout $a\in A_{G,\infty}$, il vient ainsi
\begin{align*}
J_\text{spec}^G(f_s)&=J_{\geom}^G\left(\int_{A_{G,\infty}}f^{a}|\nu_\g(a\cdot)|^
{s+\frac{n-1}{2}}\,da\right) =\int_{A_{G,\infty}}J_{\geom}^G\left(f^{ a}|\nu_\g(a\cdot)|^{s+\frac{n-1}{2}}\right)\,da\\
&=\int_{A_{G^\ast,\infty}}J_{\geom}^{G^\ast}\left(f^{\ast, a^\ast}|\nu_\g(a^\ast\cdot)|^{s+\frac{n-1}{2}}\right)\,da^\ast=J_\text{spec}^{G^\ast}(f^\ast_s),    
\end{align*}
la permutation de $\int_{A_{G,\infty}}\,da$ et $J_{\geom}^G$ (aussi celle pour $G^\ast$) résulte de l'estimation de Matz ci-dessus, on passe par la distribution $J_\geom^T$ pour $T$ assez régulier et utilise le fait qu'il existe un entier positif $b=b(G,M_0,K)$ tel que $J_\geom^{T_0}$ soit une combinaison linéaire de $T^iJ_\geom^{T+kE_j}$ avec $0\leq i\leq b, -b\leq k\leq b$, et $0\leq j\leq \dim a_0$, où $(E_1,\dots,E_{\dim a_0})$ est une $\R$-base de $a_0$.
\end{proof}

\subsection{Isolement des représentations : théorie générale}

On utilise une technique de l'isolement employée  par Langlands en \cite[pp. 135-137]{Lang80}, puis formalisée par Lapid-Rogawski en \cite[lemme 4]{LaRo2000}.

Fixons $S\subseteq \V_F$ un ensemble fini contenant $\underline{S}$. Soit $U$ un ensemble fini de $\V_F$, disjoint de $S$, et contenant au moins deux places finies, disons $v_1$ et $v_2$, telles que les caractéristiques des corps résiduels de $F_{v_1}$ et de $F_{v_2}$ soient distincts.

On pose $\widehat{M_{U,0}}\eqdef \prod_{v\in U}\widehat{M_{v,0}}=\prod_{v\in U}\widehat{M_{v,0^\ast}}$, pour rappel $M_{v,0}$ (resp. $M_{v,0^\ast}$) est un sous-groupe de Levi minimal de $G_{v}$ (resp. $G_v^\ast$). On pose $\mathcal{X}^{G_U}\eqdef \prod_{v\in U} \mathcal{X}^{G_v}=\prod_{v\in U} \mathcal{X}^{G_v^\ast}$. On a un diagramme commutatif 
\[\begin{tikzcd}
	{a_{0^\ast,\C}^\ast} && {\widehat{M_{U,0}}} & {} \\
	\\
	{a_{0^\ast,\C}^\ast/W_{0^\ast}^{G^\ast}} && {\mathcal{X}^{G_U}}
	\arrow["{\exp_U}", from=1-1, to=1-3]
	\arrow["{\text{quotient}}"', from=1-1, to=3-1]
	\arrow[from=3-1, to=3-3]
	\arrow["{\text{quotient}}", from=1-3, to=3-3]
\end{tikzcd}\]
Tout sous-espace $a_M^\ast$ avec $M\in\L^G$ sera identifié à un sous-espace de $a_{0^\ast}^\ast$ via le torseur intérieur $\eta$. 

Soit $w\in W_0^G$. Soit $M\in\L^G$. On note $w\cdot M\eqdef wMw^{-1}$. Soit $\lambda \in i(a_M^G)^{\ast}$. On note $w\cdot \lambda\in i(a_{w\cdot M}^G)^\ast$ le vecteur donné par l'action de $W_0^G$ sur $i(a_0^G)^\ast$. Soit $\pi\in \Pi_\unit(M(\A_F))$. On note $w\cdot \pi$ la représentation de $(w\cdot M)(\A_F)$ donnée par $(w\cdot \pi)(wm w^{-1})=\pi(m)$ agissant sur le même espace que celui sous-jacent de $\pi$, avec $m\in M(\A_F)$.

Posons $\Pi_\disc(M)^{S-\text{ur}}=\{\pi\in\Pi_\disc(M)\mid \pi^S\text{ est non-ramifiée}\}$. Soit $\mathcal{Y}\eqdef \{(M,\pi)\mid M\in \L^G,\pi\in \Pi_{\disc}(M)^{S-\text{ur}}\}$. On se donne pour chaque $(M,\pi)\in\mathcal{Y}$ une fonction continue $g_{(M,\pi)}:i(a_M^G)^\ast\to\C$. On suppose que \begin{enumerate}
    \item si $(M,\pi)\in \mathcal{Y}$ et $w\in W_0^G$, alors $g_{(M,\pi)}(\lambda)=g_{(w\cdot M,w\cdot \pi)}(w\cdot \lambda)$ ;
    \item $\sum_{(M,\pi)\in\mathcal{Y}}\int_{i(a_M^G)^\ast}|g_{(M,\pi)}(\lambda)|\,d\lambda<+\infty$.
\end{enumerate}  

Soit $\mathcal{Y}'\eqdef \{(M',\pi')\mid M'\in \L^{G^\ast},\pi'\in \Pi_{\disc}(M')^{S-\text{ur}}\}$. On se donne également pour chaque $(M',\pi')\in\mathcal{Y'}$ une fonction continue $g_{(M',\pi')}:i(a_{M'}^{G^\ast})^\ast\to\C$. On suppose que 
\begin{enumerate}
    \item si $(M',\pi')\in \mathcal{Y}'$ et $w'\in W_{0^\ast}^{G^\ast}$, alors $g_{(M',\pi')}(\lambda')=g_{(w'\cdot M',w'\cdot \pi')}(w'\cdot \lambda')$ ;
    \item $\sum_{(M',\pi')\in\mathcal{Y}'}\int_{i(a_{M'}^{G^\ast})^\ast}|g_{(M',\pi')}(\lambda')|\,d\lambda'<+\infty$.
\end{enumerate}

\begin{lemma}\label{lem:isolementtheoriegenerale}
Supposons que
\[\sum_{(M',\pi')\in\mathcal{Y}'}\int_{i(a_{M'}^{G^\ast})^\ast}g_{(M',\pi')}(\lambda')\phi_{U}(z_{\I_{P'}^{G^\ast}(\pi_{\lambda'}')_U})\,d\lambda'=\sum_{(M,\pi)\in\mathcal{Y}}\int_{i(a_M^G)^\ast}g_{(M,\pi)}(\lambda)\phi_{U}(z_{\I_P^G(\pi_{\lambda})_U})\,d\lambda\]
pour tout $U$ ensemble fini de $\V_F$, disjoint de $S$, contenant au moins deux places finies, disons $v_1$ et $v_2$, telles que les caractéristiques des corps résiduels de $F_{v_1}$ et de $F_{v_2}$ soient distinctes, et pour toute fonction $\phi_U\in C^0(\mathcal{X}_{\unit,\text{ur}}^{G_U};\C)$. Alors pour tout $M'\in \L^{G^\ast}$, tout $P' \in\P^{G^\ast}(M')$, tout $\pi_0'\in \Pi_{\disc}(M')$, et presque tous les $\lambda_0'\in i(a_{M'}^{G^\ast})^\ast$ (pour la mesure de Lebesgue), nous avons 
\begin{align*}
&\sum_{\substack{(\pi',\lambda')\in \Pi_\disc(M')^{S-\text{ur}}\times i(a_{M'}^{G^\ast})^\ast\\ z_{\Ind_{P'}^{G^\ast}(\pi_{\lambda'}')^S}=z_{\Ind_{P'}^{G^\ast}(\pi_{0,\lambda_0'}')^S}}}g_{(M',\pi')}(\lambda') \\
&=\begin{cases}
\sum_{\substack{(\pi,\lambda)\in \Pi_\disc(M)^{S-\text{ur}}\times i(a_M^G)^\ast\\  z_{\Ind_P^G(\pi_\lambda)^S}=z_{\Ind_{P'}^{G^\ast}(\pi_{0,\lambda_0'}')^S}}}g_{(M,\pi)}(\lambda) &, \text{si $M'=M^\ast$} \\
0 &, \text{si $M'$ ne se transfère pas}
\end{cases}
\end{align*} 
ici $P\in \P^G(M)$ si $M'=M^\ast$.
\end{lemma}
\begin{proof}
On fixe $M'\in \L^{G^\ast}$, $P' \in\P^{G^\ast}(M')$ et $\pi_0'\in \Pi_{\disc}(M')$. 

Pour $\sigma\in \mathcal{X}^{G_U}$, on note $1_\sigma$ la fonction indicatrice de $\{\sigma\}$ sur $\mathcal{X}^{G_U}$. Posons $\Pi_\disc(M)_{U-\text{ur}}=\{\pi\in\Pi_\disc(M)\mid \pi_U\text{ est non-ramifiée}\}$. On note $1_{\Pi_\disc(M)_{U-\text{ur}}}$ la fonction indicatrice de $\Pi_\disc(M)_{U-\text{ur}}$ sur $\Pi_\disc(M)$.

La méthode de \cite[démonstration du lemme 4, du début jusqu'à l'équation (24)]{LaRo2000} nous conduit à : pour presque tous les $\lambda_0'\in i(a_{M'}^{G^\ast})^\ast$, nous avons
\begin{align*}
&\sum_{(\pi',\lambda')\in \Pi_\disc(M')\times i(a_{M'}^{G^\ast})^\ast}g_{(M',\pi')}(\lambda') 1_{\Pi_\disc(M')_{U-\text{ur}}}(\pi')1_{z_{\Ind_{P'}^{G^\ast}(\pi_{0,\lambda_0'}')_U}}(z_{\Ind_{P'}^{G^\ast}(\pi_{\lambda'}')_U})\\
&=\begin{cases}
\sum_{(\pi,\lambda)\in \Pi_\disc(M)\times i(a_M^G)^\ast}g_{(M,\pi)}(\lambda) 1_{\Pi_\disc(M)_{U-\text{ur}}}(\pi)1_{z_{\Ind_{P'}^{G^\ast}(\pi_{0,\lambda_0'}')_U}}(z_{\Ind_{P}^{G}(\pi_{\lambda})_U})&, \text{si $M'=M^\ast$}\\
0 &, \text{si $M'$ ne se transfère pas}
\end{cases}
\end{align*} 
ici $P\in \P^G(M)$ si $M'=M^\ast$. En faisant tendre $U$ vers $\V_F\setminus S$ et en utilisant le théorème de convergence dominée, on en déduit le lemme.
\end{proof}

\subsection{Isolement des représentations}\label{subsec:isolementdesreprésentations} 
Fixons $S\subseteq \V_F$ un ensemble fini contenant $\underline{S}$. Prenons $U$ un ensemble fini de $\V_F$, disjoint de $S$, et contenant au moins deux places finies, disons $v_1$ et $v_2$, telles que les caractéristiques des corps résiduels de $F_{v_1}$ et de $F_{v_2}$ soient distinctes. Soient $f\in \S(\g(F_{S}))\subseteq\S(\g(\A_F))$, $\phi_U\in \H(G(F_U)\ssla K_U)$, et $s\in\D_{\adel}(G)$. Nous tâchons de décomposer l'égalité dans la proposition \ref{prop:innerformspecsideequality} en suivant la méthode de Lapid-Rogawski évoquée plus haute. Le point de départ est le développement fin (théorème \ref{thm:devlopfinspec}) :
\begin{equation*}
\begin{split}
J_\spec^G(\phi_U\ast (f_s))=\sum_{M\in\L^G}\sum_{\pi\in\Pi_{\disc}(M)}&\int_{i(a_{M}^G)^\ast}\sum_{L\in\L^G(M)}|W_0^M||W_0^G|^{-1}\\
&\cdot a_\disc^M(\pi)r_M^L(\pi_{\lambda})\tr\left(\mathcal{N}_{L}(P,\pi,\lambda)\I_{P}^G(\pi_{\lambda},(\phi_U\ast (f_s))^1\right)\,d\lambda  
\end{split}    
\end{equation*}
avec $P\in \P^G(M)$.

\begin{lemma}
Soient $M\in \L^G$ et $L\in \L^G(M)$. L'application
\[\lambda\in i(a_M^G)^\ast\longmapsto\sum_{\pi\in\Pi_{\disc}(M)}\sum_{L\in\L^G(M)}|W_0^M||W_0^G|^{-1}a_\disc^M(\pi)r_M^L(\pi_{\lambda})\tr\left(\mathcal{N}_{L}(P,\pi,\lambda)\I_{P}^G(\pi_{\lambda},f_s^1)\right)\]
est continue.
\end{lemma}
\begin{proof}
En vertu des majorations obtenues au cours de la démonstration du théorème \ref{thm:devlopfinspec}, il nous suffit de prouver que pour tout $\pi\in\Pi_\disc(M)$ fixé, les applications $\lambda\in i(a_M^G)^\ast\longmapsto r_M^L(\pi_{\lambda})$ et $\lambda\in i(a_M^G)^\ast\longmapsto\tr\left(\mathcal{N}_{L}(P,\pi_{\lambda})\I_{P}^G(\pi_{\lambda},f_s^1)\right)$ sont continues. On sait d'après la définition que $\lambda\mapsto r_M^L(\pi_{\lambda})$ est continue. Ensuite l'opérateur $\I_{P}^G(\pi_{\lambda},f_s^1)$ est non-nul uniquement si $\pi^{S}$ est non-ramifié. Le cas échéant d'après la construction d'Arthur (cf. \cite[théorème 2.1 (R8)]{Art89I}) et la formule de scindage pour une $(G,M)$-famille il est clair que, pour tout $\lambda\in i(a_M^G)^\ast$, 
{\allowdisplaybreaks
\begin{align*}
\tr(\mathcal{N}_L&(P,\pi_\lambda)\I_P^G(\pi_\lambda,f_s^1) )\\
&=\left(\sum_{(L_v)_{v\in S}\in \L^{G_{S}}(M_{S})}d_M^G\left((L_v)_{v\in S}\right)\prod_{v\in S} \tr\left(\mathcal{N}_{M_v}^{Q_{L_v}}(P_v,\pi_{v,\lambda})\I_{P_v}^{G_v}(\pi_
{v,\lambda}, f_{s,v}))\right)\right)\\
&\hspace{10cm}\cdot \prod_{i=1}^rL^{S}(\pi_i,s+\lambda_i),
\end{align*}
on a décomposé au passage $M\simeq \prod_{i=1}^r\GL_{m_i,D}$, $\pi=\pi_1\otimes\cdots\otimes\pi_r$ et $\lambda=(\lambda_1,\dots,\lambda_r)$. Cette expression est continue selon la construction des facteurs de normalisation. }
\end{proof}

\begin{lemma}
Soit $w\in W_0^G$. Soit $M\in\L^G$. Soit $P\in\P^G$. On note $w\cdot P\eqdef wPw^{-1}$. Soit $L\in \L^G(M)$. Soit $\lambda \in i(a_M^G)^{\ast}$. Soit $\pi\in \Pi_\unit(M(\A_F))$. Alors
\[a_\disc^M(\pi)=a_\disc^{w\cdot M}(w\cdot \pi), r_M^L(\pi_{\lambda})=r_{w\cdot M}^{w\cdot L}((w\cdot\pi)_{w\cdot \lambda}),\]
et
\[\tr\left(\mathcal{N}_{L}(P,\pi,\lambda\right)\I_{P}^G(\pi_{\lambda},f_s^1)=\tr\left(\mathcal{N}_{w\cdot  L}(w\cdot P,w\cdot \pi,w\cdot \lambda)\I_{w\cdot P}^G((w\cdot \pi)_{w\cdot \lambda},f_s^1\right).\] 
\end{lemma}
\begin{proof}
L'action de $w$ définit un $F$-isomorphisme $M\to w\cdot M$ qui envoie $\pi$ sur $w\cdot \pi$. D'où $a_\disc^M(\pi)=a_\disc^{w\cdot M}(w\cdot \pi)$. Pour les deux autres quantités on utilise la condition $(R_5)$ de \cite[théorème 2.1]{Art89I}. 
\end{proof}

Soit $\phi_U\in \H(G(F_U)\ssla K_U)$. On pose 
$\phi_U^\ast\eqdef \phi_U\circ ((\Ad x_U)\circ\eta_U|_{G_U})^{-1}\in \H(G^\ast(F_U)\ssla K_U^\ast)$. Pour rappel, on a fixé pour tout $v\in (\V_F\setminus \underline{S}_0)\cup \V_{\infty}$, un élément $x_v\in G^\ast(\overline{F_v})$, tel que $(\Ad x_v)\circ\eta_v|_{G_v}$ soit un isomorphisme défini sur $F_v$ entre $G_v$ et $G_v^\ast$. On fera convoler la fonction $f_s$ par $\phi_U$, afin d'aboutir à la décomposition voulue. Partant de la proposition \ref{prop:innerformspecsideequality} nous voyons

\begin{equation}\label{eq:foliationstep1}
\begin{split}
&\sum_{M\in\L^G}\sum_{\pi\in\Pi_{\disc}(M)}\int_{i(a_{M}^G)^\ast}\sum_{L\in\L^G(M)}|W_0^M||W_0^G|^{-1}\\
&\hspace{3cm}\cdot a_\disc^M(\pi)r_M^L(\pi_{\lambda})\tr\left(\mathcal{N}_{L}(P,\pi,\lambda)\I_{P}^G(\pi_{\lambda},f_s^1)\right)\S\phi_{U}(z_{\I_P^G(\pi_{\lambda})_U})\,d\lambda  \\
&=\sum_{M'\in\L^{G^\ast}}\sum_{\pi'\in\Pi_{\disc}(M')}\int_{i(a_{M'}^{G^\ast})^\ast}\sum_{L'\in\L^{G^\ast}(M')}|W_{0^\ast}^{M'}||W_{0^\ast}^{G^\ast}|^{-1}\\
&\hspace{3cm}\cdot a_\disc^{M'}(\pi')r_{M'}^{L'}(\pi_{\lambda'}')\tr\left(\mathcal{N}_{L'}(P',\pi',\lambda')\I_{P'}^{G^\ast}(\pi_{\lambda'}',f_s^{\ast,1})\right)\S\phi_{U}^\ast(z_{\I_{P'}^{
G^\ast}(\pi_{\lambda'}')_U})\,d\lambda'  
\end{split}    
\end{equation}
Constatons que le morphisme de Satake $\S:\H(G(F_U)\ssla K_U)\rightarrow C^0(\mathcal{X}_{\unit,\text{ur}}^{G_U};\C)$ est d'image dense. En l'occurrence l'égalité \eqref{eq:foliationstep1} reste valable quand on remplace $\S \phi_{U}$ par toute fonction dans $C^0(\mathcal{X}_{\unit,\text{ur}}^{G_U};\C)$, compte tenu du fait que les deux côtés sont des distributions sur $\mathcal{X}_{\unit,\text{ur}}^{G_U}$ au vu du théorème \ref{thm:devlopfinspec} et la compacité de l'espace.

Grâce au lemme \ref{lem:isolementtheoriegenerale}, on trouve :

\begin{lemma}[Isolement des représentations]\label{lem:isolement} 
Pour tout $M'\in \L^{G^\ast}$, tout $P' \in\P^{G^\ast}(M')$, tout $\pi_0'\in \Pi_{\disc}(M')$, et presque tous les $\lambda_0'\in i(a_{M'}^{G^\ast})^\ast$, nous avons 
\begin{align*}
&\sum_{\substack{(\pi',\lambda')\in \Pi_\disc(M')^{S-\text{ur}}\times i(a_{M'}^{G^\ast})^\ast\\ z_{\Ind_{P'}^{G^\ast}(\pi_{\lambda'}')^S}=z_{\Ind_{P'}^{G^\ast}(\pi_{0,\lambda_0'}')^S}}}\sum_{L'\in\L^{G^\ast}(M')}a^{L'}(\I_{P'\cap L'}^{L'}(\pi_{\lambda'}'))J_{L'}^{G^\ast}\left(\I_{P'\cap L'}^{L'}(\pi_{\lambda'}'),f_s^\ast\right) \\
&=\begin{cases}
\sum_{\substack{(\pi,\lambda)\in \Pi_\disc(M)^{S-\text{ur}}\times i(a_M^G)^\ast\\  z_{\Ind_P^G(\pi_\lambda)^S}=z_{\Ind_{P'}^{G^\ast}(\pi_{0,\lambda_0'}')^S}}}\sum_{L\in\L^{G}(M)}a^L(\I_{P\cap L}^{L}(\pi_\lambda))J_L^G\left(\I_{P\cap L}^{L}(\pi_\lambda),f_s\right) &, \text{si $M'=M^\ast$} \\
0 &, \text{si $M'$ ne se transfère pas}
\end{cases}
\end{align*} 
ici $P'\in \P^{G^\ast}(M')$, et $P\in \P^G(M)$ si $M'=M^\ast$.
\end{lemma}

\subsection{Correspondance de Jacquet-Langlands globale et les fonctions \texorpdfstring{$L$}{L} et \texorpdfstring{$\epsilon$}{e}}


Soient toujours $G=\GL_{m,D}$ avec $D$ une algèbre à division centrale sur $F$ de degré $d$, et $G^\ast=\GL_{n,F}$. Soit $v$ une place. On écrit $G_v=\GL_{m_v,D_v}$ avec $D_v$ une algèbre à division centrale sur $F_v$ de degré $d_v$. Clairement $n=md=m_vd_v$. On écrit $\bLJ_v:\mathcal{R}(F_v)\to \mathcal{R}(D_v)$ le morphisme d'anneaux obtenu dans le théorème \ref{thm:appLJlocal}.

Une représentation dans $\Pi_{\disc}(G(\A_F)^1)$ est dite $D$-compatible si elle est $d_v$-compatible en toute place $v$. \'{E}crivons $\Pi_{\disc}^{D}(G^\ast(\A_F)^1)$ l’ensemble des classes d'isomorphismes des représentations automorphes au sens $L^2$ discrètes $D$-compatibles.

Prenons $M'=G^\ast$ dans le lemme \ref{lem:isolement} on obtient une identité d'une sous-somme de la partie discrète du développement fin spectral de la formule des traces non-invariantes de $G^\ast(\A_F)$ et celle correspondante de $G(\A_F)$. Or la partie discrète de la formule des traces non-invariantes coïncide avec celle de la formule des traces invariantes. On reprend les idées de \cite[théorème 5.1]{Badu08} pour prouver le résultat suivant : la correspondance globale de Jacquet-Langlands.


\begin{theorem}[{{\cite[théorème 18.1]{Badu09}}}]\label{thm:corresglobaleJLBadu} ~{}
\begin{enumerate}
    \item Il existe une unique application $\bG:\Pi_{\disc}(G(\A_F)^1)\rightarrow \Pi_{\disc}(G^\ast(\A_F)^1)$ telle que pour toute $\pi\in \Pi_{\disc}(G(\A_F)^1)$, si $\pi^\ast=\bG(\pi)$ alors $|\bLJ|_v(\pi_v^\ast)=\pi_v$ pour toute place $v\in \V_F$. L'application $\bG$ est injective, d'image $\Pi_{\disc}^{D}(G^\ast(\A_F)^1)$.
    \item Pour toute représentation $\pi\in\Pi_{\disc}(G(\A_F)^1)$, $m_\disc(\pi)=1$.
\end{enumerate}
\end{theorem}

Esquissons la démonstration : on fait une récurrence sur le rang de $G$, puis on utilise le théorème fort de multiplicité un pour $G^\ast$ et le fait que $M_{P|P}(s,0)$ (resp. $M_{P'|P'}(s',0)$) est un opérateur unitaire pour estimer les coefficients $a^{G^\ast}(\pi')$ et $a^G(\pi)$, ensuite on utilise le lemme \ref{YDLspeclem:tracetransferiden} pour savoir que $J_{G^\ast}^{G^\ast}\left(\pi^{\ast,S},f_s^{\ast,S}\right)=J_G^G\left(\pi^S,f_s^S\right)$ et $J_{G^\ast}^{G^\ast}\left(\pi_S',f_{S,s}\right)=\pm J_G^G\left(|\bLJ_S|(\pi_S'),f_{S,s}\right)$, et on conclut par l'indépendance linéaire des caractères sur le groupe $G(F_S)$. Dans la preuve nous pouvons omettre la sommation partielle — utilisée chez Arthur — suivant les modules des caractères infinitésimaux à l'infini désormais inutile puisque nous savons que le développement spectral de la formule des traces non-invariante est absolument convergent.

Il est possible d'étendre l'application $\bG$ sur le spectre $L^2$ tout entier, résultat découlant du théorème \ref{thm:irredunitparaind}. On identifiera dans la suite $a_{\GL_{m,D},\C}^\ast$ à $\C$ en envoyant $\nu_{\gl_{m,D}}\otimes 1$ sur $1$.

\begin{proposition}\label{prop:rigidificationL2spec} Soient $\pi_j^{(t)}\in \Pi_{\disc}(\GL_{m_j^{(t)},D}(\A_F)^1)$, $\lambda_j^{(t)}\in i a_{\GL_{m_j^{(t)},D}}^\ast$,
$j\in \{1,\dots,k_{t}\}$, et $t\in\{1,2\}$, avec $\sum_{j=1}^{k_t} \lambda_j^{(t)}=0$ et $\sum_{j=1}^{k_1}m_j^{(1)}=\sum_{j=1}^{k_2}m_j^{(2)}$. S'il existe $S'$ un sous-ensemble fini de $\V_F$ contenant les places archimédiennes tel que $\times_{j=1}^{k_t}\pi_{j,\lambda_j^{(t)}}^{(t)}$ pour $t=1,2$ possèdent les mêmes composantes locales non-ramifiées hors $S'$, alors $k_1=k_2$ et il existe une permutation $\tau$ de $\{1,\dots,k_1\}$ telle que $\pi_j^{(2)}=\pi_{\tau(j)}^{(1)}$ et $\lambda_j^{(2)}=\lambda_{\tau(j)}^{(1)}$ pour tout $j$.
\end{proposition}

\begin{proof}
Supposons d'abord que $D=F$. La même preuve qu'en \cite[proposition 4.1]{Badu08} donne $k_1=k_2$ et l'existence d'une permutation $\tau$ telle que $\pi_{j,\lambda_j^{(2)}}^{(2)}=\pi_{\tau(j),\lambda_{\tau(j)}^{(1)}}^{(1)}$. Il est claire que   $\pi_j^{(2)}=\pi_{\tau(j)}^{(1)}$ et $\lambda_j^{(2)}=\lambda_{\tau(j)}^{(1)}$ pour tout $j$. Pour $D\not=F$ on prend $m=\sum_{j=1}^{k_1}m_j^{(1)}$, $G=\GL_{m,D}$, et on reprend les notations précédentes. La forme intérieure induit un isomorphisme $\eta:a_M\xrightarrow{\sim} a_{M^\ast}$ pour tout $M\in \L^G$. La propriété de $\bG$ donne que $\times_{j=1}^{k_t}\bG\left(\pi_{j}^{(t)}\right)_{\eta\left(\lambda_j^{(t)}\right)}$ pour $t=1,2$ possèdent les mêmes composantes locales non-ramifiées hors $S'\cup S$ (pour rappel $S\supseteq \underline{S}$), d'où $k_1=k_2$ et l'existence d'une permutation $\tau$ telle que $\bG\left(\pi_{j}^{(2)}\right)=\bG\left(\pi_{\tau(j)}^{(1)}\right)$ et $\eta\left(\lambda_j^{(2)}\right)=\eta\left(\lambda_{\tau(j)}
^{(1)}\right)$, on conclut alors par l'injectivité de $\bG$ et $\eta$.
\end{proof}

La propriété ci-dessus détermine complètement $\Pi_{L^2}(G(\A_F)^1)$, l’ensemble des classes d’isomorphismes des sous-quotients irréductibles de $L^2(G(F)\backslash G(\A_F)^1)$. On prolonge la définition de $\bG$ sur $\Pi_{L^2}(G(\A_F)^1)$ comme suit : supposons que $M=\prod_{j=1}^k \GL_{m_j,D}$ un sous-groupe de Levi de $G=\GL_{m,D}$, $\pi_j\in \Pi_{\disc}(\GL_{m_j,D}(\A_F)^1)$ et $\lambda_j\in ia_{\GL_{m_j,D}}^\ast$ des nombres complexes avec $\sum_{j=1}^k \lambda_j=0$. Alors on définit
\begin{equation}\label{eq:defJ-LglobalL2}
\begin{split}
 \bG : \Pi_{L^2}(G(\A_F)^1)&\longrightarrow \Pi_{L^2}(G^\ast(\A_F)^1)\\
 \times_{j=1}^k \pi_{j,\lambda_j}&\longmapsto \times_{j=1}^k\bG(\pi_{j})_{\eta(\lambda_j)}.   
\end{split}    
\end{equation}
L'application est bien définie d'après la proposition ci-dessus. 

De même une représentation dans $\Pi_{L^2}(G(\A_F)^1)$ est dite $D$-compatible si elle est $d_v$-compatible en toute place. \'{E}crivons $\Pi_{L^2}^{D}(G^\ast(\A_F)^1)$ l’ensemble des classes d'isomorphismes des représentations automorphes $D$-compatible. Il est facile de voir, grâce à la proposition \ref{prop:d-compcommutewithind}, que si $M'=\prod_{j=1}^k \GL_{n_j,F}$ est un sous-groupe de Levi de $G^\ast=\GL_{n,F}$, $\pi_j'\in \Pi_{\disc}(\GL_{n_j,F}(\A_F)^1)$ et $\lambda_j'\in ia_{\GL_{n_j,F}}^\ast$ des nombres complexes avec $\sum_{j=1}^k \lambda_j'=0$, alors $\times_{j=1}^k \pi_{j,\lambda_j'}'$ est $D$-compatible si et seulement si $M'=M^\ast$ se transfère et $\pi_j'\in\Pi_{\disc}^D(\GL_{n_j,F}(\A_F)^1)$ pour tout $j$. On déduit de ces discussions le théorème suivant.
\begin{theorem}\label{thm:globalJLL2}~{}
\begin{enumerate}
    \item L'application $\bG:\Pi_{L^2}(G(\A_F)^1)\rightarrow \Pi_{L^2}(G^\ast(\A_F)^1)$ satisfait, pour toute $\pi\in \Pi_{L^2}(G(\A_F)^1)$, si $\pi^\ast=\bG(\pi)$ alors $|\bLJ_v|(\pi_v^\ast)=\pi_v$ pour toute place $v\in \V_F$. Elle est unique pour cette propriété. L'application $\bG$ est injective, d'image $\Pi_{L^2}^{D}(G^\ast(\A_F)^1)$.
    \item Nous avons $L\left(\bG\left(\times_{j=1}^k\pi_{j,\lambda_j}\right),s\right)=\prod_{j=1}^k L\left(\bG\left(\pi_{j}\right)_{\eta(\lambda_j)},s\right)$ et $\epsilon\left(\bG\left(\times_{j=1}^k\pi_{j,\lambda_j}\right),s\right)=\prod_{j=1}^k\epsilon\left(\bG\left(\pi_{j}\right)_{\eta(\lambda_j)},s\right)$.
\end{enumerate}    
\end{theorem}

\begin{corollary}\label{coro:LandgammaunderglobalJL}
Soient $\pi\in \Pi_{L^2}(G(\A_F)^1)$ et $\pi^\ast=\bG(\pi)$ alors 
\[L(\pi^\ast,s)\mid L(\pi,s)\,\,\,\,\text{et}\,\,\,\,\epsilon(\pi^\ast,s)=\epsilon(\pi,s).\]    
\end{corollary}
\begin{proof}
Ceci est une conséquence directe du théorème précédent et la proposition \ref{prop:Lepsilongammalocal}.    
\end{proof}

\begin{remark}~{}
En général une égalité $L(\pi,s)=L(\pi^\ast,s)$ est attendue lorsque $\pi^\ast$ est dans le spectre cuspidal, auquel cas $\pi$ l'est aussi. Cela résultera de la conjecture de Ramanujan-Petersson pour les groupes du type GL (toute représentation cupidale est localement partout tempérée) et de la proposition \ref{prop:Lepsilongammalocal}. Toutefois quand on passe à l'autre extrémité du spectre discret la situation dégénère compte tenu du fait que la fonction $L$ n'est pas respectée par la formation d'un quotient. \`{A} titre d'exemple soient $D$ une algèbre à division centrale sur $F$ de degré $d$ quasi-déployée aux places archimédiennes, $G=\GL_{1,D}$ et $G^\ast=\GL_{d,F}$. On vérifie aisément (en manipulant \cite[proposition 3.7]{Badu08}) 
que la représentation triviale $1_{G(\A_F)^1}$ est envoyée sur la représentation triviale $1_{G^\ast(\A_F)^1}$ par $\bG$. Nous avons
\[\frac{L(1_{G(\A_F)^1},s)} {L(1_{G^\ast(\A_F)^1},s)}=\prod_{v\in \V_{\fin}}\prod_{j\in\{1,\dots,m_v-1\}}(1-q_v^{-s+d_vj-\frac{d-1}{2}})\not=1,\] 
avec $d_v$ l'ordre de $D\otimes_FF_v$ dans le groupe de Brauer $\text{Br}(F_v)$ et $m_v=d/d_v$.

\end{remark}


Pour $L\in\L^G$, on notera $\bG^L:\Pi_{L^2}(L(\A_F)^1)\to \Pi_{L^2}(L^\ast(\A_F)^1)$ la correspondance de Jacquet-Langlands globale de $L(\A_F)^1$. Elle donne une bijection $\Pi_{L^2}(L(\A_F)^1)\to \Pi_{L^2}^D(L^\ast(\A_F)^1)$. 
\begin{definition}
On dit que $\pi'\in \Pi_{L^2}(L^\ast(\A_F)^1)$ se transfère à $L$ si $\pi'\in\Pi_{L^2}^D(L^\ast(\A_F)^1)$. On note $\pi'=\pi^\ast$ et $\pi\arr \pi^\ast$ avec $\pi$ l'image réciproque de $\pi^\ast$ par $\bG^L$. On dit que $\pi^\ast$ se transfère à $L$ en $\pi$.     
\end{definition}

\begin{proposition}\label{prop:globalJ-Lindtrans}~{}
Soient $M'\in\L^{G^\ast}$, $ P'\in\P^{G^\ast}(M')$, $L^\ast\in \L^{G^\ast}(M')$ qui se transfère à $G$ en $L$, $\pi'\in \Pi_\disc(M')$, et $\lambda'\in i(a_{M'}^{L^\ast})^\ast$. Alors $\I_{P'\cap L^\ast}^{L^\ast}(\pi_{\lambda'}')$ se transfère à $L$ si et seulement si $M'=M^\ast$ et $\pi'=\pi^\ast$ via $\bG^M$. Le cas échéant, on aura $\I_{P^\ast\cap L^\ast}^{L^\ast}(\pi_{\lambda^\ast}^\ast)=\I_{P\cap L}^L(\pi_\lambda)^\ast\eqdef \bG^L(\I_{P\cap L}^L(\pi_\lambda))$.
\end{proposition}
\begin{proof}Prouvons d'abord l'équivalence. \'{E}bauchons par l'implication directe. Soient $A'\in \L^{M'}$, $Q'\in \P^{G^\ast}(A')$, et $\sigma'\in \Pi_\disc(A'(\A_F)^1)$ tels que $\pi'= \I_{Q'\cap M^\ast}^{M^\ast}(\sigma')$. L'hypothèse dit que $\I_{Q'\cap L^\ast}^{L^\ast}(\sigma_{\lambda'}')$ se transfère. De la proposition \ref{prop:rigidificationL2spec} et la définition de $\bG^L$ nous tirons que $A'=A^\ast$ se transfère et  $\sigma'=\sigma^\ast$ se transfère via $\bG^A$. Ainsi $M'=M^\ast$. Puis par la définition de $\bG^M$ on a $\pi'=\I_{Q^\ast\cap M^\ast}^{M^\ast}(\sigma^\ast)=\bG^M(\I_{Q\cap M}^{M}(\sigma))$. Ce qu'il fallait.

Procédons à la réciproque. Supposons que $M'=M^\ast$ et $\pi'=\pi^\ast$. Soient $A\in \L^M$, $Q\in \P^G(A)$, et $\sigma\in \Pi_\disc(A(\A_F)^1)$ tels que $\pi= \I_{Q\cap M}^M(\sigma)$. Soient $A'\in \L^{M^\ast
}$, $Q'\in \P^{G^\ast}(A')$, et $\sigma'\in \Pi_\disc(A'(\A_F)^1)$ tels que $\pi'= \I_{Q'\cap M^\ast}^{M^\ast}(\sigma')$. Par la définition de $\bG^M$ on a $\I_{Q'\cap M^\ast}^{M^\ast}(\sigma')=\pi'=\pi^\ast=\I_{Q^\ast\cap M^\ast}^{M^\ast}(\sigma^\ast)$. Par l'intermédiaire de la proposition \ref{prop:rigidificationL2spec} on voit que l'on peut supposer, quitte à conjuger, que $A'=A^\ast$ et $\sigma'=\sigma^\ast$. Enfin par on voit que $\I_{P^\ast\cap L^\ast}^{L^\ast}(\pi_{\lambda^\ast}^\ast)=\I_{Q^\ast\cap L^\ast}^{L^\ast}(\sigma_{\lambda^\ast}^\ast)=\bG^L(\I_{Q\cap L}^{L}(\sigma_{\lambda}))$. Ce qu'il fallait.
\end{proof}

\section{Identités spectrales non-invariantes de transfert}\label{sec:Idenspec}

\subsection{Identité globale}\label{subsec:idenspecglo}
Soit $F$ un corps de nombres. Pour $L\in\L^G$, on a défini un ensemble $\Pi(L)$ (équation \eqref{eq:defPi(G)triplets}) et une fonction $a^L(-):\Pi(L)\to\C$ (équation \eqref{eq:defcoeffspeca^Gtriplets}). Soient $L\in \L^G$, $P\in \P^G(L)$, $Q\in \F^G(L)$, $\lambda\in a_{L,\C}^\ast$ et $\pi\in \Pi_\unit(L(\A_F)^1)$. On déduit de l'équation \eqref{eq:opN} le poids $\mathcal{N}_L^Q(P,\pi,\lambda)$, et cette fonction méromorphe de $\lambda$ est analytique sur $ia_L^\ast$. On écrit $\mathcal{N}_L^Q(P,\pi)=\mathcal{N}_L^Q(P,\pi,0)$. On défini les caractères pondérés globaux en suivant l'équation \eqref{eqdef:weightedchara} :
\begin{equation}
J_L^Q(\pi,h)\eqdef\tr(\mathcal{N}_L^Q(P,\pi)\I_P^G(\pi,h^1))
\end{equation}
pour une fonction $h\in \mathcal{C}(G(\A_F))$. Cette expression ne dépend pas de $P\in \P^G(L)$. On peut aussi le considérer comme une fonction $J_L^Q(-,h):\Pi_\unit(L(\A_F)^1)\to\C$.

On rappelle deux formules utiles qui résultent des combinatoires des $(G,L)$-familles (cf. numéro \ref{subsec:GMfamilles}).

\begin{proposition}Soient $L\in\L^G$, $Q\in \F^G(L)$, et $h\in \mathcal{C}(G(\A_F))$.
\begin{enumerate}
    \item (Formule de descente de l'induction) Soient $M\in\L^L$ et $\pi\in \Pi_\unit(M(\A_F)^1)$. On a 
    \[J_L^Q(\I_P^L(\pi),h)=\sum_{M_1\in\L^{L_Q}(M)}d_M^{L_Q}(L,M_1)J_M^{Q_{M_1}}(\pi,h).\]
    Ici $P\in\P^G(M)$ et $Q_{M_1}$ est la deuxième composante de $s(L,M_1)$ de l'équation \eqref{eq:formuledescentGMsec} (avec $G$ remplacé par $L_Q$).
    \item (Formule de descente parabolique) Soit $\pi\in \Pi_\unit(L(\A_F)^1)$. On a,
    \[J_L^Q(\pi,h)=J_L^{L_Q}(\pi,h_Q)\]
    avec $h_Q$ la fonction appartenant à $\mathcal{C}(L_Q(\A_F))$ définie par
    \[h_Q(l)=\gamma(Q)\int_{K}\int_{N_Q(\A_F)}h\left((\Ad k^{-1})(ln)\right)\,dn\,dk,\,\,\,\,\forall l\in L_Q(\A_F),\]
    et $\gamma(Q)>0$ est la dérivée de Radon-Nikodym de la mesure sur $G(\A_F)$ par rapport à la mesure produit sur $L_Q(\A_F)N_Q(\A_F)K=G(\A_F)$. On rappelle que la descente parabolique $h\mapsto h_Q$ définit une application linéaire continue $\mathcal{C}(G(\A_F))\to\mathcal{C}(L_Q(\A_F))$ (cf. \cite[lemme 2.1 (4)]{FiLa16}).
\end{enumerate}    
\end{proposition}

Reprenons dans cette sous-section les notations de la sous-section \ref{subsec:isolementdesreprésentations} et les normalisations des mesures fixées depuis la sous-section \ref{subsec:cotespecTFnon-inv}. 

\begin{theorem}\label{thm:idenendoglobal}~{}
\begin{enumerate}
    \item Soient $\L^G\ni M\arr M^\ast\in\L^{G^\ast}$, $\P^G(M)\ni P\arr P^\ast\in\P^{G^\ast}(M^\ast)$, $\L^G(M)\ni L\arr L^\ast\in \L^{G^\ast}(M^\ast)$, $\pi\in \Pi_\disc(M)$, et $\lambda\in i(a_M^L)^\ast$. On note $\pi^\ast=\bG^M(\pi)\in \Pi_\disc(M^\ast)$ et $\lambda^\ast=\eta(\lambda)\in i(a_{M^\ast}^{L^\ast})^\ast$. Alors nous avons
    \[a^L(\I_{P\cap L}^L(\pi_\lambda))=a^{L^\ast}(\I_{P^\ast\cap L^\ast}^{L^\ast}(\pi_{\lambda^\ast}^\ast)).\]
    \item Soient $\S(\g(\A_F),\psi)\ni f\underset{\geom}{\arr} f^\ast\in \S(\g^\ast(\A_F),\psi)$, $L'\in\L^{G^\ast}$, $Q'\in \F^{G^\ast}(L')$, $\lambda'\in ia_{L'}^\ast$ et $\pi
'\in\Pi_{L^2}(L'(\A_F)^1)$. Nous avons, pour tout $s\in \D_{\adel}(G)$, 
    \begin{equation}\label{eq:globalspecidenwant}
    J_{L'}^{Q'}(\pi_{\lambda'}',f_s^\ast) =\begin{cases}J_L^Q(\pi_\lambda,f_s) &, \text{ si $(L',Q',\lambda',\pi')=(L^\ast,Q^\ast,\lambda^\ast,\pi^\ast\eqdef\bG^L(\pi)$)} \\
    0 &, \text{ sinon}.
    \end{cases}    
    \end{equation}
    Plus précisément, « sinon » signifie que  $(L',Q',\lambda')=(L^\ast,Q^\ast,\lambda^\ast)$ et $\pi'\not\in \Pi_{L^2}^{D}(L^\ast(\A_F)^1)$, ou que $L'$ ne se transfère pas.
\end{enumerate}
\end{theorem}

\begin{proof}
Primo, prouvons le point 1. On suppose d'abord que $\pi\in \Pi_{\disc}(M(\A_F)^1)$. Selon la définition $a^L(\I_{P\cap L}^L(\pi_\lambda))=a_{\disc}^M(\pi)r_M^L(\pi_{\lambda})$. Or $a_{\disc}^M(\pi)$ est la multiplicité de $\pi$ dans $L_{\disc}^2(M(F)\backslash M(\A_F)^1)$, on a $a_{\disc}^M(\pi)=1$ (théorème \ref{thm:corresglobaleJLBadu}). Il ne reste qu'à prouver 
\begin{equation}\label{YDLspeceq:globalweicaracompinter1}
r_M^L(\pi_{\lambda})=r_{M^\ast}^{L^\ast}(\pi_{\lambda^\ast}^\ast).    
\end{equation}
Or le problème est local, on est conduit à prouver
\begin{equation}\label{eq:compglobalspeccoeff}
r_{M_v}^{L_v}(|\bLJ_v|(\pi'))=r_{M_v^\ast}^{L_v^\ast}(\pi')    
\end{equation}
pour tout $\pi'\in\Pi_{\unit}^{d_v}(G^\ast(F_v))$. On va abréger pour le moment $F_v$ en $F$. 

Nous allons introduire des notations supplémentaires. Esquissons au départ des résultat de classification de Bernstein-Zelevinsky et Tadi\'{c}. Soient $k,q\in \N_{>0}$ des entiers stricitement positifs. Soit $\rho$ une représentation de $\GL_q(F)$, on désigne par $\rho(c)$ la représentation $\rho$ tordue par $|\det|^c$ pour $c\in\C$. Soit maintenant $\rho$ supercuspidale, on appele $[\rho,\rho(1),\dots,\rho(k-1)]=\Delta$ un segment. L'induite $\rho\times\cdots\times \rho(k-1)$ est un module standard et on note $Q(\Delta)$ son quotient de Langlands, qui est en outre essentiellement de carré intégrable. Toute représentation essentiellement de carré intégrable est de la forme $Q(\Delta)$ pour un segment $\Delta$ uniquement déterminé. On écrit $\Delta(c)=[\rho(c),\rho(c+1),\dots,\rho(c+k-1)]$ pour $c\in\C$. Soient $\Delta_1=[\rho_1,\rho_1(1),\dots,\rho_1(k_1-1)]$ et $\Delta_2=[\rho_2,\rho_2(1),\dots,\rho_2(k_2-1)]$ deux segments, on dit que $\Delta_1$ et $\Delta_2$ sont liés si $\Delta_1\not\subseteq \Delta_2$, $\Delta_2\not\subseteq \Delta_1$ et $\Delta_1\cup \Delta_2$ est un segment. On dit que $\Delta_1$ précède $\Delta_2$ s'ils sont liés et $\rho_2=\rho_1(c)$ pour un $c\in\N$. Soient $\Delta_1,\dots,\Delta_c$ des segments, vérifiant que si $i<j$ alors $\Delta_i$ ne précède pas $\Delta_j$, alors $Q(\Delta_1)\times\cdots \times Q(\Delta_c)$ admet un unique quotient irréductible noté $Q(\Delta_1,\dots,\Delta_c)$. 

Afin se coller aux notations de Badulescu \cite{Badu08} on pose $Z^u(\rho,k)=Q\left(\left[\rho\left(-\frac{k-1}{2}\right)\right.\right.,\linebreak[0]\left.\left.\rho\left(-\frac{k-3}{2}\right),\dots,
\rho\left(\frac{k-1}{2}\right)\right]\right)$
pour $\rho$ supercuspidale unitaire. \'{E}galement $u(\sigma,k)\eqdef Q\left(\Delta\left(\frac{k-1}{2}\right),\Delta\left(\frac{k-3}{2}\right),\right.\linebreak[0]\dots,\left.\Delta\left(-\frac{k-1}{2}\right)\right)$ pour $\sigma=Q(\Delta)$ essentiellement de carré intégrable. Finalement $\pi(u(\sigma,k),\alpha)=|\det|^\alpha u(\sigma,k)\times |\det|^{-\alpha} u(\sigma,k)$ pour $\alpha\in ]0,\frac{1}{2}[$ les séries complémentaires, qui sont irréductibles unitaires. 

Grâce aux travaux de Tadi\'{c} \cite[théorème D, p.338]{Ta86}, on sait que tout produit $\times$ de représentations de la forme $u(\sigma,k)$ ou $\pi(u(\sigma,k),\alpha)$ pour $\sigma$ de carré intégrable, $k\in \N_{>0}$ et $\alpha\in ]0,\frac{1}{2}[$ sont irréductibles unitaires. Réciproquement toute représentation $\pi$ irréductible unitaire de $G^\ast(F)$ est un tel produit, les facteurs sont déterminés par $\pi$, à une permutation près.

Soit $i_{G^\ast}$ l'involution de Zelevinsky-Aubert de $G^\ast(F)=\GL_n(F)$. Soient $k,l,q$ trois entiers avec $klq=n$, et $\rho$ une représentation supercuspidale unitaire de $\GL_q(F)$. En exécutant l'algorithme de Mœglin-Waldspurger \cite{MW86} on trouve
{\allowdisplaybreaks
\begin{align*}
i_{G^\ast}Q&\left(\left[\rho\left(-\frac{k-1}{2}\right),\rho\left(-\frac{k-3}{2}\right),\dots,\rho\left(\frac{k-1}{2}\right)\right]\left(\frac{l-1}{2}\right)\right.,\\
&\hspace{1cm}\left[\rho\left(-\frac{k-1}{2}\right),\rho\left(-\frac{k-3}{2}\right),\dots,\rho\left(\frac{k-1}{2}\right)\right]\left(\frac{l-3}{2}\right),\\
&\hspace{2cm}\dots,\\
&\hspace{3cm}\left.\left[\rho\left(-\frac{k-1}{2}\right),\rho\left(-\frac{k-3}{2}\right),\dots,\rho\left(\frac{k-1}{2}\right)\right]\left(-\frac{l-1}{2}\right)\right)\\
=Q&\left(\left[\rho\left(-\frac{l-1}{2}\right),\rho\left(-\frac{l-3}{2}\right),\dots,\rho\left(\frac{l-1}{2}\right)\right]\left(\frac{k-1}{2}\right)\right.,\\
&\hspace{1cm}\left[\rho\left(-\frac{l-1}{2}\right),\rho\left(-\frac{l-3}{2}\right),\dots,\rho\left(\frac{l-1}{2}\right)\right]\left(\frac{k-3}{2}\right),\\
&\hspace{2cm}\dots,\\
&\hspace{3cm}\left.\left[\rho\left(-\frac{l-1}{2}\right),\rho\left(-\frac{l-3}{2}\right),\dots,\rho\left(\frac{l-1}{2}\right)\right]\left(-\frac{k-1}{2}\right)\right)
\end{align*}
}
Soit $i_{G^\ast}(u(Z^u(\rho,k),l))=u(Z^u(\rho,l),k)$.

Revenons sur la preuve de l'équation \eqref{eq:compglobalspeccoeff}. Comme $\bLJ$ est un morphisme d'anneaux, il nous suffit de vérifier l'équation \eqref{eq:compglobalspeccoeff} pour les représentations de la forme $\pi'=u(\sigma,k)$. L'image de $u(\sigma,k)$ par $\bLJ$ est explitée dans \cite[proposition 3.7]{Badu08}, on a alors seulement deux cas (a) et (b) à vérifier. En cas de (a) l'équation \eqref{eq:compglobalspeccoeff} est triviale par définition ; en cas de (b) l'équation \eqref{eq:compglobalspeccoeff} provient du fait que l'involution de Zelevinsky-Aubert commute, à un signe près à $\bLJ$ (\cite[théorème 3.15]{Badu07}) et l'équation $i_{G^\ast}(u(Z^u(\rho,k),l))=u(Z^u(\rho,l),k)$ que l'on a vérifiée. La preuve se conclut. Le cas général $\pi\in\Pi_{\disc}(M)$ se traite de la même façon, en utilisant la comparaison d'opérateurs d'entrelacement \cite[corollaire 5.4]{Badu08}.

Secundo, prouvons le point 2. De nouveau $F$ est un corps de nombres.  Constatons que si 
\[\S(\g(\A_F),\psi)\ni f\underset{\geom}{\arr}f^\ast\in \S(\g^\ast(\A_F),\psi) \,\,\,\,\text{ relativement à $\eta:\g\to\g^\ast$},\]
alors pour tout $\F^{G}\ni P\arr P^\ast\in \F^{G^\ast}$ on a 
\[\S(\m_P(\A_F),\psi)\ni f_P\underset{\geom}{\arr}f_{P^\ast}^\ast\in \S(\m_{P^\ast}(\A_F),\psi)\,\,\,\,\text{relativement à $\eta|_{\m_P}:\m_P\to\m_{P^\ast}$}.\]
Aussi, pour tout $P'$ qui ne se transfère pas,  
\[\S(\m_{P'}(\A_F),\psi)\ni f_{P'}^\ast \underset{\geom}{\arr}0 \in \S(\m_{P'}(\A_F),\psi)\,\,\,\,\text{relativement à $
\Id:\m_{P'}\to \m_{P'}$}.\]

Faisons une récurrence sur $\dim G$ pour prouver l'assertion. Si $\dim G=1$, alors $\eta:\g\to\g^\ast$ est un torseur intérieur trivial. Afin de prouver \eqref{eq:globalspecidenwant} on utilise le lemme \ref{YDLspeclem:tracetransferiden} puis le corollaire \ref{coro:CVAholocotespec}, i.e. le prolongement holomorphe des deux côtés sur $\D_\adel(G)$. Supposons que la proposition est établie pour $\dim G<m\in\mathbb{N}_{>1}$. On se place maintenant au cas où $\dim G=m$. 

Soient $\S(\g(\A_F),\psi)\ni f\underset{\geom}{\arr} f^\ast\in \S(\g^\ast(\A_F),\psi)$, $L'\in\L^{G^\ast}$, $Q'\in \F^{G^\ast}(L')$, $\lambda_0'\in ia_{L'}^\ast$, $\pi_0'\in\Pi_{L^2}(L'(\A_F)^1)$ et $s\in \D_{\adel}(G)$. On veut prouver 
\begin{equation}
J_{L'}^{Q'}(\pi_{0,\lambda_0'}',f_s^\ast) =\begin{cases}J_L^Q(\pi_{0,\lambda_0},f_s) & \text{ si $(L',Q',\lambda_0',\pi_0')=(L^\ast,Q^\ast,\lambda_0^\ast,\pi_0^\ast\eqdef\bG^L(\pi_0)$)} \\
0 & \text{ sinon}.
\end{cases}    
\end{equation} 
Procédons à une récurrence sur $\dim a_{L'}^{L_{Q'}'}$. Si $\dim a_{L'}^{L_{Q'}'}=0$, alors on peut supposer que $L'=L_{Q'}'=G^\ast$, sinon on applique la formule de descente parabolique et on utilise l'hypothèse de récurrence sur $\dim G$. L'égalité voulue résulte encore une fois du lemme \ref{YDLspeclem:tracetransferiden} et de l'holomorphie des termes sur $s\in\D_\adel(G)$.

Supposons que le théorème vaut au rang $\dim a_{L'}^{L_{Q'}'}<n\in\mathbb{N}\setminus\{0\}$. On se place dans le cas de $\dim a_{L'}^{L_{Q'}'}=n$. Par la formule de descente de l'induction, la formule de descente parabolique, et l'hypothèse de récurrence sur $\dim G$, on peut et on va supposer que $Q'=G^\ast$, $L'=L^\ast$, $\lambda_0'=\lambda_0^\ast$ et $\pi_0'=\pi_0^\ast\in\Pi_\disc(L^\ast(\A_F)^1)$.

On prend $P\in \P^{G}(L)$ et $P^\ast\in \P^{G^\ast}(L^\ast)$. On suppose pour l'instant que $\lambda_0^\ast \in i(a_{L^{\ast}}^{G^\ast})^\ast$ et que le lemme \ref{lem:isolement} est valide pour le quadruple  $(M',P',\pi_0',\lambda_0')=(L^\ast,P^\ast,\pi_0^\ast,\lambda_0^\ast)$. On applique le lemme \ref{lem:isolement} à $(M',P',\pi_0',\lambda_0')=(L^\ast,P^\ast,\pi_0^\ast,\lambda_0^\ast)$. Il est clair que la somme 
\[\sum_{\substack{(\pi',\lambda^\ast)\in \Pi_\disc(L^\ast)^{S-\text{ur}}\times i(a_{L^\ast}^{G^\ast})^\ast\\ z_{\Ind_{P^\ast}^{G^\ast}(\pi_{\lambda^\ast}')^S}=z_{\Ind_{P^\ast}^{G^\ast}(\pi_{0,\lambda_0^\ast}^\ast)^S}}}\]
dans ce lemme se réduit à
\[\sum_{(\pi^\ast,\lambda^\ast)\in \{(w^\ast\cdot\pi_0^\ast,w^\ast\cdot\lambda_0^\ast)\mid w^\ast\in W^{G^\ast}(L^\ast)\}}.\]
En effet cette propriété découle de la proposition \ref{prop:rigidificationL2spec} et le fait que si $w^\ast\in W^{G^\ast}(L^\ast)$ alors $w^\ast\cdot \Pi_\disc(L^\ast)^{S-\text{ur}}=\Pi_\disc(L^\ast)^{S-\text{ur}}$. De même, la somme 
\[\sum_{\substack{(\pi,\lambda)\in \Pi_\disc(L)^{S-\text{ur}}\times i(a_{L}^{G})^\ast\\ z_{\Ind_{P}^{G}(\pi_{\lambda})^S}=z_{\Ind_{P^\ast}^{G^\ast}(\pi_{0,\lambda_0^\ast}^\ast)^S}}}\]
dans ce lemme se réduit à
\[\sum_{(\pi,\lambda)\in \{(w\cdot\pi_0,w\cdot\lambda_0)\mid w\in W^{G}(L)\}}.\]
En conclusion, ce lemme nous dit que
\begin{align*}
&\sum_{(\pi^\ast,\lambda^\ast)\in \{(w^\ast\cdot\pi_0^\ast,w^\ast\cdot\lambda_0^\ast)\mid w^\ast\in W^{G^\ast}(L^\ast)\}}\sum_{A^\ast\in\L^{G^\ast}(L^\ast)}a^{A^\ast}(\I_{P^\ast\cap A^\ast}^{A^\ast}(\pi_{\lambda^\ast}^\ast))J_{A^\ast}^{G^\ast}\left(\I_{P^\ast\cap A^\ast}^{A^\ast}(\pi_{\lambda^\ast}^\ast),f_s^\ast\right) \\
&=\sum_{(\pi,\lambda)\in \{(w\cdot\pi_0,w\cdot\lambda_0)\mid w\in W^{G}(L)\}}\sum_{A\in\L^{G}(L)}a^A(\I_{P\cap A}^{A}(\pi_{\lambda}))J_L^G\left(\I_{P\cap A}^{A}(\pi_{\lambda}),f_s\right). 
\end{align*} 
En manipulant l'hypothèse de récurrence sur $n$ et la comparaison des coefficients $a$ établie plus haut, on en déduit que 
\begin{align*}
\frac{|W^{G^\ast}(L^\ast)|}{|\text{Stab}_{W^{G^\ast}(L^\ast)}(\pi_0^\ast,\lambda_0^\ast)|}&a^{L^\ast}(\I_{L^\ast}^{L^\ast}(\pi_{0,\lambda_0^\ast}^\ast))J_{L^\ast}^{G^\ast}\left(\pi_{0,\lambda_0^\ast}^\ast,f_s^\ast\right) \\
&=
\frac{|W^{G}(L)|}{|\text{Stab}_{W^{G}(L)}(\pi_0,\lambda_0)|}a^L(\I_{L}^{L}(\pi_{0,\lambda}))J_L^G\left(\pi_{0,\lambda_0},f_s\right). 
\end{align*} 
Ici $\text{Stab}_{W^{G^\ast}(L^\ast)}(\pi_0^\ast,\lambda_0^\ast)\eqdef\{w^\ast\in W^{G^\ast}(L^\ast)\mid (w^\ast\cdot\pi_0^\ast,w^\ast\cdot\lambda_0^\ast)=(\pi_0^\ast,\lambda_0^\ast)\}$ et $\text{Stab}_{W^{G}(L)}(\pi_0,\lambda_0)\eqdef\{w\in W^{G}(L)\mid (w\cdot\pi_0,w\cdot\lambda_0)=(\pi_0,\lambda_0)\}$. On voit que $|W^{G^\ast}(L^\ast)|=|W^{G}(L)|$. On voit aussi que $|\text{Stab}_{W^{G^\ast}(L^\ast)}(\pi_0^\ast,\lambda_0^\ast)|=|\text{Stab}_{W^{G}(L)}(\pi_0,\lambda_0)|$ en passant par le caractère de Harish-Chandra. Or $a^{L^\ast}(\I_{L^\ast}^{L^\ast}(\pi_{0,\lambda_0^\ast}^\ast))=a_\disc^{L^\ast}(\pi_0^\ast)=1=a^L(\I_{L}^{L}(\pi_{0,\lambda}))$. On en déduit que 
\[J_{L^\ast}^{G^\ast}\left(\pi_{0,\lambda_0^\ast}^\ast,f_s^\ast\right)=J_L^G\left(\pi_{0,\lambda_0},f_s\right).\]
Cette égalité a lieu pour presque tous les $\lambda_0^\ast\in i(a_{L^\ast}^{G^\ast})^\ast$. On conclut par la continuité que cette égalité a lieu pour tout les $\lambda_0^\ast\in i(a_{L^\ast}^{G^\ast})^\ast$. Finalement, soit $\lambda_0^\ast\in ia_{L^\ast}^\ast$. On le décompose en $\lambda_0^\ast=\lambda_1^\ast+\lambda_2^\ast$ avec $\lambda_1^\ast\in i(a_{L^\ast}^{G^\ast})^\ast$ et $\lambda_2^\ast\in ia_{G^\ast}^\ast$. Alors 
\[J_{L^\ast}^{G^\ast}\left(\pi_{0,\lambda_0^\ast}^\ast,f_s^\ast\right)=J_{L^\ast}^{G^\ast}\left(\pi_{0,\lambda_1^\ast}^\ast,f_{s+\lambda_2^\ast}^\ast\right)=J_L^G\left(\pi_{0,\lambda_1},f_{s+\lambda_2}\right)=J_L^G\left(\pi_{0,\lambda_0},f_{s}\right).\]
La deuxième égalité vient du fait évident $ia_{G}^\ast+\D_\adel(G)=\D_\adel(G)$. La preuve se termine.
\end{proof}

À la différence de l'approche d'Arthur-Clozel \cite{AC}, où les analyses géométriques et spectrales sont menées simultanément, et où l'égalité $I_\disc^G=I_\disc^{G^\ast}$ ainsi que la correspondance globale de Jacquet-Langlands découlent à la suite de ces analyses, notre approche non-invariante adopte une démarche différente. Tout d'abord, nous nous appuyons sur des analyses purement géométriques pour établir l'égalité du côté géométrique $J_{\geom}^G(f_s)=J_{\geom}^{G^\ast}(f_s^\ast)$. L'égalité du côté spectral $J_{\spec}^G(f_s)=J_{\spec}^{
G^\ast}(f_s^\ast)$ en est ensuite une conséquence immédiate. En mobilisant enfin des arguments d'isolement et des analyses purement spectrales, nous parvenons à la correspondance globale de Jacquet-Langlands et l'identité globale de transfert.

\subsection{Formule des traces locale}\label{subsec:TFlocal}

Préalablement à la discussion sur l'identité locale, nous procédons à un rappel de la formule des traces locale pour les groupes du type GL. Il convient de noter que la formulation présentée ici demeure spécifique à notre étude. Soit dans cette sous-section $F$ un corps local. Les normalisations des mesures ne seront pas explicitées car la validité de la formule des traces locale ne dépend d'aucun choix de mesure.

\subsubsection{Côté géométrique} Le contenu de ce numéro est emprunté à \cite[sections 4.8, 4.9]{Walds12}.

Soit $L\in \L^G$, on choisit pour chaque classe de $L(F)$-conjugaison d'un tore elliptique de $L$ un représentant, et on note $T_{\text{ell}}(L)$ l'ensemble des représentants. Soit $T\in T_{\text{ell}}(L)$, on note $W^L(T)$ le groupe de Weyl relatif de $(L,T)$. Soit $x\in G(F)$ semi-simple, on pose $D^G(x)=\det(1-\Ad(x);G/G_x)$. Soient $L\in \L^G$ et $g_1,g_2\in G(F)$, on pose $(v_P(\lambda,g_1,g_2))_{P\in \P^G(L)}$ la $(G,L)$-famille $v_P(\lambda,g_1,g_2)=e^{-\lambda(H_P(g_1)-H_{\overline{P}}(g_2))}$, avec $\overline{P}$ l'opposé de $P$. 

Pour $c\in\C$ on note $\overline{c}\in \C$ son conjugué complexe. Soient $L\in \L^G$, $T\in T_{\text{ell}}(L)$, $x\in T(F)$, et $h_1,h_2\in \mathcal{HC}(G(F))$. On introduit les intégrales orbitales pondérées doubles 
\[J_L^G(x,h_1,h_2)=|D^G(x)|\int_{T(F)\backslash G(F)}\int_{T(F)\backslash G(F)}\overline{h_1(g_1^{-1}xg_1)}h_2(g_2^{-1}xg_2)v_L^G(g_1,g_2)\,dg_1\,dg_2,\]
cette intégrale double ne dépend pas de $x$ mais uniquement de sa classe de $L(F)$-conjugaison. Puis on pose
\[J_{\geom,L}^G(h_1,h_2)=\sum_{T\in T_{\text{ell}}(L)}|W^L(T)|^{-1}\vol(A_L(F)\backslash T(F))\int_{T(F)}J_L^G(x,h_1,h_2)\,dx,\]
enfin
\[J_{\geom}^G(h_1,h_2)=\sum_{L\in\L^G(M_0)}|W_0^L||W_0^G|^{-1}(-1)^{\dim a_L^G}J_{\geom,L}^G(h_1,h_2).\]

On rappelle trois formules utiles qui résultent des combinatoires des $(G,L)$-familles (cf. numéro \ref{subsec:GMfamilles}).

\begin{proposition}Soient $L\in\L^G$, $Q\in \F^G(L)$, et $h_1,h_2\in \mathcal{HC}(G(F))$.
\begin{enumerate}
    \item Soit $x\in L(F)\cap G_{\rss}(F)$. On pose $D^\g(x)=\det(\ad(x);\g/\g_x)$. On a
    \begin{equation}\label{eq:TFlocaldescentgeom}
    J_L^G(x,h_1,h_2)=\frac{|D^G(x)|}{|D^{\g}(x)|}\sum_{(L_1,L_2)\in \L^G(L)^2} d_L^G(L_1,L_2)\overline{J_L^{Q_
    {L_1}}(x,h_1)}J_L^{Q_{L_2}}(x,h_2),    
    \end{equation}
    où $(L_1,L_2)\mapsto (Q_{L_1},Q_{L_2})$ est l'application \eqref{eq:formuledescentGMsec} (avec $M$ remplacé par $L$).
    \item (Formule de descente de l'induction) Soient $M\in\L^L$ et $x\in M(F)$. On a 
    \[J_L^Q(\Ind_M^L(x),h_1)=\sum_{M_1\in\L^{L_Q}(M)}d_M^{L_Q}(L,M_1)J_M^{Q_{M_1}}(x,h_1),\]
    où $Q_{M_1}$ est la deuxième composante de $s(L,M_1)$ de  \eqref{eq:formuledescentGMsec} (avec $G$ remplacé par $L_Q$).
    \item (Formule de descente parabolique) Soit $x\in L(F)$. On a,
    \[J_L^Q(x,h_1)=J_L^{L_Q}(x,h_{1,Q})\]
    avec $h_{1,Q}$ la fonction appartenant à $ \mathcal{HC}(L_Q(F))$ définie par
    \[h_{1,Q}(l)=\gamma(Q)\int_{K}\int_{N_Q(F)}h_1\left((\Ad k^{-1})(ln)\right)\,dn\,dk,\,\,\,\,\forall l\in L_Q(F),\]
    et $\gamma(Q)>0$ est la dérivée de Radon-Nikodym de la mesure sur $G(F)$ par rapport à la mesure produit sur $L_Q(F)N_Q(F)K=G(F)$.
\end{enumerate}    
\end{proposition}

On souligne la présence du facteur $\frac{|D^G(x)|}{|D^{\g}(x)|}$ dans la ligne \eqref{eq:TFlocaldescentgeom} à cause de notre convention pour une intégrale orbitale pondérée (\cite[numéro 3.3.2]{YDL23b}).

\subsubsection{Côté spectral}Le contenu de ce numéro est emprunté à \cite[sections 1, 3.25, 3.26]{Walds12}. 

La discussion ci-dessous prend en compte le fait que le R-groupe d'une représentation irréductible unitaire de $G(F)$ est trivial (remarque \ref{rem:R-grpistriv}).

Soit $L\in \L^G$, on note $a_{L,F}$ l'image de l'application de Harish-Chandra $H_L: L(F)\to a_L$, c'est un réseau de $a_L$. On note ensuite $a_{L,F}^\vee$ le groupe $\{\lambda\in a_L^\ast : \langle \lambda,H\rangle\in 2i\Z,\,\,\forall H\in a_{L,F}\}$. Posons enfin $a_{L,F}^\ast=a_L^\ast/a_{L,F}^\vee$ et $ia_{L,F}^\ast=ia_L^\ast/ia_{L,F}^\vee$.

Soient $L\in \L^G$, $M\in \L^L$ et $\sigma \in \Pi_2(M(F))$. Pour tout élément $s\in W_0^{L}$ on fixe un représentant $w_s\in L(F)\cap K$. 
Pour $l\in L(F)$ on note $(\Ad l)\sigma$ la représentation de $(\Ad l)M(F)$ agissant sur le même espace que celui sous-jacent de $\sigma$ donnée par $(\Ad l) \sigma( (\Ad l)m)=\sigma(m)$. On pose $W_0^L(\sigma)=\{s\in W_0^L: s\cdot M=M, w_s\cdot \sigma\simeq \sigma\}/W_0^{M}$.

Soit $L\in \L^G$, on note $E(L)=\{(M,
\sigma):M\in\L^L(M_0),\sigma\in \Pi_2(M(F))\}$. Cet ensemble est muni d'une action par $W_0^L\times ia_{L,F}^\ast$ sur chaque composante. Soit $(M,\sigma)\in E(L)$, on note $\text{Stab}(W_0^L\times ia_{L,F}^\ast,(M,\sigma))$ son stabilisateur dans $W_0^L\times ia_{L,F}^\ast$, c'est un groupe fini contenant le groupe $W_0^{M}$ comme sous-groupe distingué et son quotient $\text{Stab}(W_0^L\times ia_{L,F}^\ast,(M,\sigma))/W_0^{M}$ contient $W_0^L(\sigma)$ comme sous-groupe distingué. On pose ainsi 
\begin{align*}
\text{Stab}'(W_0^L\times ia_{L,F}^\ast,(M,\sigma))&=\left(\text{Stab}(W_0^L\times ia_{L,F}^\ast,(M,\sigma))/W_0^{M}\right)/W_0^L(\sigma).
\end{align*}

Soient $L\in \L^G$, $M\in \L^L$ et $\sigma \in \Pi_2(M(F))$. Un élément $s\in W_0^L(\sigma)$ agit naturellement sur $a_{M}$. Notons $a_{M}^s$ le sous-espace des points fixes par cette action. Il contient $a_L$. On note $W_{0,\reg}^L(\sigma)=\{s\in W_0^L(\sigma): a_{M}^s=a_L\}$. Soit $\alpha\in \Sigma(\g;A_M)$, on note $M_\alpha$ le sous-groupe de Levi semi-standard de $G$ avec $\Sigma(\m_\alpha;A_M)=\{\alpha,-\alpha\}$. L'ensemble des racines $\alpha$ telles que la mesure de Plancherel de $\sigma$ dans $M_\alpha(F)$ s'annule est un système de racines, dont $W_0^L(\sigma)$ est le groupe de Weyl (\cite[secion 1.11]{Walds12})
, on note $\epsilon_\sigma: W_0^L(\sigma)\to \{\pm 1\}$ le signe usuel. On note $E_\disc(L)=\{(M,\sigma)\in E(L) : W_{0,\reg}^L(\sigma)\not=\emptyset\}$. Il convient de noter que $(L,\sigma)\in E_\disc(L)$ pour tout $\sigma\in \Pi_2(L(F))$. L'ensemble $E_\disc(L)$ est muni d'une action par $W_0^L\times ia_{L,F}^\ast$ sur chaque composante. On note également $E_\disc(L)/\text{conj}$ l'ensemble des classes de $L(F)$-conjugaison de l'ensemble $E_\disc(L)$, il est muni d'une action par $ia_{L,F}^\ast$. On pose, lorsque $(M,\sigma)\in E_\disc(L)$,
\begin{align*}
\iota(M,\sigma)=|W_0^L(\sigma)|^{-1}\sum_{s\in W_{0,\reg}^L(\sigma)}  \epsilon_\sigma(s) |\det(s-1)_{a_{M}^L}|^{-1}. 
\end{align*}

Soient $L\in \L^G$, $M\in\L^L$, $\sigma\in \Pi_2(M(F))$, $\Lambda\in ia_L^\ast$, et $h_1,h_2\in \mathcal{HC}(G(F))$. On prend $R\in \P^G(M)$. Pour $\lambda \in ia_{L,F}^\ast$, on définit une $(G,L)$-famille $(J_Q^G(\I_{R\cap L}^L(\sigma_\lambda),h_1,h_2,\Lambda))_{Q\in \P^G(L)}$ par 
\begin{align*}
J_Q^G(&\I_{R\cap L}^L(\sigma_\lambda),h_1,h_2,\Lambda)\\
&=\overline{\tr\left(\mathcal{N}_{\overline{Q}}(P,\I_{R\cap L}^L(\sigma),\lambda,\Lambda)\I_{R}^G(\sigma_\lambda,h_1)\right)}\tr\left(\mathcal{N}_{Q}(P,\I_{R\cap L}^L(\sigma),\lambda,\Lambda)\I_{R}^G(\sigma_\lambda,h_2)\right),    
\end{align*}
où $\overline{Q}$ est l'opposé de $Q$, et $P\in \P^G(L)$. Cette $(G,L)$-famille est bien définie sur $ia_{L,F}^\ast$ et elle ne dépend ni du choix de $R$ ni du choix de $P$ (\cite[section 3.25]{Walds12}). On en déduit le poids $J_L^G(\I_{R\cap L}^L(\sigma_\lambda),h_1,h_2,\Lambda)$, et on pose 
\begin{equation*}
J_L^G(\I_{R\cap L}^L(\sigma_\lambda),h_1,h_2)=J_M^G(\I_{R\cap L}^L(\sigma_\lambda),h_1,h_2,\Lambda=0).
\end{equation*}
Puis
\begin{align*}    
J_{\spec,L}^G(h_1,h_2)=\sum_{(M,\sigma)\in(E_\disc(L)/\text{conj})/ia_{L,F}^\ast}|\text{Stab}'(W_0^{L}&\times ia_{L,F}^\ast,(M,\sigma))|^{-1}\iota(M,\sigma)\\
&\cdot\int_{ia_{L,F}^\ast}J_L^G(\I_{R\cap L}^L(\sigma_\lambda),h_1,h_2)\,d\lambda.
\end{align*}
Enfin 
\[J_{\spec}^G(h_1,h_2)=\sum_{L\in\L^G(M_0)}|W_0^L||W_0^G|^{-1}(-1)^{\dim a_L^G}J_{\spec,L}^G(h_1,h_2).\]

\begin{proposition}Soient $L\in\L^G$, $Q\in \F^G(L)$, et $h_1,h_2\in \mathcal{HC}(G(F))$.
\begin{enumerate}
    \item Soit $\pi\in \Pi_\temp(L(F))$. On a 
    \begin{equation}\label{eq:TFlocaldescentspec}
    J_L^G(\pi,h_1,h_2)=\sum_{(L_1,L_2)\in \L^G(L)^2} d_L^G(L_1,L_2)\overline{J_L^{Q_{L_1}}(\pi,h_1)}J_L^{Q_{L_2}}(\pi,h_2).   
    \end{equation}
    \item (Formule de descente de l'induction) Soient $M\in\L^L$ et $\pi\in \Pi_\temp(M(F))$. On a 
    \[J_L^Q(\I_P^L(\pi),h_1)=\sum_{M_1\in\L^{L_Q}(M)}d_M^{L_Q}(L,M_1)J_M^{Q_{M_1}}(\pi,h_1).\]
    \item (Formule de descente parabolique) Soit $\pi\in \Pi_\temp(L(F))$. On a,
    \[J_L^Q(\pi,h_1)=J_L^{L_Q}(\pi,h_{1,Q}).\]
\end{enumerate}    
\end{proposition}

\subsubsection{Formule finale}
\begin{theorem}[{{\cite[corollaire 5.3]{Art94}}}] \label{thm:localTF}
Pour tout $h_1\in \mathcal{HC}(G(F))$ et $h_2\in \mathcal{HC}(G(F))$, toutes les expressions intervenant dans la définition de $J_{\spec}^G(h_1,h_2)$ et $J_{\geom}^G(h_1,h_2)$ sont convergentes et ces expressions elles-mêmes le sont. On a l'égalité
\[J_{\spec}^G(h_1,h_2)=J_{\geom}^G(h_1,h_2).\]
Cette formule ne dépend d'aucun choix de mesure.
\end{theorem}

\subsection{Lemmes supplémentaires}
Soit $F$ un corps $p$-adique.

\begin{lemma}\label{lem:pluglocalintoglobal}
Soient $G$ un groupe du type GL défini sur $F$, $L$ un sous-groupe de Levi de $G$, et $\pi\in \Pi_\temp(L(F))$. Il existe $\dot{F}$ un corps totalement imaginaire, $\underline{S}$ un ensemble fini de places finies de $\dot{F}$, $\dot{G}$ un groupe du type GL défini sur $\dot{F}$, $\dot{L}$ un sous-groupe de Levi de $\dot{G}$, et $\dot{\pi}\in\Pi_{L^2}(\dot{L}(\A_{\dot{F}})^1)$ non-ramifié hors $\underline{S}$, tels que
\begin{enumerate}
    \item $|\underline{S}| \geq 2$ ;
    \item pour tout $v\in \underline{S}$, la donnée $(\dot{F}_v,\dot{G}_v,\dot{G}_v^\ast,\dot{L}_v,\dot{L}_v^\ast,\dot{\pi}_v)$ est isomorphe à $(F,G,G^\ast,L,L^\ast,\pi)$ ;
    \item pour tout $v\not\in \underline{S}$, $\dot{G}_v$ est quasi-déployé.
\end{enumerate}
\end{lemma}
\begin{proof}
Il résulte de la preuve de \cite[11.1 proposition]{Walds97} l'existence de $\dot{F}$ un corps totalement imaginaire, $\dot{G}$ un groupe du type GL défini sur $\dot{F}$, $\dot{L}$ un sous-groupe de Levi de $\dot{G}$, et $\underline{S}$ un ensemble fini de places finies de $\dot{F}$ tels que   
\begin{enumerate}
    \item $|\underline{S}| \geq 2$ ;
    \item pour tout $v\in \underline{S}$, $(\dot{F}_v,\dot{G}_v,\dot{G}_v^\ast,\dot{L}_v,\dot{L}_v^\ast)$ est isomorphe à $(F,G,G^\ast,L,L^\ast)$ ;
    \item pour tout $v\not\in \underline{S}$, $\dot{G}_v$ est quasi-déployé.
\end{enumerate}

Quant à la réalisation d'une représentation tempérée  en tant que composante locale d'une représentation automorphe, on suppose d'abord que $\pi\in \Pi_2(L(F))$ et $L$ est quasi-déployé. Alors il existe $h_\pi\in \H((L/Z(L))(F))$ tel que $\langle \pi, h_\pi\rangle=1$ et $\langle \pi', h_\pi\rangle=0$ pour toute représentation générique $\pi'\not\simeq \pi$. L'existence d'une telle fonction découle du théorème de Paley-Wiener (\cite[proposition 11.2]{Rog88}) et de la classification de Bernstein-Zelevinski (\cite[théorèmes 9.3 et 9.7]{Zel80}). Prenons $v_1\not\in \underline{S}$ une place finie de $\dot{F}$. Prenons $\dot{h}\in \H((\dot{L}/Z(\dot{L}))(\A_{\dot{F}}))$ une fonction test telle que 
\begin{enumerate}
    \item pour tout $v\in \underline{S}$, $\dot{h}_v=h_\pi$ ;
    \item pour $v=v_1$, $\dot{h}_{v_1}$ est supporté sur le lieu régulier elliptique de $\dot{L}(\dot{F}_{v_1})$ ;
    \item pour $v\not \in \underline{S}\cup\{v_1\}$ place finie, $\dot{h}_{v}$ est la fonction caractéristique de $\dot{L}(\O_{\dot{F}_v})$.
\end{enumerate}
Plongeons $\dot{h}$ dans la formule des traces simple de $\dot{L}/Z(\dot{L})$ (cf. \cite[A.I.d. théorème]{DKV84}, on fixe 1 comme caractère central). 
Quitte à modifier $\dot{h}_{v_1}$ et $\dot{h}_{\infty}$ on peut s'arranger pour que le côté géométrique admette un seul terme non nul. En conclusion, on vient de démontrer l'existence de $\dot{\pi}\in \Pi_{L^2}(\dot{L}(\A_{\dot{F}})^1)$ non-ramifiée hors $\underline{S}$ tel que $\dot{\pi}_v\simeq \rho$ pour tout $v\in \underline{S}$. 

Abandonnons à présent l'hypothèse que $L$ est quasi-déployé. On a $\bLJ|_{\Pi_2(L(F))}^{-1}(\pi)\in \Pi_2(L^\ast(F))$, il existe alors $\dot{\rho}^\ast\in \Pi_{L^2}(\dot{L}(\A_{\dot{F}})^1)$ non-ramifié hors $\underline{S}$ tel que $\dot{\rho}_v^\ast\simeq \bLJ^{-1}(\pi)$ pour tout $v\in \underline{S}$. Il est facile de vérifier que $\dot{\rho}^\ast$ est dans l'image de $\bG^{\dot{L}} : \Pi_{L^2}(\dot{L}(\A_{\dot{F}})^1)\to \Pi_{L^2}(\dot{L}^\ast(\A_{\dot{F}})^1)$ : sans perte de généralité on peut supposer que $\dot{L}=\GL_{m,\dot{D}}$ avec $\dot{D}$ une algèbre à division centrale sur $\dot{F}$, et pour toute place $v$ on note $d_v$ l'ordre de $\dot{D}\otimes_{\dot{F}}\dot{F}_v$ dans le groupe de Brauer $\text{Br}(\dot{F}_v)$. Alors $\dot{\rho}_v$ est bien $d_v$-compatible pour tout $v$, donc $\dot{\rho}^\ast\in \Pi_{L^2}^{\dot{D}}(\dot{L}^\ast(\A_{\dot{F}})^1)=\text{Im}(\bG^{\dot{L}})$. On note $\dot{\rho}$ l'image réciproque de $\dot{\rho}^\ast$, elle est non-ramifiée hors $\underline{S}$ et $\dot{\rho}_v\simeq \bLJ(\bLJ^{-1}(\pi))=\pi$ pour tout $v\in \underline{S}$. En conclusion on peut toujours réaliser une représentation de carré intégrable en tant que composante locale d'une représentation automorphe cuspidale. Le cas général $\pi\in \Pi_\temp(L(F))$ résulte de la décomposition de Langlands \eqref{eq:Langlandsdecomp} ainsi que des théorèmes \ref{thm:classtempbysquareint} et \ref{thm:irredunitparaind}.
\end{proof}

\begin{lemma}\label{lem:IOPcentrallimit}Soit $\mathfrak{z}$ le centre de $\g$. Pour tout $f\in \S(\g(F))$, il existe $U$ un voisinage de $0$ de $\mathfrak{z}(F)$ tel que $f(Z+X)=f(X)$ pour tout $Z\in U$ et $X\in \g(F)$.
\end{lemma}
\begin{proof}
Cela vient du fait que $\g(F)$ est localement compact et totalement discontinu. \qedhere
\end{proof}

\subsection{Identité locale}\label{subsec:idenspecloc}

Soit dans la suite $F$ un corps $p$-adique. Soit $\eta:G\to G^\ast$ une forme intérieure avec $G^\ast$ un groupe du type GL quasi-déployé sur $F$. Pour chaque tore maximal de $G$ ou de $G^\ast$ on munit son groupe des $F$-points d'une mesure de Haar, et si deux tores maximaux sont conjuguées par un élément de $G(F)$ ou de $G^\ast(F)$ alors on demande que les mesures de leur groupes des $F$-points se correspondent par conjugaison par ce même élément. On choisit aussi un produit scalaire $W_0^G$-invariant (resp. $W_{0^\ast}^{G^\ast}$-invariant) sur $a_{M_0}$ (resp. $a_{M_{0^\ast}}$). On prend sur tout sous-espace de $a_{M_0}$ (resp. $a_{M_{0^\ast}}$) la mesure engendrée par ce produit scalaire. On exige que si $T$ est un tore maximal de $G$ et $T^\ast$ un tore maximal de $G^\ast$, et $y\in G^\ast(\overline{F})$ est tel que $(\Ad y)\circ\eta|_T$ soit un isomorphisme défini sur $F$ de $T$ sur $T^\ast$, alors les mesures sur $T(F)$ et $T^\ast(F)$ se correspondent par cet isomorphisme. Pour tout $M\in\L^G$, le torseur intérieur $\eta$ induit un isomorphisme d'espaces  vectoriels réels $a_M\xrightarrow{\sim} a_{M^\ast}$, on exige également que les mesures sur $a_M$ et $a_{M^\ast}$ se correspondent par cet isomorphisme. 

\begin{definition}\label{def:loctemptrans}
Soit $L\in \L^{G}$. On dit qu'un élément $\pi'\in \Pi_\temp(L^\ast(F))$ se transfère à $L$ si $\pi\eqdef\bLJ(\pi')\in \Pi_\temp(L(F))$ (cf. corollaire \ref{coro:localJLtemp}). On écrit $\pi'=\pi^\ast$ et $\pi\arr \pi^\ast$. 
\end{definition}

Tout élément de $\Pi_2(L^\ast(F))$ se transfère (cf. théorème \ref{thm:appLJlocal} (4)). Nous abordons ensuite diverses notions de transfert d'une fonction test, en rappelant également la définition déjà posée (\cite[théorème 4.4]{YDL23b}) afin de mieux les comparer.

\begin{definition}[Transfert d'une fonction test]\label{def:transfergeomgeom'spec}~{}
\begin{enumerate}
    \item Pour tout $L\in \L^G$, $Q\in \F^G(L)$ et $\o$ une classe de $L(F)$-conjugaison dans $\mathfrak{l}(F)$ avec $(\Ad G)\o$ une classe de conjugaison semi-simple régulier, on appelle $(L,Q,\o)$ un triplet géométrique de $G$. 
    \item Soient $f\in \S(\g(F))$ et $f^\ast\in \S(\g^\ast(F))$. Si
    \begin{align*}
    J_{L'}^{Q'}(\mathfrak{o}',f^\ast)=\begin{cases}
    J_{L}^{Q}(\mathfrak{o},f) & \text{si $(L',Q',\mathfrak{o}')= (L^\ast,Q^\ast,\mathfrak{o}^\ast)$ ;} \\
    0 & \text{si $Q'$ ne se transfère pas.}
    \end{cases}
    \end{align*}
    pour tout triplet géométrique $(L',Q',\mathfrak{o}')$ de $G^\ast$, alors on note $f\underset{\geom}{\arr} f^\ast$.
    
    \item Pour tout $L\in \L^G$, $Q\in \F^G(L)$ et $\pi\in \Pi_\temp(L(F))$, on appelle $(L,Q,\pi)$ un triplet spectral de $G$.
    \item Soient $f\in \S(\g(F))$ et $f^\ast\in \S(\g^\ast(F))$, si 
    \begin{align*}
    J_{L'}^{Q'}(\pi',f_s^\ast)=\begin{cases}
    e(G) J_L^Q(\pi,f_s) & \text{si $(L',Q',\pi')=(L^\ast,Q^\ast,\pi^\ast)$ ;} \\
    0 &\text{si $Q'$ ne se transfère pas.}
    \end{cases}
    \end{align*}
    pour tout $s\in \D_\temp(G^\ast)$, et tout triplet spectral $(L',Q',\pi')$ de $G^\ast$, alors on note $f\underset{\spec}{\arr} f^\ast$.
\end{enumerate}
\end{definition}

\begin{theorem}\label{thm:transfertspecmeaning}
Pour tout torseur intérieur $\eta:G\to G^\ast$ avec $G^\ast$ un groupe du type GL quasi-déployé sur $F$, tout $f\in \S(\g(F))$ et tout $f^\ast \in \S(\g^\ast(F))$, on a l'implication « $f\underset{\geom}{\arr}f^\ast$ » $\Rightarrow$ « $f\underset{\spec}{\arr}f^\ast$ ».
\end{theorem}

\begin{proof}
Dans la preuve de ce théorème on va munir les groupes $N_P(F)$, $M_P(F)$ pour tout $P\in \P^G$ (resp. $N_{P'}(F)$, $M_{P'}(F)$ pour tout $P'\in \P^{G^\ast}$), et le groupe $K$ (resp. $K^\ast$), des mesures de Haar, de telle sorte que
\[\int_{G(F)} h(x)\,dx =\int_{M_P(F)}\int_{N_P(F)}\int_{K} h(mnk)\,dk\,dn\,dm,\,\,\,\,\forall P\in \P^G, h\in L^1(G(F)).\]
(resp. 
\[\int_{G^\ast(F)} h(x)\,dx =\int_{M_{P'}(F)}\int_{N_{P'}(F)}\int_{K^\ast} h(mnk)\,dk\,dn\,dm,\,\,\,\,\forall P'\in \P^{G^\ast}, h\in L^1(G^\ast(F)).)\]
On munit aussi $\n_P(F)$ (resp. $\n_{P'}(F)$) de la mesure de Haar qui correspond à celle de $N_P(F)$ (resp. $\n_{P'}(F)$) au sens expliqué dans le numéro \ref{YDLspecsubsec:fonctionmeropluvar}.  Ces choix supplémentaires des mesures n'altèrent en rien la validité des résultats recherchés, ces options sont imposées uniquement pour des considérations de commodité.

Constatons que si 
\[\S(\g(F))\ni f\underset{\geom}{\arr}f^\ast\in \S(\g^\ast(F))\,\,\,\, \text{(resp. $\underset{\spec}{\arr}$)}\,\,\,\,\text{ relativement à $\eta:\g\to\g^\ast$},\]
alors pour tout $\F^{G}\ni P\arr P^\ast\in \F^{G^\ast}$ on a 
\[\S(\m_P(F))\ni f_P\underset{\geom}{\arr}f_{P^\ast}^\ast\in \S(\m_{P^\ast}(F))\,\,\,\,\text{(resp. $\underset{\spec}{\arr}$)}\,\,\,\,\text{relativement à $\eta|_{\m_P}:\m_P\to\m_{P^\ast}$}.\]
Aussi, pour tout $P'$ qui ne se transfère pas,  
\[\S(\m_{P'}(F))\ni f_{P'}^\ast \underset{\geom}{\arr}0 \in \S(\m_{P'}(F))\,\,\,\,\text{(resp. $\underset{\spec}{\arr}$)}\,\,\,\,\text{relativement à $
\Id:\m_{P'}\to \m_{P'}$}.\]
Ces propriétés nous permettront de démontrer le théorème en utilisant une récurrence sur $\dim G$ et la formule de descente parabolique.

Faisons une récurrence sur $\dim G$. Si $\dim G=1$, alors $\eta :G\to G^\ast$ est un torseur trivial. On prend $x\in G^\ast(\overline{F})$ tel que $(\Ad x)\circ\eta|_G$ soit un isomorphisme défini sur $F$ de $G$ sur $G^\ast$. Il est clair que « $f\underset{\geom}{\arr}f^\ast$ » $\Leftrightarrow$ « $f\circ((\Ad x)\circ\eta|_G)^{-1}=f^\ast$ ». Puis « $f\underset{\spec}{\arr}f^\ast$ » $\Leftrightarrow$ « $f\circ((\Ad x)\circ\eta|_G)^{-1}=f^\ast$ » d'après la formule de Plancherel et le lemme \ref{lem:IOPcentrallimit}. La preuve se termine.

Supposons que la proposition est établie pour $\dim G<m\in\mathbb{N}_{>1}$. On se place maintenant au cas où $\dim G=m$. Soit $\S(\g(F))\ni f\underset{\geom}{\arr}f^\ast\in \S(\g^\ast(F))$. Soient $(L',Q',\pi')$ un triplet spectral et $s\in \D_\temp(G^\ast)$. On veut montrer que 
\begin{equation}\label{eq:geomtospecfinalwanteq}
J_{L'}^{Q'}(\pi',f_s^\ast)=\begin{cases}
    e(G)J_{L}^{Q}(\pi,f_s) & \text{si $(L',Q',\pi')= (L^\ast,Q^\ast,\pi^\ast)$ ;} \\
    0 & \text{si $Q'$ ne se transfère pas.}
    \end{cases}    
\end{equation}
Procédons ensuite à une récurrence sur $\dim a_{L'}^{L_{Q'}'}$. Si $\dim a_{L'}^{L_{Q'}'}=0$, alors on peut supposer que $L'=L_{Q'}'=G^\ast$, sinon on applique la formule de descente parabolique et on utilise l'hypothèse de récurrence sur $\dim G$. L'égalité voulue résulte du lemme \ref{YDLspeclem:tracetransferiden} et de l'holomorphie des termes en question sur $\D_\temp(G)$ (proposition \ref{prop:holotempcarainDt}).

Supposons que le théorème vaut au rang $\dim a_{L'}^{L_{Q'}'}<n\in\mathbb{N}_{\geq 1}$. On se place dans le cas de $\dim a_{L'}^{L_{Q'}'}=n$. Par la formule de descente de l'induction, la formule de descente parabolique, et l'hypothèse de récurrence sur $\dim G$, on peut et on va supposer que $Q'=G^\ast$, $L'=L^\ast$ et $\pi'=\pi^\ast\in\Pi_2(L^\ast(F))$. On note $(\dot{F},\dot{G},\dot{L},\dot{\pi},\underline{S})$ la donnée obtenue en appliquant le lemme \ref{lem:pluglocalintoglobal} à $(F,G,L,\pi)$. On exige que les mesures des groupes impliquées vérifient les normalisations données dans la sous-section \ref{subsec:cotespecTFnon-inv}. On exige aussi que pour tout $v\in \underline{S}$, la mesure sur le groupe des $\dot{F}_v$-points d'un tore maximal de $\dot{G}_v$ ou de $\dot{G}_v^\ast$ soit la même que la mesure sur le groupe des $F$-points du tore maximal correspondant de $G$ ou de $G^\ast$. On exige enfin que la mesure sur $\dot{G}_v(F_v)$ (resp. $\dot{G}_v^\ast(F_v)$) soit la même que la mesure sur $G(F)$ (resp. $G^\ast(F)$).

Supposons, pour l'instant, que $s\in\D_{\adel}(\dot{G})\mapsto J_{G^\ast}^{G^\ast}(\I_{Q_{L^\ast}}^{G^\ast}(\pi^\ast),f_s^\ast)$ n'est pas la fonction nulle, où $Q_{L^\ast}\in \P^{G^\ast}(L^\ast)$. Prenons $\dot{f}\in \S(\g(\A_{\dot{F}}))$ tel que 
\begin{enumerate}
    \item pour tout $v\in \underline{S}$, $\dot{f}_v=f$ ;
    \item pour tout $v\in \V_{\fin}\setminus \underline{S}$, $\dot{f}_v=1_{\g(\O_{\dot{F}_v})}$ ;
    \item pour tout $v\in \V_{\infty}$, $J_{\dot{G}_v}^{\dot{G}_v}(\I_{Q_{\dot{L}_v}}^{\dot{G}_v}(\dot{\pi}_v),\dot{f}_{v,s})\not=0$.
\end{enumerate}
Comme $J_{\dot{G}_v}^{\dot{G}_v}(\I_{Q_{\dot{L}_v}}^{\dot{G}_v}(\dot{\pi}_v),\dot{f}_{v,s})=\S
    1_{\g(\O_{\dot{F}_v}),s}(z_{{\I_{Q_{\dot{L}_v}}^{\dot{G}_v}(\dot{\pi}_v)}})=\det(1-q_v^{-s}z_{\I_{Q_{\dot{L}_v}}^{\dot{G}_v}(\dot{\pi}_v)})^{-1}\not=0$ pour tout $v\in \V_{\fin}\setminus \underline{S}$, on a $J_{\dot{G}_v}^{\dot{G}_v}(\I_{Q_{\dot{L}_v}}^{\dot{G}_v}(\dot{\pi}_v),\dot{f}_{v,s})\not=0$ pour tout $v\in \V_{\dot{F}}$. Appliquons ensuite la formule de scindage au théorème \ref{thm:idenendoglobal} (2) pour $(\dot{G}^\ast,\dot{L}^\ast,\bG^{\dot{L}}(\dot{\pi}))$ et $(\dot{G},\dot{L},\dot{\pi})$, nous obtenons
    \begin{align*}
    J_{\dot{L}^\ast}^{\dot{G}^\ast}(\bG^{\dot{L}}(\dot{\pi}),\dot{f}_s^\ast)=\sum_{(M_v^\ast)_v\in \prod_{v}\L^{\dot{G}_v^\ast}(\dot{L}_v^\ast)}d_{\dot{L}^\ast}^{\dot{G}^\ast}((M_v^\ast)_v)\prod_{v\in\V_F}J_{\dot{L}_v^\ast}^{Q_{M_v^\ast}}(\bG^{\dot{L}}(\dot{\pi})_v,\dot{f}_{s,v}^\ast)    
    \end{align*}
    et 
    \begin{align*}
    J_{\dot{L}}^{\dot{G}}(\dot{\pi},\dot{f}_s)=\sum_{(M_v)_v\in \prod_{v}\L^{\dot{G}_v}(\dot{L}_v)}d_{\dot{L}}^{\dot{G}}((M_v)_v)\prod_{v\in\V_F}J_{\dot{L}_v}^{Q_{M_v}}(\dot{\pi}_v,\dot{f}_{s,v}).    
    \end{align*}
    Les membres de gauche des deux égalités sont égaux. Grâce à l'hypothèse de récurrence sur $\dim a_{L^\ast}^{G^\ast}$ et la définition \ref{def:deftransfertgeomoffunctioninspecart}
    on en déduit (on fixe $v_0\in\underline{S}$ une place) 
    \begin{align*}
     |\underline{S}|J_{L^\ast}^{G^\ast}(\pi^\ast,f_{s}^\ast)\prod_{v\not=v_0}J_{\dot{L}_v^\ast}^{Q_{\dot{L}_v^\ast}}(\dot{\pi}_v^\ast,\dot{f}_{s,v}^\ast)=|\underline{S}|J_{L}^{G}(\pi,f_{s})\prod_{v\not=v_0}J_{\dot{L}_v}^{Q_{\dot{L}_v}}(\dot{\pi}_v,\dot{f}_{s,v}).   
    \end{align*}
    On obtient l'égalité voulue en divisant les deux côtés par $|\underline{S}|\prod_{v\in \V_{\dot{F}}}e(\dot{G}_v)\prod_{v\not=v_0}J_{\dot{L}_v}^{Q_{\dot{L}_v}}(\dot{\pi}_v,\dot{f}_{s,v})=|\underline{S}|\prod_{v\in \V_{\dot{F}}}e(\dot{G}_v)\prod_{v\not=v_0}J_{\dot{G}_v}^{\dot{G}_v}(\I_{Q_{\dot{L}_v}}^{\dot{G}_v}(\dot{\pi}_v),\dot{f}_{s,v})\not=0$. 

Supposons à présent que $s\in\D_{\adel}(\dot{G})\mapsto J_{G^\ast}^{G^\ast}(\I_{Q_{L^\ast}}^{G^\ast}(\pi^\ast),f_s^\ast)$ est la fonction nulle. Fixons $v_0\in S$ une place. Prenons $\dot{f}\in \S(\g(\A_{\dot{F}}))$ tel que 
\begin{enumerate}
    \item $\dot{f}_{v_0}=f$ ;
    \item pour tout $v\in \V_{\fin}\setminus \underline{S}$, $\dot{f}_v=1_{\g(\O_{\dot{F}_v})}$ ;
    \item pour tout $v\in \V_{\infty}\cup (\underline{S}\setminus\{v_0\})$, $J_{\dot{G}_v}^{\dot{G}_v}(\I_{Q_{\dot{L}_v}}^{\dot{G}_v}(\dot{\pi}_v),\dot{f}_{v,s})\not=0$.
\end{enumerate}
Alors l'hypothèse de récurrence sur $\dim a_{L^\ast}^{G^\ast}$ et les mêmes démarches nous conduisent toujours à $J_{L^\ast}^{G^\ast}(\pi^\ast,f_s^\ast)=e(G) J_L^G(\pi,f_s)$, l'égalité voulue. La preuve s'achève.
\end{proof}

\begin{lemma}\label{lem:compobjTFlocalspec}
Soient $\L^G\ni L\arr L^\ast\in \L^{G^{\ast}}$ et $\Pi_{\temp}(L(F))\ni \pi \arr \pi^\ast \in \Pi_\temp(L^\ast(F))$. Soient $M\in \L^L$ et $\sigma\in \Pi_2(M(F))$ tels que $\pi\simeq\I_P^L(\sigma)$ avec $P\in \P^L(M)$ (la classe de $W_0^G$-conjugaison de $(M,\sigma)$ est uniquement déterminée d'après le théorème \ref{thm:classtempbysquareint}). Alors $(M,\sigma)\in E_\disc(L)$ si et seulement si $(M^\ast,\sigma^\ast)\in E_\disc(L^\ast)$. Le cas échéant, on a
\[|\text{Stab}'(W_{0}^L\times ia_{L,F}^\ast,(M,\sigma))|= |\text{Stab}'(W_{0^\ast}^{L^\ast}\times ia_{L^\ast,F}^\ast,(M^\ast,\sigma^\ast))|,\text{ et }\iota(M,\sigma)=\iota(M^\ast,\sigma^\ast).\]
\end{lemma}
\begin{remark}
En général il n'y a pas de bijection naturelle entre $E_\disc(L)$ et $E_\disc(L^\ast)$. \`{A} titre d'exemple soit $L$ la forme intérieure non-triviale de $L^\ast=\GL_{2,F}$, alors $(\GL_{1,F}\times\GL_{1,F},1)\in E_\disc(L^\ast)$ ne se transfère pas.    
\end{remark}
\begin{proof}
Le torseur intérieur $\eta$ induit naturellement les isomorphismes de groupes suivant : $a_{L}\xrightarrow{\sim}a_{L^\ast}$, $a_{L}^\ast\xrightarrow{\sim}a_{L^\ast}^\ast$, $a_{L,F}\xrightarrow{\sim}a_{L^\ast,F}$, d'où $ia_{L,F}^\ast\xrightarrow{\sim}ia_{L^\ast,F}^\ast$. Aussi $W_M^L\xrightarrow{\sim} W_{M^\ast}^{L^\ast}$, et cet isomorphisme commute aux actions sur $a_M$ et $a_{M^\ast}$, autrement dit $s\cdot \lambda =\eta(s)\cdot \eta(\lambda)$ avec $s\in W_M^L$ et $\lambda\in a_M$. En regardant le caractère de Harish-Chandra on voit que $\eta$ induit les isomorphismes de groupes $W_0^L(\sigma)\xrightarrow{\sim}W_0^{L^\ast}(\sigma^\ast)$ et $W_{0,\reg}^L(\sigma)\xrightarrow{\sim}W_{0,\reg}^{L^\ast}(\sigma^\ast)$, on note $s\mapsto s^\ast$ ce dernier isomorphisme. En particulier $(M,\sigma)\in E_\disc(L)$ si et seulement si $(M^\ast,\sigma^\ast)\in E_\disc(L^\ast)$. Le cas échéant, les isomorphismes ci-dessus nous assurent que $\eta$ induit l'isomorphism de groupes $\text{Stab}'(W_{0}^L\times ia_{L,F}^\ast,(M,\sigma))\xrightarrow{\sim}\text{Stab}'(W_{0^\ast}^{L^\ast}\times ia_{L^\ast,F}^\ast(M^\ast,\sigma^\ast))$. Enfin il est clair que $\epsilon_\sigma(s)=\epsilon_{\sigma^\ast}(s^\ast)$ car $\eta$ induit un isomorphism de systèmes des racines impliqués, et $\det(s-1)_{a_M^L}=\det(s^\ast-1)_{a_{M^\ast}^{L^\ast}}$.
\end{proof}

\begin{theorem}\label{thm:transfergeom'specequiv} Pour tout torseur intérieur $\eta:G\to G^\ast$ avec $G^\ast$ un groupe du type GL quasi-déployé sur $F$, tout $f\in \S(\g(F))$ et tout $f^\ast \in \S(\g^\ast(F))$, on a l'implication « $f\underset{\spec}{\arr}f^\ast$ » $\Rightarrow$ « $f\underset{\geom}{\arr}f^\ast$ ».
\end{theorem}

\begin{proof}
La démonstration s'inspire de celle de \cite[5.5 théorème]{Walds12}. Dans la preuve de ce théorème on va munir les groupes $N_P(F)$, $M_P(F)$ pour tout $P\in \P^G$ (resp. $N_{P'}(F)$, $M_{P'}(F)$ pour tout $P'\in \P^{G^\ast}$), et le groupe $K$ (resp. $K^\ast$), des mesures de Haar, de telle sorte que
\[\int_{G(F)} h(x)\,dx =\int_{M_P(F)}\int_{N_P(F)}\int_{K} h(mnk)\,dk\,dn\,dm,\,\,\,\,\forall P\in \P^G, h\in L^1(G(F)).\]
(resp. 
\[\int_{G^\ast(F)} h(x)\,dx =\int_{M_{P'}(F)}\int_{N_{P'}(F)}\int_{K^\ast} h(mnk)\,dk\,dn\,dm,\,\,\,\,\forall P'\in \P^{G^\ast}, h\in L^1(G^\ast(F)).)\]
On munit aussi $\n_P(F)$ (resp. $\n_{P'}(F)$) de la mesure de Haar qui correspond à celle de $N_P(F)$ (resp. $\n_{P'}(F)$) au sens expliqué dans le numéro \ref{YDLspecsubsec:fonctionmeropluvar}. Ces choix supplémentaires des mesures n'altèrent en rien la validité des résultats recherchés, ces options sont imposées uniquement pour des considérations de commodité. 

Constatons que si 
\[\S(\g(F))\ni f\underset{\geom}{\arr}f^\ast\in \S(\g^\ast(F))\,\,\,\, \text{(resp. $\underset{\spec}{\arr}$)}\,\,\,\,\text{ relativement à $\eta:\g\to\g^\ast$},\]
alors pour tout $\F^{G}\ni P\arr P^\ast\in \F^{G^\ast}$ on a 
\[\S(\m_P(F))\ni f_P\underset{\geom}{\arr}f_{P^\ast}^\ast\in \S(\m_{P^\ast}(F))\,\,\,\,\text{(resp. $\underset{\spec}{\arr}$)}\,\,\,\,\text{relativement à $\eta|_{\m_P}:\m_P\to\m_{P^\ast}$}.\]
Aussi, pour tout $P'$ qui ne se transfère pas,  
\[\S(\m_{P'}(F))\ni f_{P'}^\ast \underset{\geom}{\arr}0 \in \S(\m_{P'}(F))\,\,\,\,\text{(resp. $\underset{\spec}{\arr}$)}\,\,\,\,\text{relativement à $\Id:\m_{P'}\to \m_{P'}$}.\]
Ces propriétés nous permettront de démontrer le théorème en utilisant une récurrence sur $\dim G$ et la formule de descente parabolique.


Faisons une récurrence sur $\dim G$. Si $\dim G=1$, alors $\eta :G\to G^\ast$ est un torseur trivial. On prend $x\in G^\ast(\overline{F})$ tel que $(\Ad x)\circ\eta|_G$ soit un isomorphisme défini sur $F$ de $G$ sur $G^\ast$. Il est clair que « $f\underset{\geom}{\arr}f^\ast$ » $\Leftrightarrow$ « $f\circ((\Ad x)\circ\eta|_G)^{-1}=f^\ast$ ». Puis « $f\underset{\spec}{\arr}f^\ast$ » $\Leftrightarrow$ « $f\circ((\Ad x)\circ\eta|_G)^{-1}=f^\ast$ » d'après la formule de Plancherel et le lemme \ref{lem:IOPcentrallimit}. La preuve se termine.

Supposons que le théorème vaut pour tout $G$ groupe du type GL avec $\dim G< m\in\N\setminus\{0\}$. Soit maintenant $G$ un groupe du type GL avec $\dim G=m$. Soit $\S(\g(F))\ni f \underset{\spec}{\arr}f^\ast\in\S(\g^\ast(F))$. Soit $(L',Q',\o')$ un triplet géométrique. On veut montrer que 
\begin{equation}\label{eq:spectogeomfinalwanteq}
J_{L'}^{Q'}(\mathfrak{o}',f^\ast)=\begin{cases}
    J_{L}^{Q}(\mathfrak{o},f) & \text{si $(L',Q',\mathfrak{o}')= (L^\ast,Q^\ast,\mathfrak{o}^\ast)$ ;} \\
    0 & \text{si $Q'$ ne se transfère pas.}
    \end{cases}    
\end{equation}

Procédons ensuite à une récurrence sur $\dim a_{L'}^{L_{Q'}'}$. Si $\dim a_{L'}^{L_{Q'}'}=0$, alors on peut supposer que $L'=L_{Q'}'=G^\ast$, sinon on applique la formule de descente parabolique et on utilise l'hypothèse de récurrence sur $\dim G$. De même on peut supposer que $\o'$ est semi-simple elliptique dans $\g^\ast(F)$ (\cite[proposition A.13]{YDL23b}), sinon on applique la formule de descente de l'induction et la formule de descente parabolique et on utilise de nouveau l'hypothèse de récurrence sur $\dim G$. Supposons donc maintenant que  $L'=L_{Q'}'=G^\ast$ et $\o'=\o^\ast$ est est semi-simple elliptique dans $\g^\ast(F)$ (on sait que $\o'$ se transfère à $G$ par \cite[lemme 4.6]{YDL23b}), et $\o$ est également semi-simple elliptique dans $\g(F)$.

Prenons $Z\in \mathfrak{z}(F)$ tel que $Z+\o\in G(F)$, on a aussi $Z^\ast+\o^\ast \in G^\ast(F)$. Soient $y\in Z+\o$ et $y^\ast\in Z^\ast+\o^\ast$. En manipulant un simple argument de la partition de l'unité, on dispose, pour tout $k\in\N$, de $\mathcal{U}_k$ un voisinage de $y$ stable par $G(F)$ et de $\mathcal{U}_k'$ un voisinage de $y^\ast$ stable par $G^\ast(F)$, ainsi que de  $C_c^\infty(G(F))\ni h_k\underset{\geom}{\arr}h_k^\ast\in C_c^\infty(G^\ast(F))$ (on a donc aussi $h_k\underset{\spec}{\arr}h_k^\ast$ au vu du théorème \ref{thm:transfertspecmeaning}), tels que
    \begin{enumerate}
        \item[(a)] $(\mathcal{U}_k)_{k\in\N}$ (resp. $(\mathcal{U}_k')_{k\in\N}$) est un système fondamental de voisinages de $y$ (resp. $y^\ast$) ;  
        \item[(b)] $\supp(h_k)\subseteq \mathcal{U}_k$, $\supp(h_k^\ast)\subseteq \mathcal{U}_k'$, $h_k(y)=1$, et $h_k^\ast(y^\ast)=1$ ; 
        \item[(c)] pour tout $M\in \L^G$, $\supp(h_k)\cap (\Ad G(F))M\not=\emptyset$ si et seulement si $M=G$, de même, pour tout $M'\in \L^{G^\ast}$, $\supp(h_k^\ast)\cap (\Ad G^\ast(F))M'\not=\emptyset$ si et seulement si $M'=G^\ast$.  
    \end{enumerate}
    Si $R\in \F^G(M_0)\setminus\{G\}$ alors la condition (c) nous dit que le terme constant $h_{k,R}=0$. De même si $R^\ast\in \F^{G^\ast}(M_{0^\ast})\setminus\{G^\ast\}$ alors $h_{k,R^\ast}^\ast=0$. Ainsi, pour tout $s\in\D_\temp(G)$, $H \in\L^G$ et $\pi\in \Pi_\temp(H(F))$, nous avons
    \begin{align*}
    J_H^G(\pi,f_s,h_k)=\overline{J_H^{Q_{H}}(\pi,f_s)}J_H^{G}(\pi,h_k)=\overline{J_{H^\ast}^{Q_{H^\ast}^\ast}(\pi^\ast,f_s^\ast)}J_{H^\ast}^{G^\ast}(\pi^\ast,h_k^\ast)=J_{H^\ast}^{G^\ast}(\pi^\ast,f_s^\ast,h_k^\ast)
    \end{align*} 
    avec $Q_H$ la deuxième composante de $s(G,H)$ de l'équation \eqref{eq:formuledescentGMsec} (avec $M$ remplacé par $H$). \`{A} l'aide du lemme \ref{lem:compobjTFlocalspec} on en déduit que $J_\spec^G(f_s,h_k)=J_\spec^{G^\ast}(f_s^\ast,h_k^\ast)$. Le théorème \ref{thm:localTF} dit alors $J_\geom^G(f_s,h_k)=J_\geom^{G^\ast}(f_s^\ast,h_k^\ast)$. La condition (c) nous dit encore que si $H\in \L^G\setminus\{G\}$ alors $J_{\geom,H}^G(f_s,h_k)=0$, donc $J_\geom^G(f_s,h_k)=J_{\geom,G}^G(f_s,h_k)$. De même $J_\geom^{G^\ast}(f_s^\ast,h_k^\ast)=J_{\geom,G^\ast}^{G^\ast}(f_s^\ast,h_k^\ast)$. De ce fait,
    \begin{align*}
    0=&\sum_{T\in T_{\text{ell}}(G)}|W^G(T)|^{-1}\vol(A_G(F)\backslash T(F))\int_{T(F)}\frac{|D^G(x)|}{|D^\g(x)|}\overline{J_G^G(x,f_s)}J_G^G(x,h_k)\,dx\\
    &-\sum_{T^\ast\in T_{\text{ell}}(G^\ast)}|W^{G^\ast}(T^\ast)|^{-1}\vol(A_{G^\ast}(F)\backslash T^\ast(F))\int_{T^\ast(F)}\frac{|D^{G^\ast}(x^\ast)|}{|D^{\g^\ast}(x^\ast)|}\overline{J_{G^\ast}^{G^\ast}(x^\ast,f_s^\ast)}J_{G^\ast}^{G^\ast}(x^\ast,h_k^\ast)\,dx^\ast.    
    \end{align*}
    En faisant $k\to \infty$ et en tenant compte du fait qu'il n'y a qu'un nombre fini de classes de conjugaison de tores maximaux, on voit que 
    \[J_{L}^Q(y,f)=J_{G}^G(y,f)=J_{G^\ast}^{G^\ast}(y^\ast,f^\ast)=J_{L^\ast}^{Q^\ast}(y^\ast,f^\ast).\]
    Grâce au lemme \ref{lem:IOPcentrallimit} on voit que $J_{L}^Q(\o,f)=\lim_{\substack{Z\in\mathfrak{z}(F)\to 0\\Z+\o\in G(F)}} J_{L}^Q(Z+\o,f)=J_{L^\ast}^{Q^\ast}(Z^\ast+\o^\ast,f^\ast)=J_{L^\ast}^{Q^\ast}(\o^\ast,f^\ast)$, ce qu'il fallait. La preuve du théorème dans le cas où $\dim G=m$ et $n=0$ se termine.

Supposons que le théorème vaut au rang $\dim a_{L'}^{L_{Q'}'}<n\in\mathbb{N}\setminus\{0\}$. On se place dans le cas de $\dim a_{L'}^{L_{Q'}'}=n$. Pour rappel on cherche à justifier l'équation \eqref{eq:spectogeomfinalwanteq}. Par la formule de descente de l'induction, la formule de descente parabolique, et l'hypothèse de récurrence sur $\dim G$, on peut et on va supposer que $Q'=G^\ast$ et $\o'$ semi-simple elliptique dans $\mathfrak{l}'(F)$. On n'a qu'à traiter le cas où $L=L^\ast$ et $\o'=\o^\ast$. 

Comme précédemment, quitte à translater les orbites par des éléments dans $\mathfrak{z}(F)$ et $\mathfrak{z}^\ast(F)$ qui se correspondent, on peut et on va supposer que $\o^\ast\subseteq L^\ast(F)$ et $\o\subseteq L(F)$. Soient $y\in \o$ et $y^\ast\in \o^\ast$. En manipulant un simple argument de la partition de l'unité, on dispose, pour tout $k\in\N$, de $\mathcal{U}_k$ un voisinage de $y$ stable par $G(F)$ et de $\mathcal{U}_k'$ un voisinage de $y^\ast$ stable par $G^\ast(F)$, ainsi que de  $C_c^\infty(G(F))\ni h_k\underset{\geom}{\arr}h_k^\ast\in C_c^\infty(G^\ast(F))$ (on a donc aussi $h_k\underset{\spec}{\arr}h_k^\ast$  au vu du théorème \ref{thm:transfertspecmeaning}), tels que
\begin{enumerate}
        \item[(a)] $(\mathcal{U}_k)_{k\in\N}$ (resp. $(\mathcal{U}_k')_{k\in\N}$) est un système fondamental de voisinages de $y$ (resp. $y^\ast$) ;  
        \item[(b)] $\supp(h_k)\subseteq \mathcal{U}_k$, $\supp(h_k^\ast)\subseteq \mathcal{U}_k'$, $h_k(y)=1$, et $h_k^\ast(y^\ast)=1$ ; 
        \item[(c)] pour tout $M\in \L^G$, $\supp(h_k)\cap (\Ad G(F))M\not=\emptyset$ si et seulement si $M$ est conjugué par $W_0^G$ à un élément de $\L^{G}(L)$ ; de même, pour tout $M'\in \L^{G^\ast}$, $\supp(h_k^\ast)\cap (\Ad G^\ast(F))M'\not=\emptyset$ si et seulement si $M^\ast$ est conjugué par $W_0^G$ à un élément de $\L^{G^\ast}(L^\ast)$ (cf. \cite[lemme 4.9]{YDL23b}).  
    \end{enumerate}

On tâche maintenant de prouver 
\[J_\spec^G(f_s,h_k)=J_\spec^{G^\ast}(f_s^\ast,h_k^\ast).\]
On décompose la démonstration en trois parties. On prouve en premier lieu que pour tout $M \in\L^G$, $(M_1,M_2)\in \L^G(M)^2$, et $\pi\in \Pi_\temp(M(F))$, il y a l'égalité 
    \begin{equation}\label{eq:TFlocalcompeq1}
        d_M^G(M_1,M_2)\overline{J_M^{Q_{M_1}}(\pi,f_s)}J_M^{Q_{M_2}}(\pi,h_k)=d_{M^\ast}^{G^\ast}(M_1^\ast,M_2^\ast)\overline{J_{M^\ast}^{Q_{M_1^\ast}}(\pi^\ast,f_s^\ast)}J_{M^\ast}^{Q_{M_2^\ast}}(\pi^\ast,h_k^\ast).  
    \end{equation}
En effet, cette égalité est triviale au vu de $d_M^G(M_1,M_2)=d_{M^\ast}^{G^\ast}(M_1^\ast,M_2^\ast)$, $f\underset{\spec}{\arr} f^\ast$ et $h_k\underset{\spec}{\arr} h_k^\ast$. 
    
En second lieu, on démontre que pour tout $M^\ast \in\L^{G^\ast}$, $(M_1^\ast,M_2^\ast)\in \L^G(M^\ast)^2$, et $\pi'\in \Pi_\temp(M^\ast(F))$ qui ne se transfère pas à $M$, il y a l'égalité 
    \begin{equation}\label{eq:TFlocalcompeq1'}
    d_{M^\ast}^{G^\ast}(M_1^\ast,M_2^\ast)\overline{J_{M^\ast}^{Q_{M_1^\ast}}(\pi',f_s^\ast)}J_{M^\ast}^{Q_{M_2^\ast}}(\pi',h_k^\ast)=0.  
    \end{equation}
Soient $A'\in \L^{M^\ast}$ et $\sigma'\in \Pi_2(A'(F))$ tels que $\pi'=\I_{R'}^{M^\ast}(\sigma')$, où $R'\in \P^{M^\ast}(A')$. Puisque $\pi'$  ne se transfère pas, on sait que $A'$ ne se transfère pas non plus. La formule de descente de l'induction implique
\[J_{M^\ast}^{Q_{M_1^\ast}}(\I_{R'}^{M^\ast}(\sigma'),f_s^\ast)=\sum_{B'\in\L^{M_1^\ast}(A')}d_{A'}^{M_1^\ast}(M^\ast,B')J_{A'}^{Q_{B'}}(\sigma',f_s^\ast).\]
D'un côté si $B'=B^\ast\in\L^{M_1^\ast}(A')$ se transfère à $M_1^\ast$, alors $d_{A'}^{M_1^\ast}(M^\ast,B^\ast)=0$ comme  $A'$ ne se transfère pas. De l'autre côté si $B'\in\L^{M_1^\ast}(A')$ ne se transfère pas à $M_1^\ast$, alors $J_{A'}^{Q_{B'}}(\sigma',f_s^\ast)=0$ par l'hypothèse $f\underset{\spec}{\arr} f^\ast$. Ainsi, l'égalité \eqref{eq:TFlocalcompeq1'} est confirmée. 

    En troisième lieu, on démontre que pour tout 
    $M' \in\L^{G^\ast}$ qui ne se transfère pas, $(M_1',M_2')\in \L^G(M')^2$, et $\pi'\in \Pi_\temp(M'(F))$, il y a l'égalité 
    \begin{equation}\label{eq:TFlocalcompeq1''}
    d_{M'}^{G^\ast}(M_1',M_2')\overline{J_{M'}^{Q_{M_1'}}(\pi',f_s^\ast)}J_{M'}^{Q_{M_2'}}(\pi',h_k^\ast)=0.  
    \end{equation}
    Si $d_{M'}^{G^\ast}(M_1',M_2')=0$ alors l'égalité est bien sûr vérifiée. Ensuite si $d_{M'}^{G^\ast}(M_1',M_2')\not=0$, alors au moins l'un des $M_1',M_2'$ ne se transfère pas à $G$. On a ou bien $J_{M'}^{Q_{M_1'}}(\pi',f_s^\ast)=0$, ou bien $J_{M'}^{Q_{M_2'}}(\pi',h_k^\ast)$, l'égalité est aussi vérifiée. Ainsi, l'égalité \eqref{eq:TFlocalcompeq1'} est confirmée.
    
    En combinant \eqref{eq:TFlocalcompeq1}, \eqref{eq:TFlocalcompeq1'}, \eqref{eq:TFlocalcompeq1''} et le lemme \ref{lem:compobjTFlocalspec}, on établit $J_\spec^G(f_s,h_k)=J_\spec^{G^\ast}(f_s^\ast,h_k^\ast)$. La formule des traces locale (théorème \ref{thm:localTF}) donne de ce fait \begin{equation}\label{eq:TFlocalcompeq2}
    J_\geom^G(f_s,h_k)=J_\geom^{G^\ast}(f_s^\ast,h_k^\ast).
    \end{equation}
    Comparons à présent les termes géométriques. Déjà si $M'\in\L^{G^\ast}$ ne se transfère pas alors
    \begin{equation}\label{eq:TFlocalcompeq3}
    J_{\geom,M'}^{G^\ast}(f_s^\ast,h_k^\ast)=0
    \end{equation} 
    en raison de l'équation \eqref{eq:TFlocaldescentgeom} et de la propriété $h_{k}\underset{\geom}{\arr}h_{k}^\ast$. 
    On veut ensuite montrer que 
    \begin{equation}\label{eq:TFlocalcompeq4}
    J_{\geom,M}^{G}(f_s,h_k)=J_{\geom,M^\ast}^{G^\ast}(f_s^\ast,h_k^\ast),\,\,\,\,{\text{si $M$ n'est pas conjugué à $L$ sous $W_0^G$}}.
    \end{equation} 
    Soit donc $M\in \L^G$ qui n'est pas conjugué à $L$ sous $W_0^G$. D'une part si $x'\in M^\ast(F)\cap G_{\rss}^\ast(F)$ ne se transfère pas, alors on peut écrire $x'=\Ind_{H'}^{M^\ast}(x_e')$ avec $H'\in\L^{M^\ast}$, $Q_{H'}\in\L^{G^\ast}(M_1')$, et $x_e'$ semi-simple elliptique dans $\mathfrak{h}'(F)$. On sait que $H'$ ne se transfère pas (\cite[proposition A.14]{YDL23b}). Au moyen de la formule de descente de l'induction, on obtient, pour tout $M_2^\ast\in\L^{G^\ast}(M^\ast)$, l'égalité $J_{M^\ast}^{Q_{M_2^\ast}}(x',h_{k}^\ast)=\sum_{L'\in \L^{M_2^\ast}(H')}d_{H'}^{M_2^\ast}(M^\ast,L')J_{H'}^{Q_L'}(x_e',h_k^\ast)$. Mais si $d_{H'}^{M_2^\ast}(M^\ast,L')\not =0$ alors $L'$ ne se transfère pas (\cite[corollaire 4.13]{YDL23b}). De la propriété $h_{k}\underset{\geom}{\arr}h_{k}^\ast$ on tire que $J_{M^\ast}^{Q_{M_2^\ast}}(x',h_{k}^\ast)=0$. Ainsi $J_{M^\ast}^{G^\ast}(x',f_s^\ast,h_{k}^\ast)=0$ en raison de l'équation \eqref{eq:TFlocaldescentgeom}. En définitive, si $x'\in M^\ast(F)\cap G_{\rss}^\ast(F)$ ne se transfère pas, alors $J_{M^\ast}^{G^\ast}(x',f_s^\ast,h_k^\ast)=0$.
    D'autre part si $x'=x^\ast$ se transfère, on va prouver que pour tout $(M_1,M_2)\in \L^G(M)^2$ il y a l'égalité 
    \begin{equation}\label{eq:TFlocalcompeq5}
        d_M^G(M_1,M_2)\overline{J_M^{Q_{M_1}}(x,f_s)}J_M^{Q_{M_2}}(x,h_k)=d_{M^\ast}^{G^\ast}(M_1^\ast,M_2^\ast)\overline{J_{M^\ast}^{Q_{M_1^\ast}}(x^\ast,f_s^\ast)}J_{M^\ast}^{Q_{M_2^\ast}}(x^\ast,h_k^\ast).  
    \end{equation}
    La preuve est presque la même que celle de \eqref{eq:TFlocalcompeq1}. Si $\dim a_{M}^{M_1}< n$ alors on a bien l'égalité \eqref{eq:TFlocalcompeq5} par l'hypothèse de récurrence sur $n$ et la propriété $h_{k}\underset{\geom}{\arr}h_{k}^\ast$. Si $d_M^G(M_1,M_2)=0$ alors $d_{M^\ast}^{G^\ast}(M_1^\ast,M_2^\ast)=0$, on a aussi l'égalité voulue. Si $\dim a_{M}^{M_1}> n$ et $d_M^G(M_1,M_2)\not=0$, on a  $h_{k,Q_{M_2}}=h_{k,Q_{M_2^\ast}}^\ast=0$.    Enfin si $\dim a_M^{M_1}=n$ et $d_M^G(M_1,M_2)\not=0$, on peut et on va supposer que $Q_{M_1}=G$ car sinon l'égalité \eqref{eq:TFlocalcompeq5} résulte de la formule de descente parabolique et l'hypothèse de récurrence sur $\dim G$. On a $a_M^G=a_{G}^G\oplus a_{M_2}^G$, ainsi $M=M_2$. Mais par l'hypothèse de récurrence sur $n$ imposée au début on a aussi $\dim a_{L'}^{Q_{L'}'}=\dim a_L^G=n$. Puisque $\dim a_M^G=n=a_L^G$ et $M$ n'est pas conjugué à $L$ sous $W_0^G$, on sait que $M$ n'est conjugué par $W_0^G$ à aucun élément de $\L^G(L)$. La condition (c) sur les fonctions $h_k$ et $h_k^\ast$ nous garantit $h_{k,Q_{M_2}}=h_{k,Q_{M_2^\ast}}^\ast=0$. En résumé l'égalité \eqref{eq:TFlocalcompeq5} a toujours lieu. 

    En combinant les équations \eqref{eq:TFlocalcompeq3}, \eqref{eq:TFlocalcompeq4}, \eqref{eq:TFlocalcompeq5}, on voit que 
    \begin{equation}\label{eq:TFlocalcompeq6}
    J_{\geom,L}^{G}(f_s,h_k)=J_{\geom,L^\ast}^{G^\ast}(f_s^\ast,h_k^\ast).
    \end{equation} 
    C'est-à-dire
    \begin{align*}
    0=&\sum_{T\in T_{\text{ell}}(L)}|W^L(T)|^{-1}\vol(A_G(F)\backslash T(F))\int_{T(F)}J_L^G(x,f_s,h_k)\,dx\\
    &-\sum_{T^\ast\in T_{\text{ell}}(L^\ast)}|W^{L^\ast}(T^\ast)|^{-1}\vol(A_{G^\ast}(F)\backslash T^\ast(F))\int_{T^\ast(F)}J_{L^\ast}^{G^\ast}(x^\ast,f_s^\ast,h_k^\ast)\,dx^\ast.    
    \end{align*}
    Par le biais de l'équation \eqref{eq:TFlocaldescentgeom}, de la formule de descente parabolique, et de l'hypothèse de récurrence sur $n$, on retient
    \begin{align*}
    0=&\sum_{T\in T_{\text{ell}}(G)}|W^G(T)|^{-1}\vol(A_G(F)\backslash T(F))\int_{T(F)}\frac{|D^G(x)|}{|D^\g(x)|}\overline{J_L^G(x,f_s)}J_L^{Q_L}(x,h_k)\,dx\\
    &-\sum_{T^\ast\in T_{\text{ell}}(G^\ast)}|W^{G^\ast}(T^\ast)|^{-1}\vol(A_{G^\ast}(F)\backslash T^\ast(F))\int_{T^\ast(F)}\frac{|D^{G^\ast}(x^\ast)|}{|D^{\g^\ast}(x^\ast)|}\overline{J_{L^\ast}^{G^\ast}(x^\ast,f_s^\ast)}J_{L^\ast}^{Q_{L^\ast}}(x^\ast,h_k^\ast)\,dx^\ast.    
    \end{align*}
    En faisant $k\to \infty$ et en tenant compte du fait qu'il n'y a qu'un nombre fini de classes de conjugaison de tores maximaux, on voit que $J_{L}^G(y,f)=J_{L^\ast}^{Q^\ast}(y^\ast,f^\ast),$ ce qu'il fallait. La preuve du théorème \ref{thm:transfergeom'specequiv} est d'ores et déjà terminée.
\end{proof}

En combinant les théorèmes \ref{thm:transfertspecmeaning} et \ref{thm:transfergeom'specequiv}, nous voyons que si $f\in \S(\g(F))$ et $f^\ast\in\S(\g^\ast(F))$, alors $f\underset{\geom}{\arr} f^\ast$ si et seulement si $f\underset{\spec}{\arr} f^\ast$. 



\begin{corollary}\label{coro:kazhdan0noninv} Soit $G$ un groupe du type GL sur $F$. Soit $f\in \S(\g(F))$. On a l'équivalence entre
\begin{itemize}
    \item $J_L^Q(\o,f)=0$ pour tout triplet géométrique $(L,Q,\o)$.
    \item $J_L^Q(\pi,f_s)=0$ pour tout triplet spectral $(L,Q,\pi)$ et tout $s\in \D_\temp(G)$.
\end{itemize}
\end{corollary}

On peut se restreindre à travailler uniquement sur le groupe.

\begin{corollary}
\label{coro:kazhdan0noninvgrp}~{}
Soit $G$ un groupe du type GL sur $F$. Soit $f\in C_c^\infty(G(F))$. On a l'équivalence entre
    \begin{itemize}
    \item $J_L^Q(\o,f)=0$ pour tout triplet géométrique $(L,Q,\o)$ avec $\o\subseteq L(F)$.
    \item $J_L^Q(\o,f)=0$ pour tout $L\in \L^G$, $Q\in \F^G(L)$ et $\o$ classe de $L(F)$-conjugaison dans $\mathfrak{l}(F)$.
    \item $J_L^Q(\pi,f)=0$ pour tout triplet spectral $(L,Q,\pi)$.
    \item $J_L^Q(\pi,f)=0$ pour tout $L\in \L^G$, $Q\in \F^G(L)$ et $\pi\in \Pi_{\unit}(L(F))$.
    \end{itemize}
\end{corollary}
    
\begin{proof}
L'équivalence entre le premier point et le deuxième point est \cite[corollaire 4.21]{YDL23b} ; celle entre le premier point et le troisième point est le corollaire \ref{coro:kazhdan0noninv} ; celle entre le troisième point et le quatrième point est une conséquence directe de \cite[proposition 6.1]{Art89I}.
\end{proof}

Le corollaire précédent constitue une généralisation non-invariante de \cite[théorème 0, (b) $\Leftrightarrow$ (c) $\Leftrightarrow$ (d)]{Kaz86}. Une généralisations non-invariante de l'énoncé (a) dans l'article de Kazhdan devrait être une version plus précise d'un lemme de Labesse (\cite[I.7.1.lemme]{Lab95}).

\newcommand{\noop}[1]{}

\end{document}